\documentclass[11pt, reqno]{amsart}
\usepackage{amsfonts,latexsym,amsmath,amsthm, amscd,geometry, amssymb}
\usepackage{amssymb}
\usepackage{latexsym}

\newcommand \nc{\newcommand}
\newtheorem{theorem}{Theorem}[section]
\newtheorem{lemma}[theorem]{Lemma}
\newtheorem{proposition}[theorem]{Proposition}
\newtheorem{corollary}[theorem]{Corollary}

\nc{\ba}{\begin{array}}\nc{\ea}{\end{array}}
\nc{\be}{\begin{eqnarray}}\nc{\ee}{\end{eqnarray}}
\nc{\beq}{\begin{equation}}\nc{\eeq}{\end{equation}}
\nc{\bex}{\begin{eqnarray*}}\nc{\eex}{\end{eqnarray*}}
\nc{\btm}{\begin{theorem}} \nc{\etm}{\end{theorem}}
\nc{\blm}{\begin{lemma}} \nc{\elm}{\end{lemma}}
\nc{\R}{\mathbb{R}}  \nc{\ld}{\lambda}
\nc{\va}{\varphi}
\nc{\ve}{\varepsilon}

\def\pf{\noindent{\bf Proof.\quad}}

\begin{document}
\title{On static and hydrodynamic biaxial nematic liquid crystals}
\author{Junyu Lin, Yimei Li, Changyou Wang}
\date{}
%\maketitle
% Address where research was carried out
\address{Department of Mathematics\\ South China  Univ. of Technology\\ Guangzhou,  Guangdong 510640, China}
\address{School of Mathematical Sciences, Beijing Normal University, Beijing 100875, China}
\address{Department of Mathematics\\ Purdue University\\ West Lafayette, IN 47907, USA}
\email{scjylin@scut.edu.cn,  lym@mail.bnu.edu.cn,  wang2482@purdue.edu} % Email
%\date{}
%\keywords{Biaxial harmonic map, Global existence, Partial regularity} %keywords optional
%\subjclass[2000]{}
%for subject classification codes see http://www.ams.org/mathscinet/msc/msc2010.html

% Abstract should be 150 words for a short paper and 300 words for a long paper.
\begin{abstract} In the first part of this paper, we will consider minimizing configurations of the Oseen-Frank energy functional $E(n, m)$ for a biaxial nematics $(n, m):\Omega\to \mathbb S^2\times \mathbb S^2$ with $n\cdot m=0$ in dimension three, and establish that it is smooth off a closed set of  $1$-dimension Hausdorff measure zero.
In the second part, we will consider a simplified Ericksen-Leslie system for biaxial nematics $(n, m)$ in a two dimensional domain and establish the existence of a unique global weak solution $(u, n, m)$ that is smooth off
at most finitely many singular times for any initial and boundary data of finite energy. They extend to biaxial nematics of earlier results
corresponding to minimizing uniaxial nematics by Hardt-Kindelerherer-Lin \cite{HKL} and 
a simplified hydrodynamics of uniaxial liquid crystal by Lin-Lin-Wang \cite{LLW10} respectively.
\end{abstract}
\maketitle %arranges what you've just done in a nice way

\section{Introduction on static and hydrodynamic biaxial nematics}

\subsection{Brief overview of static uniaxial nematics}
The problem of liquid crystals has been a fascinating subject ever since its inception. The macroscopic
theory of nematics describes these anisotropic materials in terms of a continuum field, and there are
a number of competitive and complementary mathematical models available in the literature. Perhaps
the earliest one was the Oseen-Frank model where the nematic liquid crystal is considered
at rest at constant temperature in the absence of applied electromagnetic fields, which is 
occupied in a container $\Omega\subset\mathbb R^3$ and  is described 
by a unit-vector field $n:\Omega\to\mathbb S^2$, the unit sphere in $\mathbb R^3$. 
In the Oseen-Frank model the total elastic energy is given by
\begin{equation}\label{uniaxial}
E(n)=\int_\Omega W(n,\nabla n)\,dx, 
\end{equation}
where
\begin{equation}
2W(n, \nabla n)=K_1({\rm{div}}\ n)^2+K_2(n\cdot{\rm{curl}}\ n)^2+K_2|n\times{\rm{curl}}\ n|^2
+(K_2+K_4)( {\rm{tr}}(\nabla n)^2-({\rm{div}}\ n)^2).
\end{equation} 
Here $K_1, \cdots, K_4$ are Frank constants, which are generally assumed to satisfy (\cite{Ericksen1})
$$K_1>0,\ K_2>0, \ K_3>0, \ K_2\ge |K_4|, \ {\rm{and}}\ \ 2K_1\ge K_2+K_4.$$

The principle mathematical question concerns the existence and partial regularity of a 
unit vector field $n\in H^1(\Omega,\mathbb S^2)=\big\{v\in H^1(\Omega, \mathbb R^3):
\ v(x)\in \mathbb S^2 \ {\rm{a.e.}}\ x\in\Omega\big\}$, that minimizes the Oseen-Frank energy
$E(n)$ among all $n\in H^1(\Omega,\mathbb S^2)$ having prescribed boundary values
$n_0$ on $\partial\Omega$. The Euler-Lagrange equation for an Oseen-Frank energy
minimizer $n\in H^1(\Omega,\mathbb S^2)$ is
\begin{equation}\label{uniaxial-EL}
{\rm{div}}\frac{\partial W}{\partial\nabla n}-\frac{\partial W}{\partial n}=\lambda(x) n,
\end{equation}
where $\lambda(x)$ is the Lagrange multiplier. 

In a pioneering paper \cite{HKL},  Hardt-Kinderlehrer-Lin have
established the existence of minimizer $n$ of $E$, which is smooth away from a closed
set $\Sigma\subset\Omega$ with Hausdorff dimension less than $1$.  Note that
when the Frank constants $K_1=K_2=K_3=K_2+K_4=1$, the Oseen-Frank energy
$E(n)=\frac12\int_\Omega|\nabla n|^2\,dx$ reduces to the Dirichlet energy so that
the minimizer becomes a minimizing harmonic map into $\mathbb S^2$, which is smooth
away from a discrete set of singularity (see Schoen-Uhlenbeck \cite{SU1}). We would like
to point out that it has been an open problem
whether the singular set of an energy minimizer of general Oseen-Frank energy is a discrete set. 

\subsection{Brief overview of hydrodynamics of uniaxial nematics}

When the nematic liquid crystal is at  motion, the relevant dynamic equation is the Ericksen-Leslie 
system proposed by Ericksen and Leslie back in 1960's (see \cite{Ericksen} \cite{Leslie}). For an 
incompressible underlying fluid, Lin \cite{Lin} first proposed the following simplified equation for the fluid velocity field $v$ and the uniaxial nematic liquid crystal director field $n$:
\begin{equation}\label{uniaxial_LCF1}
\begin{cases}\displaystyle
v_t+v\cdot\nabla v-\nu\Delta v+\nabla P=-{\rm{div}}((\nabla n)^T\frac{\partial W}{\partial\nabla n}) ,\\
{\rm{div}}\ v=0,\\ 
\displaystyle(I-n\otimes n)\big(n_t+v\cdot \nabla n+\frac{\partial W}{\partial n}-{\rm{div}}\frac{\partial W}{\partial\nabla n}\big)=0,
\end{cases}
\end{equation}
where $P$ stands for the hydrostatic pressure function and $\nu>0$ denotes the fluid viscosity constant. When we take the one constant approximation for
the Oseen-Frank energy density function $W(n,\nabla n)$, \eqref{uniaxial_LCF1} reduces to 
\begin{equation}\label{uniaxial_LCF2}
\begin{cases}
v_t+v\cdot\nabla v-\nu\Delta v+\nabla P=-{\rm{div}}(\nabla n\odot\nabla n) ,\\
{\rm{div}}\ v=0,\\ 
n_t+v\cdot \nabla n=\Delta n+|\nabla n|^2 n,
\end{cases}
\end{equation}
where $\nabla n\odot\nabla n=(\langle \frac{\partial n}{\partial x_i}, \frac{\partial n}{\partial x_j}\rangle)_{1\le i, j\le 3}$
denotes the Ericksen stress tensor. The system \eqref{uniaxial_LCF1} or \eqref{uniaxial_LCF2} is a macroscopic
continuum description of the time evolution of the materials under the influence of both the flow field $v$  and
the macroscopic average of the microscopic orientation configuration field $n$ of rod-like liquid crystals, which
reduces to the Oseen-Frank model in the static case. Recently, there have been many mathematical studies 
on the simplified Ericksen-Leslie system \eqref{uniaxial_LCF2}. In dimension two, Lin-Lin-Wang \cite{LLW10} and
Lin-Wang \cite{LW2} have established the existence of a unique global weak solution to the initial and boundary value problem of \eqref{uniaxial_LCF2} that has at most finitely many singular times, see also Hong \cite{Hong} and Xu-Zhang \cite{XZ} for related works.  In dimension three, Lin-Wang \cite{LW3} have made some partial progress on the
existence of global weak solutions of \eqref{uniaxial_LCF2}. The existence of finite time singularity of short time smooth solutions to \eqref{uniaxial_LCF2} has also been constructed by Huang-Lin-Liu-Wang \cite{HLLW}
and Lai-Lin-Wang-Wei-Zhou \cite{LLWWZ} in dimensions three and two respectively. We would also like to mention
the initial works by Lin-Liu \cite{LL1, LL2} on the simplified Ericksen-Leslie system on nematic liquid crystals
with variable degrees of orientation. The interested readers can also consult with the review articles by Ball
\cite{Ball} and Lin-Wang \cite{LW4}. 
\subsection{Introduction and main theorem on static biaxial nematics} 
The so-called $Q$-tensor model to describe the nematic liquid crystals in a uniform fashion
was later proposed by De Gennes in 1970's, in which a $Q$-tensor is symmetric $3\times 3$ matrix
of trace zero that describes both the amount of order and the orientation of the materials. The advantage of
$Q$-tensor model vs the Oseen-Frank vectorial field model is that it can describe more structures of the defect
and can also resolve the issue of orientability of the director field.  By using the eigenvalue decomposition, the
$Q$-tensor order parameter of a biaxial nematic can be written as
$$
Q=S(n\otimes n-\frac13 \mathbb I_3) + T(m\otimes m-\frac13 \mathbb I_3), \ \ n, m\in\mathbb S^2, \ n\perp m.
$$
The static Landau-De Gennes theory seeks minimizers of the the Landau-De Gennes free energy functional
\begin{equation}\label{LG}
I(Q)=\int_\Omega \big(\psi_{LdG}(Q,\nabla Q)+\psi_{b}(Q)\big)\,dx, \ Q\in H^2(\Omega, \mathcal{S}_0),
\end{equation}
subject to the strong anchoring condition $Q\big|_{\partial\Omega}=Q_0$, where $\mathcal{S}_0$ denotes
the space of all $Q$-tensors. It is well-known that when restricted to the subspace of rank-one $Q$-tensors,
the Landau-De Gennes theory reduces to the Ericksen model with variable degrees of orientations (see \cite{Ericksen1}), which further reduces to the Oseen-Frank model when the degree of orientations is
constant (or  $S=T$ is constant). 

When $S, T$ are distinct constants, Govers and Vertogen \cite{GV1, GV2} proposed around 1984
the elastic continuum theory of biaxial nematics based on the Landau-De Gennes $Q$-tensor theory.
The Govers-Vertogen model uses a pair of orthogonal unit vector fields $(n, m)\in \mathbb S^2\times
\mathbb S^2$, $n\perp m$, to describe the orientation field of a nematic liquid crystal, and  considers the elastic
energy of $(n, m)$ to be of the Oseen-Frank type, namely,
\begin{equation}
\mathcal{E}(n,m)=\int_\Omega \mathcal{W}(n,m,\nabla n, \nabla m)\,dx,
\end{equation}
where
\begin{eqnarray}
\mathcal{W}(n,m,\nabla n, \nabla m)&=& \frac{k_1}{2}(\nabla\cdot n)^2+\frac{k_2}{2}(n\cdot (\nabla\times n))^2+\frac{k_3}{2}|n\times(\nabla\times n)|^2\nonumber\\
&+&\frac{k_4}{2}(\nabla\cdot m)^2+\frac{k_5}{2}(m\cdot (\nabla\times m))^2+\frac{k_6}{2}|m\times(\nabla\times m)|^2\nonumber\\
&+&\frac{k_7}{2}(n\cdot(m\times(\nabla\times m)))^2+
\frac{k_8}{2}(m\cdot(n\times(\nabla\times n)))^2\nonumber\\
&+&\frac{k_9}{2}(m\cdot(\nabla\times (n\times m)))^2+
\frac{k_{10}}{2}(n\cdot(\nabla\times(m\times n)))^2\nonumber\\
&+&\frac{k_{11}}2|\nabla\times(n\times m)|^2+\frac{k_{12}}{2}(\nabla\cdot(n\times m))^2.
\end{eqnarray}
Here all the coefficients $k_i$, $1\le i\le 12$ stand for the Frank constants, which are assumed to be all nonnegative constants.  Similar to the Oseen-Frank theory of uniaxial nematics, it is a natural question to study
the existence and partial regularity of minimizing configurations of biaxial nematics within the Govers-Vertogen
theory. This is one of our motivations in this paper, and we are able to establish a partial regularity theorem
analogous to the one on uniaxial nematics by \cite{HKL}. 

Set $\alpha_1=\min\big\{k_1, k_2, k_3\big\}>0$ and $\alpha_2=\min\big\{k_4, k_5, k_6\big\}>0$. Define
the modified Oseen-Frank energy density function of $(n, m)$ by adding the null-Lagrange terms for $n$ and
$m$, i.e.,
\begin{equation}\label{OFD_2}
\widetilde{\mathcal{W}}(n, m,\nabla n, \nabla m)=\mathcal{W}(n, m,\nabla n, \nabla m)+\alpha_1 ({\rm{tr}}(\nabla n)^2-(\nabla\cdot n)^2)
+\alpha_2 ({\rm{tr}}(\nabla m)^2-(\nabla\cdot m)^2),
\end{equation}
and the modified Oseen-Frank energy functional then reads as
\begin{equation}\label{OF_2}
\widetilde{\mathcal{E}}(n,m)=\int_\Omega \widetilde{\mathcal{W}}(n, m,\nabla n,\nabla m)\,dx.
\end{equation}
It is readily seen that there exists a constant $C>0$ depending on $k_i$ for $1\le i\le 12$ such that 
\begin{equation}\label{coer1}
\frac12\big(\alpha_1 |\nabla n|^2+\alpha_2 |\nabla m|^2\big)
\le \widetilde{\mathcal{W}}(n, m,\nabla n, \nabla m)\le C(|\nabla n|^2+|\nabla m|^2).
\end{equation}
Furthermore, it follows from the null Lagrangian property of the two added terms
$({\rm{tr}}(\nabla n)^2-(\nabla\cdot n)^2)$ and $({\rm{tr}}(\nabla m)^2-(\nabla\cdot m)^2)$
(see, Hardt-Kinderleherer-Lin \cite{HKL}) that 
\begin{equation}\label{null1}
\mathcal{E}(n, m)=\widetilde{\mathcal{E}}(n,m)
\end{equation}
holds for any $(n,m)\in \mathcal{A}(\Omega; (n_0, m_0))$, a configuration function space 
to be defined below. 

Set
$$
\mathcal{A}(\Omega)=\Big\{(n, m)\in H^1(\Omega, \mathbb S^2)\times H^1(\Omega, \mathbb S^2)\ \big|\ 
\ n(x)\cdot m(x)=0 \ {\rm{a.e.}}\ x\in\Omega\Big\}.
$$
For a given pair of maps $(n_0, m_0)\in \mathcal{A}(\Omega)$, define 
$$
\mathcal{A}(\Omega; (n_0,m_0))=\Big\{(n,m)\in\mathcal{A}(\Omega)\  \big|\ (n, m)=(n_0, m_0)\ {\rm{on}}\ \partial\Omega\Big\}.
$$
We call a pair of maps $(n,m)\in\mathcal{A}(\Omega; (n_0, m_0))$ is a minimizer of the
Oseen-Frank energy  $\mathcal{E}$, if
\begin{equation}\label{OFM1}
\mathcal{E}(n, m)\le \mathcal{E}(\tilde{n}, \tilde{m}), \ \forall (\tilde{n}, \tilde{m})\in \mathcal{A}(\Omega; (n_0, m_0)).
\end{equation}
It follows from \eqref{null1} that $(n,m)\in\mathcal{A}(\Omega; (n_0, m_0))$ is a minimizer of
the Oseen-Frank energy $\mathcal{E}$ if and 
only if it is a minimizer of the modified Oseen-Frank energy $\widetilde{\mathcal{E}}$, i.e.
\begin{equation}\label{OFM2}
\widetilde{\mathcal{E}}(n, m)\le 
\widetilde{\mathcal{E}}(\tilde{n}, \tilde{m}), \ \forall (\tilde{n}, \tilde{m})\in \mathcal{A}(\Omega; (n_0, m_0)).
\end{equation}

Now we are ready to state the first of our main theorems.
 
\begin{theorem} \label{partial_reg} For any $(n_0, m_0)\in\mathcal{A}(\Omega)$, if $(n, m)\in \mathcal{A}(\Omega; (n_0, m_0))$ is a minimizer of the biaxial Oseen-Frank energy functional $\mathcal{E}$, then there exists a closed subset $\Sigma\subset\Omega$, with $\mathcal{H}^1(\Sigma)=0$, such that
$(n,m)\in C^\infty(\Omega\setminus\Sigma, \mathbb S^2\times\mathbb S^2).$
\end{theorem}

The ideas of proof of Theorem \ref{partial_reg} are: (i) we employ a Luckhaus type extension Lemma
\cite{Luck1, Luck2} in place of the Hardt-Lin extension Lemma, since unlike $\mathbb{S}^2$ the target manifold
$\mathcal{N}$ given by \eqref{N-space} below is not simply connected, (ii) we perform a blowup argument at point
$x$ where the renormalized energy is small, which depends on the smoothness of the blowing up limit equation 
of the Euler-Lagrange equation of $(n, m)$:
\begin{equation}\label{biaxial_EL}
\displaystyle
\frac{\partial\mathcal{W}}{\partial\nabla n}-\frac{\partial\mathcal{W}}{\partial n}=\lambda_1 n+\mu m,
\ \ 
\frac{\partial\mathcal{W}}{\partial\nabla m}-\frac{\partial\mathcal{W}}{\partial m}=\lambda_2 m+\mu n,
\end{equation}
where $\lambda_1, \lambda_2,$ and $\mu$ are Lagrange multipliers. A careful examination of  the blowing up limit 
equation of \eqref{biaxial_EL} shows that it a constant coefficient strongly elliptic system, which ensures the smoothness of its solutions.  

We would like to remark that by an argument similar to the interior one can also show a boundary partial  regularity
of minimizers $(n,m)$ of $\mathcal{E}$, provided the boundary value $(n_0,m_0)\in C^\infty(\partial\Omega)$,
and we leave this to the interested readers. Because $\Pi_1(\mathcal{N})=\mathbb Z_2$, it is unknown whether
the dimension of the singular set $\Sigma$ can be less than one. 
We also point out that in the one constant approximation case
$k_1=\cdots=k_6=1$ and $k_7=\cdots=k_{12}=0$, the modified Oseen-Frank energy functional reduces
to the Dirichelet energy functional for maps $(n, m):\Omega\to \mathcal{N}$ so that the minimizer becomes
a minimizing harmonic map, which also solves the equation
\begin{equation}\label{biaxial_HM}
\begin{cases}
-\Delta n=|\nabla n|^2n+\langle\nabla n, \nabla m\rangle m,\\
-\Delta m=|\nabla m|^2m+\langle\nabla n, \nabla m\rangle n.
\end{cases}
\end{equation}
In this case, the set of defects is an discrete set, and it would be an interesting question to study the
blowup profile at a singular point.

\subsection{Introduction and main theorem on hydrodynamics of biaxial nematics}  In this subsection,
we will consider a simplified Ericksen-Leslie system of biaxial nematics $(n, m)$ in which the Oseen-Frank
energy density $\mathcal{W}$ is replaced by $e(n, m)=\frac12(|\nabla n|^2+|\nabla m|^2)$ so that the
resulting equation  for $u:\Omega\times (0,\infty)\to \mathbb R^d, n, m: \Omega\times (0,\infty)
\to \mathbb S^2$ with $n\perp m$, takes the form
 \begin{equation}\label{biaxial-nlcf0}
\begin{cases}
u_t+u\cdot\nabla u-\nu\Delta u+\nabla P=-{\rm{div}}(\nabla n\odot\nabla n+\nabla m\odot\nabla m) ,\\
{\rm{div}}\ u=0,\\ 
n_t+u\cdot \nabla n=\Delta n+|\nabla n|^2 n+\langle\nabla n, \nabla m\rangle m,\\
m_t+u\cdot \nabla m=\Delta m+|\nabla m|^2 m+\langle\nabla n, \nabla m\rangle n,
\end{cases}
\end{equation}
subject to the initial and boundary condition:
\begin{equation}\label{IBC}
\begin{cases}
(u, n, m)\big|_{t=0}=(u_0, n_0, m_0),\\
(u, n, m)\big|_{\partial\Omega}=(0, n_0, m_0).
\end{cases}
\end{equation}
Here $u_0$ satisfies ${\rm{div}} u_0=0$ and $u_0=0$ on $\partial\Omega$,
and $(n_0, m_0)\in \mathcal{A}(\Omega)$.

Set
$$\mathbf{H}=L^2{\rm{-closure\ of }} \big\{ v\in C^\infty_0 (\Omega, \mathbb R^d): \ {\rm{div}}\ v=0 \big\},$$
and
$$\mathbf{J}=H^1{\rm{-closure\ of }} \big\{ v\in C^\infty_0 (\Omega, \mathbb R^d): \ {\rm{div}}\ v=0 \big\}.$$
When $\Omega\subset\mathbb R^2$ is a bounded, smooth domain, we are able to extend the main
result of \cite{LLW10} and obtain the global existence of a Leray-Hopf weak solution of \eqref{biaxial-nlcf0}, which is smooth away from finitely many singular times. More precisely, we have

\begin{theorem}\label{LH} For $u_0\in {\bf H}$, ${n}_0, {m}_0\in H^1(\Omega, \mathbb S^2)\cap C^{2+\alpha}(\partial\Omega, \mathbb S^2)$ for some $\alpha\in (0,1)$, with
${n}_0\cdot{m}_0=0$ in $\Omega$, there exists a global weak solution $u\in L^\infty(\mathbb R_+, {\bf H})\cap L^2(\mathbb R_+, {\bf J})$ and
${n}, {m}\in L^\infty(\mathbb R_+, H^1(\Omega,\mathbb S^2))$ to \eqref{biaxial-nlcf0} and \eqref{IBC}. Moreover, the
following properties hold:
\begin{itemize}
\item[(i)] There exists $L\in\mathbb N$ depending only on $(u_0, {n}_0, {m}_0)$ and $0<T_1< ...<T_L<\infty$ such that
$$(u, {n}, {m})\in C^\infty(\Omega\times ((0,\infty)\setminus\{T_i\}_{i=1}^L))\cap C^{2+\alpha,1+\frac{\alpha}2}(\overline\Omega\times ((0,\infty)\setminus\{T_i\}_{i=1}^L)).$$
\item [(ii)] Each singular time $T_i$, $1\le i\le L$,  can be characterized by
\begin{equation}\label{concentration0}
\limsup_{t\uparrow T_i^{-}}\max_{x\in\overline\Omega}\int_{\Omega\cap B_r(x)}(|u|^2+|\nabla{n}|^2+|\nabla{m}|^2)(y,t)\,dy\ge {\mathcal{C}}_0,\  \forall r>0.
\end{equation}
Here ${\mathcal{C}}_0>0$ denotes the infimum of energies of biaxial bubbles in $\mathbb S^2$ given by
\begin{equation}\label{energy_sp}
{\mathcal{C}}_0=\inf\big\{\int_{\mathbb S^2}(|\nabla{n}|^2+|\nabla{m}|^2)\,d\sigma: ({n}, {m})\in \mathcal{B}\big\},
\end{equation}
where $\mathcal{B}$ denotes the set of all biaxial bubbles, i.e. ${n}, {m}\in C^\infty(\mathbb S^2,\mathbb S^2)$,
with $n\perp m$, is a nontrivial solution of
\begin{equation}\label{biaxial-bubble}
\begin{cases}-\Delta{n}=|\nabla {n}|^2{n}+\langle\nabla{n},\nabla{m}\rangle {m}, \\
-\Delta{m}=|\nabla {m}|^2{m}+\langle\nabla{n},\nabla{m}\rangle {n}.
\end{cases}
\end{equation}
Moreover, there exist $x_k^i\rightarrow x_0^i\in\overline\Omega$, $t_k^i\uparrow T_i$, $r_k^i\downarrow 0$ and a biaxial bubble $(\omega_i^1, \omega_i^2)\in \mathcal{B}$
such that
$$(u_k^i, {n}_k^i, {m}_k^i)\rightarrow (0, \omega_i^1, \omega_i^2) \ {\rm{in}}\ C^2_{\rm{loc}}(\mathbb R^2\times \mathbb R_{-}),$$
where
$$(u_k^i, {n}_k^i, {m}_k^i)(x,t)=(r_k^iu, {n}, {m})(x_k^i+r_k^i x, t_k^i+(r_k^i)^2 t).$$
\item [(iii)] For any $\epsilon>0$ and $0\le i\le L$ (with $T_0=0$), it holds
$$(|{n}_t|+|\nabla^2{n}|)+ (|{m}_t|+|\nabla^2{m}|)\in L^2(\Omega\times [T_i, T_{i+1}-\epsilon]), 
\ |u_t|+|\nabla^2 u|\in L^\frac43(\Omega\times [T_i, T_{i+1}-\epsilon]),$$
and for any $T\in (T_L, \infty)$, it holds that 
$$(|{n}_t|+|\nabla^2{n}|)+ (|{m}_t|+|\nabla^2{m}|)\in L^2(\Omega\times [T_L, T]), 
\ |u_t|+|\nabla^2 u|\in L^\frac43(\Omega\times [T_L, T]).$$
\item [(iv)] There exist $t_k\uparrow\infty$ and a biaxial harmonic map $({n}_\infty, {m}_\infty)\in \big(C^\infty(\Omega, \mathbb S^2)\cap C^{2+\alpha}(\overline\Omega, \mathbb S^2)\big)^{\otimes 2}$, with ${n}_\infty\perp {m}_\infty$ in $\Omega$  and $({n}_\infty, {m}_\infty)=({n}_0, {m}_0)$ on $\partial\Omega$, and a nonnegative integer
$K$ and biaxial bubbles $(\omega_i^1, \omega_i^2)\in\mathcal{B}$ for $1\le i\le K$ such that
$u(t_k)\rightarrow 0$ in $L^2(\Omega)$, $({n}(t_k), {m}(t_k))\rightharpoonup  ({n}_\infty, {m}_\infty)$ in $H^1(\Omega,\mathbb S^2)$, and 
\begin{equation}\label{EI1}
\lim_{k\rightarrow\infty} E({n}(t_k), {m}(t_k); \Omega)=E({n}_\infty, {m}_\infty; \Omega)+ \sum_{i=1}^K E(\omega_i^1, \omega_i^2; \mathbb S^2).
\end{equation}
Here $E(n, m; U)=\frac12\int_U(|\nabla n|^2+|\nabla m|^2)\,dx$ for $U=\Omega$ or $\mathbb S^2$.
\item [(v)] If $E(u_0, {n}_0, {m}_0; \Omega)=\frac12\int_U(|u_0|^2+|\nabla n_0|^2+|\nabla m_0|^2)\,dx\le {\mathcal{C}}_0$, then $(u, {n}, {m})\in C^\infty(\Omega\times (0,\infty))\cap C^{2+\alpha,1+\frac{\alpha}2}(\overline\Omega\times (0,\infty))$.
Moreover, there exist $t_k\uparrow\infty$ and a biaxial harmonic map $({n}_\infty, {m}_\infty)\in \big(C^\infty(\Omega, \mathbb S^2)\cap C^{2+\alpha}(\overline\Omega, \mathbb S^2)\big)^{\otimes 2}$, with ${n}_\infty\perp {m}_\infty$ in $\Omega$  and $({n}_\infty, {m}_\infty)=({n}_0, {m}_0)$ on $\partial\Omega$
such that
$$(u(t_k), {n}(t_k), {m}(t_k))\rightarrow (0, {n}_\infty, {m}_\infty) \ {\rm{in}}\ C^2(\overline\Omega).$$

\item [(vi)] The solution is unique in the class of all weak solutions $(u, {n}, {m})\in L^2_tH^1_x\cap L^\infty_t L^2_x(\Omega\times [0, \infty),\mathbb R^2)\times L^\infty_t H^1_x(\Omega\times [0,\infty),\mathbb S^2)\times  L^\infty_t H^1_x(\Omega\times [0,\infty),\mathbb S^2)$, that satisfy
\begin{itemize}
\item [(a)] $E(u(t), {n}(t), {m}(t);\Omega)$ is monotone nonincreasing for $t\in [0, \infty)$, and
\item [(b)] ${n}_t+u\cdot\nabla{n},\ {m}_t+u\cdot\nabla{m}\in L^2([0,T), L^2(\Omega))$
for any $0<T<\infty$.
\end{itemize}
\end{itemize}
\end{theorem}

The ideas of proof of the existence part of Theorem \ref{LH} are similar to that by \cite{LLW10}.
The crucial ingredient utilizes a priori estimate under the small $L^4$-norm condition, of which we
give a different proof  in this paper that is based on  a blowing up argument. Our assumption for the
uniqueness as well its proof is slightly different from that by \cite{LW2}, namely we first show the
$L^{2}H^2$-regularity of $(n, m)$, which is needed for the uniqueness, as a consequence of this
new condition that is a seemingly weaker one.   We would like to point out that Theorem \ref{LH}
likely holds for time-dependent boundary values for $(n, m)$. We refer the interested readers
to the paper \cite{LWZZ} where the uniaxial case has been considered. 

The paper is organized as follows. In Section 2, we will prove Theorem \ref{partial_reg}.  
In Section 3, we will provide some preliminary results on \eqref{biaxial-nlcf0}. In Section 4,
we will sketch the existence of local smooth solutions of \eqref{biaxial-nlcf0}.
In Section 5, we will establish a priori estimate of \eqref{biaxial-nlcf0} under the small $L^4$-condition.
In Section 6, we will prove the existence part of Theorem \ref{LH}. Finally, in Section 7
we will prove the uniqueness part of Theorem \ref{LH}.

\section{Proof of Theorem \ref{partial_reg}}
\setcounter{section}{2} \setcounter{equation}{0}

This section is devoted to proof of Theorem \ref{partial_reg}. 
First we need to recall the following  extension Lemma, originally due to Luckhaus \cite{Luck1} (see also \cite{HKL}).

\begin{lemma} \label{lucks} For any $\Lambda>0$ and $0<\lambda<1$, there exists an $\epsilon=\epsilon(\Lambda,\lambda)>0$ such that for $i=1,2$, if
$(n_i, m_i)\in H^1(\partial B_1, \mathbb S^2)$, with $n_i\cdot m_i=0$ a.e. in $\partial B_1$, satisfies
\begin{equation}\label{bdry_cond}
\sum_{i=1}^2\int_{\partial B_1} (|\nabla_{\rm{tan}} n_i|^2+|\nabla_{\rm{tan}} m_i|^2)\,dH^2\le \Lambda,
\  \int_{\partial B_1} (|n_1-n_2|^2+|m_1-m_2|^2)\,dH^2\le\epsilon,
\end{equation}
then there exists a map $(N, M)\in H^1(B_1\setminus B_{1-\lambda},
\mathbb S^2\times \mathbb S^2)$ such that
\begin{equation}\label{extension_p}
 \begin{cases}
 N\cdot M=0 \ {\rm{in}}\ B_1\setminus B_{1-\lambda},\\
 (N, M)(x)=(n_2, m_2)(x) \ {\rm{on}}\ \partial B_1,\\
  (N, M)(x)=(n_1, m_1)(\frac{x}{1-\lambda}) \ {\rm{on}}\ \partial B_{1-\lambda},
  \end{cases}
  \end{equation}
  and
  \begin{eqnarray}\label{extension_est}
 \int_{B_1\setminus B_{1-\lambda}} (|\nabla N|^2+|\nabla M|^2)\,dx&\le& 
 C\lambda \sum_{i=1}^2\int_{\partial B_1} (|\nabla_{\rm{tan}} n_i|^2+|\nabla_{\rm{tan}} m_i|^2)\,dH^2\nonumber\\
 &+&C\lambda^{-1} \int_{\partial B_1} (|n_1-n_2|^2+|m_1-m_2|^2)\,dH^2.
\end{eqnarray}
\end{lemma}
\pf First define the submanifold $\mathcal{N}\subset\mathbb S^2\times\mathbb S^2$ by letting 
\begin{equation}\label{N-space}
\mathcal{N}=\Big\{(p,q)\in \mathbb S^2\times \mathbb S^2\ \big|\ p\cdot q=0 \Big\}.
\end{equation}
Let $SO(3)$ denote the special orthogonal group of order $3$. Then it is readily seen that $\mathcal{N}$ is differeomorphic to $SO(3)$
through  the map $\Phi(p,q)=(p, q, p\times q):\mathcal{N}\to SO(3)$. In particular, the first homotopy group of $\mathcal{N}$, $\Pi_1(\mathcal{N})=\Pi_1(SO(3))=\mathbb Z_2$.
It is easy to check that for any bounded smooth domain $O\subset\mathbb R^3$, let $E=\partial O$ or $O$,  
then $(n, m)\in H^1(E,\mathbb S^2)\times H^1(E,\mathbb S^2)$ satisfies $n\cdot m=0$  in $E$ if and only if
$U=(n, m)\in H^1(E, \mathcal{N})$.  Since $U_i=(n_i, m_i)\in H^1(\partial B_1, \mathcal{N})$ for $i=1,2$, it follows from Luckhaus \cite{Luck1, Luck2} that
there exists a map $V\in H^1(B_1\setminus B_{1-\lambda}, \mathbb R^6)$ such that
$$V(x)=U_2(x) \ {\rm{on}}\ \partial B_1, \ \ V(x)=U_1(\frac{x}{1-\lambda}) \ {\rm{on}}\ \partial B_{1-\lambda},$$
$$\int_{B_1\setminus B_{1-\lambda}} |\nabla V|^2\,dx\le
C\lambda \int_{\partial B_1} (|\nabla_{\rm{tan}} U_1|^2+|\nabla_{\rm{tan}} U_2|^2)\,dH^2
+C\lambda^{-1}\int_{\partial B_1}|U_1-U_2|^2\,dH^2,$$
and
$${\rm{dist}}(V(x), \mathcal{N})\le C\Lambda \big(\frac{\epsilon}{\lambda^2}\big)^\frac14, 
\ {\rm{a.e.}}\ x\in B_1\setminus B_{1-\lambda}.$$
Recall that there exist $\delta=\delta(\mathcal{N})>0$ and a smooth nearest point projection 
map 
$$\Pi:\mathcal{N}_\delta=\big\{y\in\mathbb R^6: {\rm{dist}}(y, \mathcal{N})<\delta\big\}\to \mathcal{N}.$$
Thus by choosing a sufficiently small $\epsilon=\epsilon(\Lambda,\lambda)>0$, we have that
$V(B_1\setminus B_{1-\lambda})\subset \mathcal{N}_\delta$ and hence 
$(N, M)=U=\Pi(V):B_1\setminus B_{1-\lambda}\to \mathcal{N}$ will be the desired extension map.
\qed

\bigskip

A crucial step in the proof of Theorem \ref{partial_reg} is  the following $\epsilon_0$-decay Lemma, which will be proved by a blowing up argument here.
\begin{lemma}\label{epsilon-reg} There exist an $\epsilon_0>0$ and $\theta_0\in (0,\frac12)$, depending only on $k_i$ for $i=1,\cdots, 12$, such that 
if $(n,m)\in \mathcal{A}(\Omega; (n_0, m_0))$ is a minimizer of the modified Oseen-Frank energy
$\widetilde{\mathcal{E}}$, that satisfies, for a ball $B_r(x_0)\subset \Omega$,
\begin{equation}\label{epsilon1} r^{-1}\int_{B_r(x_0)} (|\nabla n|^2+|\nabla m|^2)\,dx\le \epsilon_0^2,
\end{equation}
then
\begin{equation}\label{epsilon2} 
(\theta_0 r)^{-1}\int_{B_{\theta_0 r}(x_0)} (|\nabla n|^2+|\nabla m|^2)\,dx\le \theta_0 r^{-1}\int_{B_r(x_0)} (|\nabla n|^2+|\nabla m|^2)\,dx.
\end{equation}
\end{lemma}

\pf  We argue it by contradiction. Suppose that for any $\theta\in (0,\frac12)$, there exists a minimizer $(n, m)\in \mathcal{A}(\Omega; (n_0, m_0))$ of the modified Oseen-Frank energy functional
$\widetilde{\mathcal{E}}$, and a sequence $\{x_i\}\subset\Omega$ and $\{r_i\}\in (0, {\rm{dist}}(x_i,\partial\Omega))$ such that 
\begin{equation}\label{epsilon3} r_i^{-1}\int_{B_{r_i}(x_i)} (|\nabla n|^2+|\nabla m|^2)\,dx=\epsilon_i^2\rightarrow 0,
\end{equation}
but 
\begin{equation}\label{epsilon4} 
(\theta r_i)^{-1}\int_{B_{\theta r_i}(x_i)} (|\nabla n|^2+|\nabla m|^2)\,dx >\theta r_i^{-1}\int_{B_{r_i}(x_i)} (|\nabla n|^2+|\nabla m|^2)\,dx=\theta \epsilon_i^2.
\end{equation}
Define the rescaled sequence of maps $(n_i, m_i)\in \mathcal{A}(B_1)$ by letting
$$(n_i, m_i)(x)=(n,m)(x_i+r_i x), \ x\in B_1.$$
Then it is easy to see that $(n_i, m_i)\in \mathcal{A}(B_1)$ is a sequence of minimizers of the modified Oseen-Frank energy in the sense that
$$\int_{B_1}\widetilde{{\mathcal{W}}}(n_i, m_i,\nabla n_i, \nabla m_i)\,dx\le\int_{B_1} \widetilde{\mathcal{W}}(n, m,\nabla n, 
\nabla m)\,dx,$$
whenever $(n, m)\in \mathcal{A}(B_1)$ satisfies $(n, m)=(n_i, m_i)$ on $\partial B_1$.
Furthermore, it holds that 
\begin{equation}\label{epsilon3} \int_{B_1} (|\nabla n_i|^2+|\nabla m_i|^2)\,dx=\epsilon_i^2\rightarrow 0,
\end{equation}
and
\begin{equation}\label{epsilon4}
\theta^{-1}\int_{B_\theta} (|\nabla n_i|^2+|\nabla m_i|^2)\,dx
>\theta  \int_{B_1} (|\nabla n_i|^2+|\nabla m_i|^2)\,dx=\theta\epsilon_i^2.
\end{equation}

For $f:B_1\to\mathbb R^3$, let $\overline{f}=\frac{1}{|B_1|}\int_{B_1} f\,dx$ denote the average of
$f$ over $B_1$.  Since $(n_i, m_i): B_1\to \mathcal{N}$,  it follows from the Poincar\`e inequality that
$${\rm{dist}}^2\big((\overline{n_i}, \overline{m_i}), \mathcal{N}\big)
\le \frac{1}{|B_1|} \int_{B_1} (|n_i-\overline{n_i}|^2+|m_i-\overline{m_i}|^2)\,dx
\le C \int_{B_1}(|\nabla n_i|^2+|\nabla m_i|^2)\le C\epsilon_i^2.
$$
Thus there exists a point $(p_i, q_i)\in\mathcal{N}$ such that
\begin{equation}\label{proj}
|\overline{n_i}-p_i|^2+|\overline{m_i}-q_i|^2=
{\rm{dist}}^2\big((\overline{n_i}, \overline{m_i}), \mathcal{N}\big)\le C\epsilon_i^2.
\end{equation}
Define the blow-up sequence
$$u_i=\frac{n_i-p_i}{\epsilon_i},  \ v_i=\frac{m_i-q_i}{\epsilon_i}\ \ \ {\rm{in}}\ \ \ B_1.$$
It follows from \eqref{epsilon3}, the Poincar\`e inequality, and \eqref{proj}  that
$$\sup_{i}\big\{\|u_i\|_{H^1(B_1)}+\|v_i\|_{H^1(B_1)}\big\}\le  C.$$
From \eqref{epsilon4}, we also have
\begin{equation}\label{epsilon5}
\theta^{-1}\int_{B_\theta} (|\nabla u_i|^2+|\nabla v_i|^2)\,dx>\theta.
\end{equation}
Hence we may assume that there exists a pair of maps $(u,v)\in H^1(B_1, \mathbb R^6)$
such that, after passing to a subsequence, 
$$(u_i, v_i)\rightharpoonup (u,v)\ {\rm{in}}\ H^1(B_1,\mathbb R^6),\ \   (u_i, v_i)\rightarrow (u,v)\ {\rm{in}}\ L^2(B_1,\mathbb R^6) \ {\rm{and\ a.e.\ in}}\ B_1.
$$
We may assume that there exists a point $(p,q)\in \mathcal{N}$, i.e. $p, q\in \mathbb S^2$ with
$p\perp q$, such that $(p_i, q_i)\rightarrow (p,q)$. 
Since $T_{(p_i, q_i)}\mathcal{N}\rightarrow T_{(p, q)}\mathcal{N}$, it is not hard to see
that $(u,v): B_1\to T_{(p,q)}\mathcal{N}$. This is equivalent to
that $u\in T_p\mathbb{S}^2$, $v\in T_q \mathbb{S}^2$
in $B_1$, and $u\cdot q+p\cdot v=0$ in $B_1$. 

Next we want to show that $(u,v)$ is a minimizer of the following energy functional, which can be viewed as the 
limit of the modified Oseen-Frank energy functional  $\widetilde{\mathcal{E}}$ under the blowing up process. 
Define
\begin{equation}\label{blowup_f}
\mathcal{F}(u,v)=\int_{B_1} e(\nabla u,\nabla v)\,dx,
\end{equation}
where
\begin{eqnarray}\label{new_density1}
e(\nabla u,\nabla v)&=& \frac{k_1}{2}(\nabla\cdot u)^2+\frac{k_2}{2}(p\cdot (\nabla\times u))^2+\frac{k_3}{2}|p\times(\nabla\times u)|^2\nonumber\\
&+&\frac{k_4}{2}(\nabla\cdot v)^2+\frac{k_5}{2}(q\cdot (\nabla\times v))^2+\frac{k_6}{2}|q\times(\nabla\times v)|^2\nonumber\\
&+& \frac{\alpha_1}2 ({\rm{tr}}(\nabla u)^2-(\nabla\cdot u)^2)+
 \frac{\alpha_2}2 ({\rm{tr}}(\nabla v)^2-(\nabla\cdot v)^2)\nonumber\\
&+&\frac{k_7}{2}(p\cdot(q\times(\nabla\times v)))^2+
\frac{k_8}{2}(q\cdot(p\times(\nabla\times u)))^2\nonumber\\
&+&\frac{k_9}{2}(q\cdot(\nabla\times (p\times v+u\times q)))^2+
\frac{k_{10}}{2}(p\cdot(\nabla\times(q\times u+v\times p)))^2\nonumber\\
&+&\frac{k_{11}}2|\nabla\times(p\times v+u\times q)|^2+\frac{k_{12}}{2}(\nabla\cdot(p\times v+u\times q))^2.
\end{eqnarray}
For $(p,q)\in \mathcal{N}$, or $p, q\in \mathbb S^2$ with $p\cdot q=0$, 
we define the corresponding configuration space 
\begin{eqnarray*}
\mathcal{B}(B_1)&=& H^1(B_1, T_{(p,q)} \mathcal{N})\\
&=&\Big\{(u,v)\in H^1(B_1, T_p\mathbb S^2)\times H^1(B_1, T_q\mathbb S^2)\ \big|
\ u\cdot q+p\cdot v=0 \ {\rm{a.e.\ in}}\ B_1\Big\}.
\end{eqnarray*}
Then we have
\begin{lemma} \label{strong-conv} For any fixed $0<\lambda<1$, it holds that
\begin{equation}\label{OF-convergence}
\lim_{i\rightarrow\infty}\int_{B_{1-\lambda}} \widetilde{\mathcal{W}}(p_i+\epsilon u_i, q_i+\epsilon v_i,\nabla u_i, \nabla v_i)\,dx
=\int_{B_{1-\lambda}} e(\nabla u,\nabla v)\,dx.
\end{equation}
Moreover, 
$(u,v)\in \mathcal{B}(B_1)$ is a minimizer of the energy functional $\mathcal{F}$ over the space
$\mathcal{B}(B_1)$ in the sense that 
\begin{equation}\label{mini}
\mathcal{F}(u,v)=\int_{B_1} e(\nabla u,\nabla v)\,dx\le \mathcal{F}(\tilde{u}, \tilde{v})=
\int_{B_1} e(\nabla \tilde{u}, \nabla \tilde{v})\,dx
\end{equation}
holds for any comparison map $(\tilde{u}, \tilde{v})\in \mathcal{B}(B_1)$ 
satisfying $(\tilde{u}, \tilde{v})=(u,v)$ in $B_1\setminus B_{1-\lambda}$.
\end{lemma}
\pf For any small $0<\lambda<1$, take a comparison map $(\tilde{u}, \tilde{v})\in \mathcal{B}(B_1)$ satisfying 
$(\tilde{u}, \tilde{v})=(u,v)$ in $B_1\setminus B_{1-\lambda}$. By Fatou's lemma and Fubini's theorem, there exists
$\rho_0\in (1-\lambda, 1)$ such that
$$\delta_i^2=\int_{\partial B_{\rho_0}} (|u_i-\tilde{u}|^2+|v_i-\tilde{v}|^2)\,dH^2\rightarrow 0,
\ \int_{\partial B_{\rho_0}} (|\nabla u_i|^2+|\nabla v_i|^2+|\nabla \tilde{u}|^2+|\nabla \tilde{v}|^2)\,dH^2\le C.
$$
Choose $R_i\rightarrow \infty$ such that $R_i\epsilon_i\rightarrow 0$. Define
$$
(\widehat{u_i}, \widehat{v_i})=\frac{(R_i \tilde{u}, R_i\tilde{v})}{\max\{R_i, |\tilde{u}|+|\tilde{v}|\}} \ \ {\rm{in}} \ \ B_1,
$$
and
$$(\widetilde{n_i}, \widetilde{m_i})=\Psi_{i}((p_i,q_i)+\epsilon_iT_{(p,q), (p_i,q_i)}(\widehat{u_i}, \widehat{v_i}))
\  \ {\rm{in}}\ \ B_1,$$
where $T_{(p,q), (p_i,q_i)}: T_{(p,q)}\mathcal{N}\to T_{(p_i,q_i)}\mathcal{N}$ is a smooth diffeoermorphic map
satisfying
$$\begin{cases} \big|T_{(p,q),(p_i, q_i)} (v,w)-(v,w)\big|\le O(i^{-1})(|v|+|w|),\\
\big|\nabla T_{(p,q),(p_i, q_i)} (v,w)-{\rm{Id}}\big|\le O(i^{-1})
\end{cases}
$$ for all $(v,w)\in T_{(p,q)}\mathcal{N}$, 
and $\Psi_i: \big(\{(p_i,q_i)\}+T_{(p_i,q_i)}\mathcal{N}\big)\cap B_{R_i\epsilon_i}(p_i,q_i)\to \mathcal{N}$ is a smooth map 
satisfying
$$|\Psi_i(v,w)-(v,w)|\le O(i^{-1})(|v|+|w|), \ \ |\nabla\Psi_i (v,w)-{\rm{Id}}|\le O(i^{-1})$$
for all $(v,w)\in (\{(p_i,q_i)\}+T_{(p_i,q_i)}\mathcal{N})\cap B_{R_i\epsilon_i}(p_i,q_i).$

It is not hard to verify that
\begin{equation}
\begin{cases}\displaystyle
\int_{\partial B_{\rho_0}} (|n_i-\widetilde{n_i}|^2+|m_i-\widetilde{m_i}|^2)\,dH^2\le C\delta_i^2 \epsilon_i^2,\\
\displaystyle
\int_{\partial B_{\rho_0}} (|\nabla n_i|^2+|\nabla\widetilde{n_i}|^2+|\nabla m_i|^2+|\nabla\widetilde{m_i}|^2)\,dH^2\le C\Lambda \epsilon_i^2.
\end{cases}
\end{equation}
Now we can apply Lemma \ref{lucks} to $(n_i, m_i)$ and $(\widetilde{n_i}, \widetilde{m_i})$ and conclude that there exists an
extension map $(\widehat{n_i}, \widehat{m_i})\in H^1(B_1, \mathcal{N})$ such that 
\begin{equation}
\begin{cases}\displaystyle
(\widehat{n_i}, \widehat{m_i})(x)=(\widetilde{n_i}, \widetilde{m_i})(\frac{x}{1-\delta_i}), \ \forall |x|\le\rho_0(1-\delta_i),\\
\displaystyle(\widehat{n_i}, \widehat{m_i})(x)=(n_i, m_i)(x), \  \forall \rho_0\le |x|\le 1,
\end{cases}
\end{equation}
and
\begin{eqnarray}\label{ext-est}
&&\int_{B_{\rho_0}\setminus B_{\rho_0(1-\delta_i)}}(|\nabla\widehat{n_i}|^2+|\nabla\widehat{m_i}|^2)\,dx\nonumber\\
&&\le C\delta_i  \int_{\partial B_{\rho_0}} (|\nabla n_i|^2+|\nabla\widetilde{n_i}|^2+|\nabla m_i|^2+|\nabla\widetilde{m_i}|^2)\,dH^2\nonumber\\
&&\ \ \ +C\delta_i^{-1}\int_{\partial B_{\rho_0}} (|n_i-\widetilde{n_i}|^2+|m_i-\widetilde{m_i}|^2)\,dH^2\nonumber\\
&&\le C\delta_i\epsilon_i^2.
\end{eqnarray}
From the minimality of $(n_i, m_i)$, we conclude that
\begin{eqnarray}
&&\int_{B_{\rho_0}} e(\nabla u,\nabla v)\,dx\nonumber\\
&&\le\liminf_{i\rightarrow\infty} \int_{B_{\rho_0}} \widetilde{\mathcal{W}}(p_i+\epsilon_i{u_i}, q_i+\epsilon_i {v_i}, \nabla {u_i}, \nabla {v_i})\,dx\nonumber\\
&&= \liminf_{i\rightarrow\infty}\epsilon_i^{-2} \int_{B_{\rho_0}} \widetilde{\mathcal{W}}(n_i,m_i, \nabla n_i, \nabla m_i)\,dx\nonumber\\
&&\le \liminf_{i\rightarrow\infty}\epsilon_i^{-2} \int_{B_{\rho_0}} \widetilde{\mathcal{W}}(\widehat{n_i},\widehat{m_i}, \nabla \widehat{n_i}, \nabla \widehat{m_i})\,dx\nonumber\\
&&=\liminf_{i\rightarrow\infty}\epsilon_i^{-2} \Big(\int_{B_{\rho_0(1-\delta_i)}} \widetilde{\mathcal{W}}(\widehat{n_i},\widehat{m_i}, \nabla \widehat{n_i}, \nabla \widehat{m_i})\,dx
+\int_{B_{\rho_0}\setminus B_{\rho_0(1-\delta_i)}} \widetilde{\mathcal{W}}(\widehat{n_i},\widehat{m_i}, \nabla \widehat{n_i}, \nabla \widehat{m_i})\,dx\Big)\nonumber\\
&&\le \liminf_{i\rightarrow\infty}\epsilon_i^{-2} \Big((1-\delta_i)\int_{B_{\rho_0}} \widetilde{\mathcal{W}}(\widetilde{n_i},\widetilde{m_i}, \nabla \widetilde{n_i}, \nabla \widetilde{m_i})\,dx
+C\int_{B_{\rho_0}\setminus B_{\rho_0(1-\delta_i)}} (|\nabla \widehat{n_i}|^2+|\nabla \widehat{m_i}|^2)\,dx\Big)\nonumber\\
&&\le \liminf_{i\rightarrow\infty}\epsilon_i^{-2} \Big((1-\delta_i)\int_{B_{\rho_0}} \widetilde{\mathcal{W}}(\widetilde{n_i},\widetilde{m_i}, \nabla \widetilde{n_i}, \nabla \widetilde{m_i})\,dx
+C\delta_i\epsilon_i^2\Big)\nonumber\\
&&= \liminf_{i\rightarrow\infty}\epsilon_i^{-2}\int_{B_{\rho_0}} \widetilde{\mathcal{W}}(\widetilde{n_i},\widetilde{m_i}, \nabla \widetilde{n_i}, \nabla \widetilde{m_i})\,dx. \label{lsc1}
\end{eqnarray}
Observe that
$$\big\|\widetilde{n_i}-(p_i+\epsilon_i \widehat{u_i})\big\|_{L^\infty(B_{\rho_0})}+\big\|\widetilde{m_i}-(q_i+\epsilon_i \widehat{v_i})\big\|_{L^\infty(B_{\rho_0})}\rightarrow 0,$$
and
$$|\nabla(\epsilon_i^{-1}\widetilde{n_i})-\nabla\widehat{u_i}|\le O(i^{-1})|\nabla\widehat{u_i}|,
\  |\nabla(\epsilon_i^{-1}\widetilde{m_i})-\nabla\widehat{v_i}|\le O(i^{-1})|\nabla\widehat{v_i}|.$$
We can estimate
$$
\Big|\epsilon_i^{-2}\widetilde{\mathcal{W}}(\widetilde{n_i},\widetilde{m_i}, \nabla \widetilde{n_i}, \nabla \widetilde{m_i})
-\widetilde{\mathcal{W}}(p_i+\epsilon_i\widehat{u_i}, q_i+\epsilon_i\widehat{v_i}, \nabla\widehat{u_i}, \nabla\widehat{v_i})\Big|
\le O(i^{-1})\big( |\nabla\widehat{u_i}|^2+|\nabla\widehat{v_i}|^2\big),
$$
which yields
\begin{equation}\label{lsc2}
\liminf_{i\rightarrow\infty}\Big|\epsilon_i^{-2}\int_{B_{\rho_0}} \widetilde{\mathcal{W}}(\widetilde{n_i},\widetilde{m_i}, \nabla \widetilde{n_i}, \nabla \widetilde{m_i})\,dx
-\int_{B_{\rho_0}} \widetilde{\mathcal{W}}(p_i+\epsilon_i\widehat{u_i}, q_i+\epsilon_i\widehat{v_i}, \nabla\widehat{u_i}, \nabla\widehat{v_i})\,dx\Big|=0.
\end{equation}
Observe also that
\begin{equation}\label{lsc3}
\lim_{i\rightarrow\infty} 
\int_{B_{\rho_0}} \widetilde{\mathcal{W}}(p_i+\epsilon_i\widehat{u_i}, q_i+\epsilon_i\widehat{v_i}, \nabla\widehat{u_i}, \nabla\widehat{v_i})\,dx
=\int_{B_{\rho_0}}\widetilde{\mathcal{W}}(p, q, \nabla \widetilde{u},\nabla \widetilde{v})\,dx
=\int_{B_{\rho_0}} e(\nabla \widetilde{u},\nabla\widetilde{v})\,dx.
\end{equation}
Combining \eqref{lsc3} together with \eqref{lsc2} and then substituting it into \eqref{lsc1}, we obtain that
\begin{equation}\label{lsc4}
\int_{B_{\rho_0}} e(\nabla u,\nabla v)\,dx\le\liminf_{i\rightarrow\infty}\epsilon_i^{-2} \int_{B_{\rho_0}} 
\widetilde{\mathcal{W}}(n_i,m_i, \nabla n_i, \nabla m_i)\,dx \le \int_{B_{\rho_0}} e(\nabla \widetilde{u},\nabla\widetilde{v})\,dx
\end{equation}
This clearly implies \eqref{mini}. Observe that \eqref{OF-convergence} 
follows from \eqref{lsc4}, if we choose $(\widetilde{u}, \widetilde{v})=(u,v)$. 
\qed

\bigskip
Now we want to show that the above limit map $(u,v)\in C^\infty(B_1)$.
More precisely, we have
\begin{lemma}\label{blowup-smooth} If $(u,v)\in \mathcal{B}$ is a critical point of the functional $\mathcal{F}$, then
$(u,v)\in C^\infty(B)$, and
\begin{equation}\label{smooth_estimate}
\|(u,v)\|_{C^k(B_r))}\le C\big(k, r, \|(u,v)\|_{H^1(B_1)}\big), \ \forall 0<r<1, \ \forall k\ge 1.
\end{equation}
\end{lemma}
\pf Without loss of generality, assume $p=(0,0,1)$ and $q=(1,0,0)$. Since $u\in T_p\mathbb S^2$, $v\in T_q\mathbb S^2$, and $u\cdot q+p\cdot v=0$ in $B_1$, we can write
$$u=(u^1, u^2, 0), \ v=(0, v^1, v^2), \ {\rm{and}}\ u^1+v^2=0\ \ {\rm{in}}\ B_1.$$
By direct calculations, we see the energy density function $e(\nabla u,\nabla v)$ can be written as
\begin{eqnarray}\label{new_density2}
e(\nabla u,\nabla v)&=&\frac{k_1}2(u_x^1+u_y^2)^2+\frac{k_2}2(u^2_x-u^1_y)^2+\frac{k_3}2 [(u_z^1)^2+(u^2_z)^2]\nonumber\\
&+&\frac{k_4}2(v_y^1+v_z^2)^2+\frac{k_5}2(v^2_y-v^1_z)^2+\frac{k_6}2 [(v_x^1)^2+(v^2_x)^2]\nonumber\\
&+&\frac{\alpha_1}2[(u_x^1)^2+(u^2_y)^2+2u_y^1u^2_x-(u_x^1+u_y^2)^2]\nonumber \\
&+&\frac{\alpha_2}2[(v_y^1)^2+(v^2_z)^2+2v_z^1v^2_y-(v_y^1+v_z^2)^2]\nonumber\\
&+& \frac{k_7}2 (v_x^2)^2+\frac{k_8}2 (u^1_z)^2+\frac{k_9}2 (u^2_y)^2+\frac{k_{10}}2(v^1_y)^2\nonumber\\
&+&\frac{k_{11}}2 [(u^2_y)^2+(v^1_y)^2+(v^1_z-u^2_x)^2]+\frac{k_{12}}2 (v_x^1+u^2_z)^2.
\end{eqnarray}
If we set $U=(u^1, u^2, v^1):B_1\to \mathbb R^3$ and use the fact that $v^2=-u^1$ in $B_1$,
then we can write $e(\nabla u,\nabla v)$ by 
$$e(\nabla u,\nabla v)=\sum_{i,j=1}^3\sum_{\alpha\beta=1}^3 A_{ij\alpha\beta} U^i_{\alpha} U^j_{\beta}$$
for some constant coefficients $A_{ij\alpha\beta}$ for $1\le i, j,\alpha,\beta\le 3$. Moreover, since
$$e(\nabla u,\nabla v)\ge \frac{\alpha_1}{2} |\nabla u|^2+\frac{\alpha_2}{2}|\nabla v|^2$$
holds for any $u\in H^1(B_1, \mathbb R^3)$ and $v\in H^1(B_1, \mathbb R^3)$, we can see 
$(A_{ij\alpha\beta})$ is strongly elliptic in the sense that 
\begin{equation}\label{strong-elliptic}
\sum_{i,j,\alpha,\beta=1}^3 A_{ij\alpha\beta}\xi_\alpha^i \xi_\beta^j \ge \frac{\alpha_3}2 |\xi|^2,
\ \forall \xi\in \mathbb R^{3\times 3},
\end{equation}
where $\alpha_3=\min\{\alpha_1,\alpha_2\}>0.$ 

By direct calculations, we know that if $(u,v)\in \mathcal{B}$ is a critical point of 
the functional $\mathcal{F}$, then $U\in H^1(B_1, \mathbb R^3)$ solves the following
constant coefficient strongly elliptic system:
\begin{equation}\label{linear-elliptic}
\sum_{j,\alpha,\beta=1}^3 A_{ij\alpha\beta} U^j_{\alpha\beta}=0, \ 1\le i\le 3, \ {\rm{in}}\ B_1.
\end{equation}
It follows from the standard theory on constant coefficient strongly elliptic systems that
$U\in C^\infty(B_1, \mathbb R^3)$, and \eqref{smooth_estimate} holds. \qed

\bigskip
Now we return to the proof of Lemma \ref{epsilon-reg}. It follows from  \eqref{smooth_estimate} 
that for any $0<\theta<\frac12$, 
\begin{equation}\label{limit_estimate}
\theta^{-1}\int_{B_\theta} e(\nabla u,\nabla v)\,dx\le C\theta^2\|(\nabla u,\nabla v)\|_{L^\infty(B_\frac12)}^2
\le C\theta^2.
\end{equation}
This, combined with Lemma \ref{strong-conv}, implies that
\begin{eqnarray}\label{gap}
\lim_{i\rightarrow\infty}
\int_{B_\theta} \widetilde{\mathcal{W}}(p_i+\epsilon u_i, q_i+\epsilon v_i,\nabla u_i, \nabla v_i)\,dx
=\int_{B_\theta} e(\nabla u,\nabla v)\,dx\le C\theta^3.
\end{eqnarray}
Observe that \eqref{coer1} implies 
\begin{equation}\label{coer10}
\widetilde{\mathcal{W}}(p_i+\epsilon u_i, q_i+\epsilon v_i,\nabla u_i, \nabla v_i)
\ge \frac{\alpha_1}{2} |\nabla u_i|^2+\frac{\alpha_2}{2}|\nabla v_i|^2\ge \frac{\alpha_3}{2}(|\nabla u_i|^2+|\nabla v_i|^2) \ \ {\rm{in}}\ \ B_1,
\end{equation}
where $\alpha_3=\min\{\alpha_1, \alpha_2\}>0$. Substituting \eqref{coer10} into \eqref{gap},
we obtain that if $i$ sufficiently large then
$$
\theta^{-1}\int_{B_{\theta}} (|\nabla u_i|^2+|\nabla v_i|^2) \,dx
\le \frac{4C}{\alpha_3} \theta^2\le \frac12\theta,
$$
provided $\theta$ is chosen sufficiently small. This contradicts to \eqref{epsilon5}.
The proof of Lemma \ref{epsilon-reg} is complete. \qed

\begin{theorem}\label{partial reg} If $(n,m)\in \mathcal{A}(\Omega; (n_0, m_0))$ is a minimizer of 
Oseen-Frank energy $\mathcal{E}$, then there exists a closed set $\Sigma\subset\Omega$,
with ${\mathcal{H}}^1(\Sigma)=0$, such that $(n,m)$ is analytic in $\Omega\setminus \Sigma$. 
\end{theorem}
\pf Define 
$$\Sigma=\Big\{a\in\Omega: \limsup_{r\rightarrow 0} r^{-1}\int_{B_r(a)} (|\nabla n|^2+|\nabla m|^2)\,dx>0\Big\}.$$
Since $\int_\Omega (|\nabla n|^2+|\nabla m|^2)\,dx<\infty$, it is well-known that $\mathcal{H}^1(\Sigma)=0$.

For any $a\in\Omega\setminus\Sigma$, there exists $R>0$ such that $B_{2R}(a)\subset\Omega$ and
$$R^{-1}\int_{B_{2R}(a)}(|\nabla n|^2+|\nabla m|^2)\,dx\le \epsilon_0^2,$$
which implies that
\begin{equation}\label{small_condition}
R^{-1}\int_{B_R(b)} (|\nabla n|^2+|\nabla m|^2)\,dx\le \epsilon_0^2, \ \forall b\in B_R(a).
\end{equation}
Thus, by repeatedly applying Lemma \ref{epsilon-reg}, we can conclude that
$$
r^{-1}\int_{B_r(b)}(|\nabla n|^2+|\nabla m|^2)\,dx\le \theta_0^{-2} \epsilon_0^2 R^{-1} r,\ 
\ \forall 0<r\le R,\ \forall b \in B_R(a).
$$
From this we see that $B_{R}(a)\subset\Omega\setminus \Sigma$, and by Morrey's decay
Lemma that $(n,m)\in C^{\frac12}(B_R(a))$. 

To show the higher order regularity of $(n,m)$ near $a$, we assume that $n(a)=e_3=(0,0,1)$ and
$m(a)=e_1=(1,0,0)$, and choose $r<R$ so that $n(B_r(a))\subset \mathbb S^2\cap B_{\frac12}(e_3)$
and $m(B_r(a))\subset \mathbb S^2\cap B_{\frac12}(e_1)$. On $B_r(a)$, for $n=(n^1, n^2, n^3)$
and $m=(m^1, m^2, m^3)$, we can solve
$$n^3=\phi(n^1, n^2)=\sqrt{1-(n^1)^2-(n^2)^2}, \ {\rm{and}}\ m^1=\psi(m^2, m^3)=\sqrt{1-(m^2)^2-(m^3)^2}.$$
It follows that
$$\nabla n^3=-\frac{n^1\nabla n^1+n^2\nabla n^2}{\sqrt{1-(n^1)^2-(n^2)^2}},
\ \nabla m^1=-\frac{m^2\nabla m^2+m^3\nabla m^3}{\sqrt{1-(m^2)^2-(m^3)^2}}, \ {\rm{in}}\ B_r(a).$$
Moreover, since $n\cdot m=0$, i.e.,
$$n^1\psi(m^2, m^3)+n^2 m^2+\phi(n^1, n^2) m^3=n^1m^1+n^2 m^2+n^3 m^3=0\ \ {\rm{in}}\ \ B_r(a),$$
we can uniquely solve $m^3=\mu(n^1, n^2, m^2)$ in $B_r(a)$ so that
\begin{eqnarray*}
&&\big(\sqrt{1-(n^1)^2-(n^2)^2}-\frac{n^1m^3}{\sqrt{1-(m^2)^2-(m^3)^2}}\big)\nabla m^3\\
&&=\frac{n^1m^2\nabla m^2}{\sqrt{1-(m^2)^2-(m^3)^2}}-\sqrt{1-(m^2)^2-(m^3)^2} \nabla n^1\\
&&\ \ \ -m^2\nabla n^2-n^2\nabla m^2+\frac{m^3(n^1\nabla n^1+n^2\nabla n^2)}{\sqrt{1-(n^1)^2-(n^2)^2}}
\ \ \ \ \ \ \ \ \ \ \ {\rm{in}}\  B_r(a).
\end{eqnarray*}
Testing the Euler-Lagrange equation of the modified Oseen-Frank functional by 
$(\eta^1,\eta^2, 0, 0, \eta^3, 0)\in C^\infty_0(B_r(a),\mathbb R^6)$, we can infer that
$u=(n^1, n^2, m^2)$ is a critical point for the functional $F(u)=\int_{B_r(a)} J(u,\nabla u)$, where,
for $p=(p^1, p^2, p^3)\in (\mathbb R^3)^3$ and $u=(u^1, u^2, u^3)\in\mathbb R^3$,
\begin{eqnarray*}
J(u,p)&=&\widetilde{\mathcal{W}}\Big[u^1, u^2, \sqrt{1-(u^1)^2-(u^2)^2}; \sqrt{1-(u^3)^2-(\mu(u^1, u^2, u^3))^2}, u^3, 
\mu(u^1, u^2, u^3); \\
&& p^1, p^2, -\frac{u^1p^1+u^2 p^2}{\sqrt{1-(u^1)^2-(u^2)^2}}; q^1, p^3, q^2\Big],
\end{eqnarray*}
where
$$
q^1=-\frac{u^3p^3+\mu(u^1,u^2,u^3) q^2}{\sqrt{1-(u^3)^2-\mu^2(u^1, u^2, u^3)}},
$$
and $q^2$ is given by 
\begin{eqnarray*}
&&\big(\sqrt{1-(u^1)^2-(u^2)^2}-\frac{u^1\mu(u^1,u^2,u^3)}{\sqrt{1-(u^3)^2-\mu^2(u^1,u^2,u^3}}\big)q^2\\
&&=\frac{u^1u^3 p^3}{\sqrt{1-(u^3)^2-\mu^2(u^1,u^2, u^3)}}-\sqrt{1-(u^3)^2-\mu^2(u^1,u^2,u^3)} p^1\\
&&\ \ \ -u^3p^2-u^2 p^3+\frac{\mu(u^1,u^2,u^3)(u^1 p^1+u^2 p^2)}{\sqrt{1-(u^1)^2-(u^2)^2}}.
\end{eqnarray*}
It is readily seen that $J$ is analytic on $\big\{(u,p): |u|<1\big\}$ and $J(u,\cdot)$ is a quadratic polynomial for each
$u$ with $|u|<1$. Moreover, it follow from \eqref{coer1} that
$$J(0,p)\ge \frac12{\alpha_3} |p|^2,\  \ \forall p\in (\mathbb R^3)^3.$$
Hence there exists $0<\delta<\frac12$ such that 
$$J(u,p)\ge \frac14{\alpha_3} |p|^2\  \ {\rm{on}}\   \big\{(u,p): \ |u|<\delta\big\}.$$
We can now conclude that $u\big|_{B_r(a)}$ satisfies a strongly elliptic system with analytic coefficients.
By \cite{Morrey}, $u$ is analytic on $B_r(a)$. This proves Theorem \ref{partial_reg}. \qed

\setcounter{section}{2} \setcounter{equation}{0}
\section{Preliminary results on \eqref{biaxial-nlcf0}}

In this section, we will derive both global and local energy inequality 
for regular solutions to \eqref{biaxial-nlcf0} in dimensions $d=2, 3$.

The first property for \eqref{biaxial-nlcf0} is the energy dissipation inequality. For any given pair of maps ${n}_0, {m}_0:\Omega\to\mathbb S^2$ satisfying
${n}_0\cdot{m}_0=0$ on $\Omega$.
\begin{lemma}\label{global_energy_ineq}
Assume $(u, {n}, {m})$ are smooth solutions to \eqref{biaxial-nlcf0} and \eqref{IBC}
%\begin{equation}\label{bc}
%(u, {n}, {m})=(0, {n}_0, {m}_0) \ {\rm{on}}\ \partial\Omega\times (0, T).
%\end{equation}
Then 
\begin{equation}\label{ED0}
\frac{d}{dt}\int_\Omega (|u|^2+|\nabla{n}|^2+|\nabla{m}|^2)=-2\int_\Omega (|{n}_t+u\cdot\nabla{n}|^2+|{m}_t+u\cdot\nabla{m}|^2).
\end{equation}
In particular, for any $t>0$ it holds that
\begin{eqnarray}\label{ED01}
&&\int_\Omega (|u|^2+|\nabla{n}|^2+|\nabla{m}|^2)(t)+2\int_0^t\int_\Omega (|{n}_t+u\cdot\nabla{n}|^2+|{m}_t+u\cdot\nabla{m}|^2)\nonumber\\
&&= \int_\Omega (|u_0|^2+|\nabla{n}_0|^2+|\nabla{m}_0|^2).
\end{eqnarray}

\end{lemma}
\pf Multiplying the first equation, the third, and fourth equation of \eqref{biaxial-nlcf0} by $u$, (${n}_t+u\cdot\nabla{n}$),
and (${m}_t+u\cdot\nabla{m}$) respectively, and integrating the resulting equations over $\Omega$, we can obtain 
\begin{equation}\label{ED1}
\frac{d}{dt}\int_\Omega \frac12|u|^2+\int_\Omega |\nabla u|^2=\int_\Omega (\nabla {n}\odot \nabla {n}+\nabla {m}\odot \nabla {m}):\nabla u,
\end{equation}
\begin{eqnarray}\label{ED2}
\int_\Omega |{n}_t+u\cdot\nabla{n}|^2&=&\int_\Omega (\Delta{n}+|\nabla{n}|^2{n})\cdot({n}_t+u\cdot\nabla{n})\nonumber\\
&&+\int_\Omega (\nabla {n}\cdot \nabla {m}){m}\cdot({n}_t+u\cdot\nabla{n}),
\end{eqnarray}
and
\begin{eqnarray}\label{ED3}
\int_\Omega |{m}_t+u\cdot\nabla{m}|^2&=&\int_\Omega (\Delta{m}+|\nabla{m}|^2{m})\cdot({m}_t+u\cdot\nabla{m})\nonumber\\
&&+\int_\Omega (\nabla {n}\cdot \nabla {m}){n}\cdot({m}_t+u\cdot\nabla{m}).
\end{eqnarray}
Adding \eqref{ED2} and \eqref{ED3}, and using $|{n}|^2=|{m}|^2=1$ and ${n}\cdot{m}=0$, we obtain
\begin{eqnarray}\label{ED4}
&&\int_\Omega \big(|{n}_t+u\cdot\nabla{n}|^2+ |{m}_t+u\cdot\nabla{m}|^2\big)\nonumber\\
&&=\int_\Omega (\Delta{n}+|\nabla{n}|^2{n})\cdot({n}_t+u\cdot\nabla{n})+(\Delta{m}+|\nabla{m}|^2{m})\cdot({m}_t+u\cdot\nabla{m})\nonumber\\
&&\quad+\int_\Omega (\nabla {n}\cdot \nabla {m})[({m}\cdot{n})_t+u\cdot\nabla({n}\cdot{m})]\nonumber\\
&&=\int_\Omega [\Delta{n}\cdot({n}_t+u\cdot\nabla{n})+\Delta{m}\cdot({m}_t+u\cdot\nabla{m})]\nonumber\\
&&=-\frac12\frac{d}{dt}\int_\Omega (|\nabla{n}|^2+|\nabla{m}|^2)+\int_\Omega {\rm{div}}[(u\cdot\nabla {n})\nabla{n}+(u\cdot\nabla {m})\nabla{m}]\nonumber\\
&&\ \ \ -\int_\Omega ({\nabla n}\odot{\nabla n}+{\nabla m}\odot{\nabla m}):\nabla u-\frac12\int_\Omega u\cdot\nabla (|\nabla{n}|^2+|\nabla{m}|^2)\nonumber\\
&&=-\frac12\frac{d}{dt}\int_\Omega (|\nabla{n}|^2+|\nabla{m}|^2) -\int_\Omega ({\nabla n}\odot{\nabla n}+{\nabla m}\odot{\nabla m}):\nabla u,
\end{eqnarray}
where we have used the boundary condition $u=0$ on $\partial\Omega$ and the condition ${\rm{div}}u=0$ in $\Omega$ in the last step.
It is readily seen that \eqref{ED0} follows by adding \eqref{ED1} and \eqref{ED4}. \qed

\smallskip
The second property is the local energy inequality for smooth solution of \eqref{biaxial-nlcf0}. 
\begin{lemma} \label{interior_local}
Under the same assumption as Lemma \ref{global_energy_ineq}, it holds that for any $0\le\phi\in C_0^\infty(\Omega)$
and any $0\le s<t\le T$, 
\begin{eqnarray}\label{LED0}
&&\int_\Omega \phi(|u|^2+|\nabla{n}|^2+|\nabla{m}|^2)(t)+2\int_s^t\int_\Omega \phi(|\nabla u|^2
+|{n}_t+u\cdot\nabla{n}|^2+|{m}_t+u\cdot\nabla{m}|^2)\nonumber\\
&&\le C\int_s^t\int_\Omega |\nabla\phi|\big[(|u|^2+|P-\overline{P}|+|\nabla u|+|\nabla{n}|^2+|\nabla{m}|^2)|u|
+|\nabla{n}||{n}_t|+|\nabla{m}||{m}_t|\big]\nonumber\\
&&\quad+ \int_\Omega \phi(|u|^2+|\nabla{n}|^2+|\nabla{m}|^2)(s), 
\end{eqnarray}
where $\overline{P}=\frac{1}{|{\rm{spt}}\phi|}\int_{{\rm{spt}}\phi} P$ denotes the average of $P$ over ${\rm{spt}}\phi$.
\end{lemma}
\pf
Multiplying \eqref{biaxial-nlcf0}$_1$ by $u\phi$ and integrating over $\Omega$ gives
\begin{eqnarray}\label{LED1}
&&\frac{d}{dt}\int_\Omega |u|^2\phi+2\int_\Omega |\nabla u|^2\phi\\
&&=2\int_\Omega \big[(\frac{|u|^2}2 u-\nabla u\cdot u +(P-P_\Omega) u)\cdot\nabla\phi+(\nabla {n}\odot \nabla {n}+\nabla {m}\odot \nabla {m}):(\nabla u\phi+u\otimes\nabla\phi)\big].\nonumber
\end{eqnarray}
Multiplying \eqref{biaxial-nlcf0}$_3$ by $({n}_t+u\cdot\nabla{n})\phi$ and \eqref{biaxial-nlcf0}$_4$ by $({m}_t+u\cdot\nabla{m})\phi$, adding the resulting equation togethers, and
using $|{n}|=|{m}|=1$ and ${n}\cdot{m}=0$, we obtain
\begin{eqnarray}\label{LED2}
&&\int_\Omega \big(|{n}_t+u\cdot\nabla{n}|^2+ |{m}_t+u\cdot\nabla{m}|^2\big)\phi\nonumber\\
&&=\int_\Omega [\Delta{n}\cdot({n}_t+u\cdot\nabla{n})+\Delta{m}\cdot({m}_t+u\cdot\nabla{m})]\phi\nonumber\\
&&=-\frac12\frac{d}{dt}\int_\Omega (|\nabla{n}|^2+|\nabla{m}|^2)\phi-\int_\Omega (\nabla{n} \cdot{n}_t+\nabla{m} \cdot{m}_t)\cdot\nabla\phi\nonumber\\
&&\ \ \ -\int_\Omega (\nabla{n}\odot\nabla{n}+\nabla{m}\odot\nabla{m}): (u\otimes\nabla\phi+\nabla u\phi)\nonumber\\
&&\ \ \ +\frac12\int_\Omega (|\nabla{n}|^2+|\nabla{m}|^2)u\cdot\nabla\phi.
\end{eqnarray}
Adding \eqref{LED1} and \eqref{LED2} yields
\begin{eqnarray*}
&&\frac{d}{dt}\int_\Omega (|u|^2+|\nabla{n}|^2+|\nabla{m}|^2)\phi+2\int_\Omega (|\nabla u|^2+|{n}_t+u\cdot\nabla{n}|^2+ |{m}_t+u\cdot\nabla{m}|^2)\phi\\
&&=\int_\Omega \nabla\phi\cdot\big[{|u|^2}u-2\nabla u\cdot u +2(P-P_\Omega)u-2(\nabla{n} \cdot{n}_t+\nabla{m} \cdot{m}_t)+(|\nabla{n}|^2+|\nabla{m}|^2)u\big].
\end{eqnarray*}
Integrating this equality over $[s, t]$ clearly yields \eqref{LED0}. \qed

\medskip
We also recall the following boundary version of local energy inequality. More precisely,

\begin{lemma}\label{bdry_local} Under the same assumption as Lemma \ref{global_energy_ineq}, there exists $r_0>0$ depending on
$\partial\Omega$ such that for any $x_0\in\partial\Omega$ and $0<r<r_0$, if $0\le\phi\in C_0^\infty(B_{r}(x_0))$
and $0\le s<t\le T$, then
\begin{eqnarray}\label{LED1}
&&\int_{\Omega\cap B_r(x_0)} \phi(|u|^2+|\nabla{n}|^2+|\nabla{m}|^2)(t)\nonumber\\
&& +2\int_s^t\int_{\Omega\cap B_r(x_0)}
\phi(|\nabla u|^2
+|{n}_t+u\cdot\nabla{n}|^2+|{m}_t+u\cdot\nabla{m}|^2)\nonumber\\
&&\le C\int_s^t\int_{\Omega\cap B_r(x_0)} |\nabla\phi|\big[(|u|^2+|\nabla u|+|P-\bar{P}|+|\nabla{n}|^2+|\nabla{m}|^2)|u|
+|\nabla{n}||{n}_t|+|\nabla{m}||{m}_t|\big]\nonumber\\
&&\quad+ \int_{\Omega\cap B_r(x_0)} \phi(|u|^2+|\nabla{n}|^2+|\nabla{m}|^2)(s), 
\end{eqnarray}
where $\bar{P}=\frac{1}{|{\rm{spt}}\phi|}\int_{{\rm{spt}}\phi}  P$.
\end{lemma}
\pf The argument is almost identical to that of Lemma \ref{interior_local}. Namely, multiply\eqref{biaxial-nlcf0}$_1$ by $u\phi$ and 
integrate over $\Omega\cap B_r(x_0)$, and observe that there is no boundary contribution because $u\phi=0$
on $\partial(\Omega\cap B_r(x_0))$. Multiply\eqref{biaxial-nlcf0}$_{2,3}$ by $(n_t+u\cdot\nabla n)\phi$ and
$(m_t+u\cdot\nabla m)\phi$ respectively and integrate over $\Omega\cap B_r(x_0)$, and observe that there is no boundary contributions because $n_t\phi=m_t\phi=0$ on $\partial(\Omega\cap B_r(x_0))$. 
\qed

\bigskip
We close this section by giving an elementary proof of the smooth approximation property of initial data.

\begin{lemma} \label{approx} For any bounded smooth domain $\Omega\subset\mathbb R^2$, if
${n}, {m}\in H^1(\Omega, \mathbb S^2)\cap C^{2+\alpha}(\overline\Omega, \mathbb S^2)$ satisfy ${n}\cdot{m}=0$ in $\overline\Omega$, then
there exists a sequence of maps ${n}^\epsilon, {m}^\epsilon \in C^\infty(\Omega,\mathbb S^2)\cap C^{2+\alpha}(\overline\Omega, \mathbb S^2)$ such that
$({n}^\epsilon, {m}^\epsilon)=(n, {m})$ on $\partial\Omega$, ${n}^\epsilon\cdot{m}^\epsilon=0$ in $\overline\Omega$, and
\begin{equation}\label{smooth-approx}
\lim_{\epsilon\rightarrow 0}\Big[\big\|{n}^\epsilon-{n}\big\|_{H^1(\Omega)}+\big\|{m}^\epsilon-{m}\big\|_{H^1(\Omega)}\Big]=0.
\end{equation}
\end{lemma}
\pf It follows from the standard mollification process that
there  exist ${n}_\epsilon, {m}_\epsilon \in C^\infty(\Omega,\mathbb R^3)
\cap C^{2+\alpha}(\overline\Omega,\mathbb R^3)$ such that ${n}_\epsilon, {m}_\epsilon)=(n, {m})$ on $\partial\Omega$, and
\begin{equation*}\label{smooth-approx1}
\lim_{\epsilon\rightarrow 0}\Big[\big\|{n}_\epsilon-{n}\big\|_{H^1(\Omega)}+\big\|{m}_\epsilon-{m}\big\|_{H^1(\Omega)}\Big]=0.
\end{equation*}
Moreover, since dim($\Omega)=2$, it follows from a modified Poincar\'e inequality that for a.e. $x\in\overline\Omega$, it holds 
$${\rm{dist}}^2({n}_\epsilon(x), \mathbb S^2)\le |{n}_\epsilon(x)-{n}(x)|^2\le C\int_{B_{2\epsilon}(x)\cap\Omega}|\nabla{n}|^2, $$
and
$${\rm{dist}}^2({m}_\epsilon(x), \mathbb S^2)\le |{m}_\epsilon(x)-{m}(x)|^2\le C\int_{B_{2\epsilon}(x)\cap\Omega}|\nabla{m}|^2.
$$
Thus there exists $\epsilon_0>0$, depending only on ${n}$ and ${m}$, such that for $\epsilon\le\epsilon_0$,
it holds that for a.e. $x\in\overline\Omega$, 
$$
|{n}_\epsilon(x)-{n}(x)|\le \frac1{10} \ {\rm{and}}\ |{m}_\epsilon(x)-{m}(x)|\le \frac1{10}.
$$
This, combined with $|{n}|=|{m}|=1$ and ${n}\cdot{m}=0$ in $\Omega$, implies that  for any $x\in\overline\Omega$,
$$
|{n}_\epsilon(x)|\ge \frac{9}{10}, \ |{m}_\epsilon(x)|\ge \frac{9}{10},\ {\rm{and}}\ |({n}_\epsilon\cdot{m}_\epsilon)(x)|\le \frac15.
$$
Now we set ${n}^\epsilon=\frac{{n}_\epsilon}{|{n}_\epsilon|}$. Then it is easy to see that ${n}^\epsilon\in C^\infty(\Omega,\mathbb S^2)\cap C^{2+\alpha}(\overline\Omega, \mathbb S^2)$ satisfies
${n}^\epsilon={n}$ on $\partial\Omega$, and
$$\lim_{\epsilon\rightarrow 0}\big\|{n}^\epsilon-{n}\big\|_{H^1(\Omega)}=0.$$
Next we set $\widehat{{m}_\epsilon}={m}_\epsilon-({n}^\epsilon\cdot {m}_\epsilon){n}^\epsilon$. Then we have
that $\widehat{{m}_\epsilon}={m}$ on $\partial\Omega$, and $$
{n}^\epsilon\cdot\widehat{{m}_\epsilon}=0, \forall\ x\in\Omega,
$$
and
$$|\widehat{{m}_\epsilon}(x)|\ge |{m}_\epsilon|-|{n}^\epsilon\cdot {m}_\epsilon||{n}^\epsilon|
=|m_\epsilon|-|n_\epsilon|^{-1} |m_\epsilon\cdot n_\epsilon|\ge \frac{9}{10}-\frac29=\frac{61}{90},
\ \forall\ x\in\Omega.$$
Finally we define ${m}^\epsilon=\frac{\widehat{{m}_\epsilon}}{|\widehat{{m}_\epsilon}|}$. Then 
${m}^\epsilon\in C^\infty(\Omega,\mathbb S^2)\cap C^{2+\alpha}(\overline\Omega, \mathbb S^2)$ satisfies ${m}^\epsilon={m}$ on $\partial\Omega$, 
${m}^\epsilon\cdot{n}^\epsilon=0$ in $\Omega$, and
$$\lim_{\epsilon\rightarrow 0}\big\|{m}^\epsilon-{m}\big\|_{H^1(\Omega)}=0.$$
This completes the proof. \qed

\section{Existence of local smooth solutions}
\setcounter{equation}{0}

In this section, we will sketch the proof of existence of short time smooth solutions by the fixed point argument.
First, we will prove a geometry preserving property for smooth solutions of \eqref{biaxial-nlcf0}. 

\begin{lemma} \label{constraint} For any $T>0$ and a pair of maps ${n}_0, {m}_0\in C^\infty(\Omega,\mathbb S^2)\cap C^{2+\alpha}(\overline\Omega,\mathbb S^2)$, with ${n}_0\cdot {m}_0=0$,
if ${n}, {m}\in C^\infty(\Omega\times [0,T), \mathbb R^3)\cap C^{2+\alpha,1+\frac{\alpha}2}(\overline\Omega\times [0,T), \mathbb R^3)$ solves
\begin{equation}\label{BLCF1}
\begin{cases}
{n}_t+u\cdot\nabla {n}-\Delta {n}=|\nabla {n}|^2 {n}+\langle\nabla{n}, \nabla{m}\rangle {m}, \\
{m}_t+u\cdot\nabla {m}-\Delta {m}=|\nabla {m}|^2 {m}+\langle\nabla{n}, \nabla{m}\rangle {n}
\end{cases} 
\ {\rm{in}}\  \Omega\times (0, T),
\end{equation}
for some $u\in C^\infty_0(\Omega\times (0,T),\mathbb R^2)$ with ${\rm{div}} u=0$ in $\Omega$, and
\begin{equation}\label{IBC1}
({n}, {m})=({n}_0, {m}_0) \ {\rm{on}}\ \partial_p(\Omega\times (0,T)).
\end{equation}
Then $|{n}|=|{m}|=1$ and ${n}\cdot {m}=0$ in $\Omega\times (0,T)$. Here $\partial_p(\Omega\times (0,T))$ denotes
the parabolic boundary.
\end{lemma}
\pf Multiplying \eqref{BLCF1}$_1$ by ${m}$ and \eqref{BLCF1}$_2$ by ${n}$, and adding these
two equations, we obtain
\begin{eqnarray}
&&({n}\cdot{m})_t+u\cdot\nabla ({n}\cdot{m})-\Delta({n}\cdot{m})
+2\langle\nabla{n},\nabla{m}\rangle\nonumber\\
&&=(|\nabla{n}|^2+|\nabla{m}|^2){n}\cdot{m}
+\langle\nabla{n}, \nabla{m}\rangle (|{n}|^2+|{m}|^2).
\end{eqnarray}
Multiplying this equation by ${n}\cdot{m}$ and integrating the resulting equation over $\Omega$, 
we obtain that
\begin{eqnarray*}
&&\frac{d}{dt}\int_\Omega ({n}\cdot{m})^2+\int_\Omega u\cdot\nabla ({n}\cdot{m})^2-2\int_\Omega \Delta({n}\cdot{m}) ({n}\cdot{m})\nonumber\\
&&=2\int_\Omega (|\nabla{n}|^2+|\nabla{m}|^2)({n}\cdot{m})^2\nonumber\\
&&+\int_\Omega \langle\nabla {n}, \nabla{m}\rangle (|{n}|^2+|{m}|^2-2)({n}\cdot{m}).
\end{eqnarray*}
Applying ${\rm{div}}u=0$ in $\Omega$, ${n}\cdot{m}={n}_0\cdot{m}_0=0$ on  $\partial\Omega$, and the integration by parts, 
we obtain
\begin{eqnarray}\label{nm-equation}
&&\frac{d}{dt}\int_\Omega ({n}\cdot{m})^2+2\int_\Omega |\nabla({n}\cdot{m})|^2\nonumber\\
&&=2\int_\Omega (|\nabla{n}|^2+|\nabla{m}|^2)({n}\cdot{m})^2
+\int_\Omega \langle\nabla {n}, \nabla{m}\rangle (|{n}|^2+|{m}|^2-2)({n}\cdot{m})\nonumber\\
&&\le C\int_\Omega \big[(|{n}|^2-1)^2+(|{m}|^2-1)^2+({n}\cdot{m})^2\big].
\end{eqnarray}
Multiplying  \eqref{BLCF1}$_1$ by ${n}$, we can get
$$
(|{n}|^2-1)_t+u\cdot\nabla (|{n}|^2-1)-\Delta (|{n}|^2-1)=2|\nabla {n}|^2 (|{n}|^2-1)
+\langle\nabla{n}, \nabla{m}\rangle {n}\cdot{m}.
$$
Multiplying both sides of this equation by $(|{n}|^2-1)$, integrating the resulting equation over $\Omega$,
and applying  ${\rm{div}}u=0$ in $\Omega$ and $|{n}|=|{n}_0|=1$ on $\partial\Omega$, we obtain
\begin{eqnarray}\label{n-equation}
&&\frac{d}{dt}\int_\Omega (|{n}|^2-1)^2+\int_\Omega |\nabla({|n}|^2-1)|^2\nonumber\\
&&=4\int_\Omega |\nabla {n}|^2 (|{n}|^2-1)^2
+2\int_\Omega \langle\nabla{n}, \nabla{m}\rangle ({n}\cdot{m})(|{n}|^2-1)\nonumber\\
&&\le C\int_\Omega \big[(|{n}^2-1)^2+({n}\cdot{m})^2\big].
\end{eqnarray}
Similarly, we have
\begin{eqnarray}\label{m-equation}
&&\frac{d}{dt}\int_\Omega (|{m}|^2-1)^2+\int_\Omega |\nabla({|m}|^2-1)|^2\nonumber\\
&&=4\int_\Omega |\nabla {m}|^2 (|{m}|^2-1)^2
+2\int_\Omega \langle\nabla{m}, \nabla{n}\rangle ({m}\cdot{n})(|{m}|^2-1)\nonumber\\
&&\le C\int_\Omega \big[(|{m}^2-1)^2+({n}\cdot{m})^2\big].
\end{eqnarray}
Now, by adding \eqref{nm-equation}, \eqref{n-equation}, and \eqref{m-equation} together, we obtain that
\begin{eqnarray}\label{n-m-nm-equation}
&&\frac{d}{dt}\int_\Omega \big[(|{n}|^2-1)^2+(|{m}|^2-1)^2+({n}\cdot{m})^2\big]\nonumber\\
&&+\int_\Omega \big(|\nabla(|{n}|^2-1)|^2+|\nabla(|{n}|^2-1)|^2+|\nabla({n}\cdot{m})|^2\big)\nonumber\\
&&\le C\int_\Omega \big[(|{n}|^2-1)^2+(|{m}|^2-1)^2+({n}\cdot{m})^2\big].
\end{eqnarray}
Since 
\begin{eqnarray*}
&&\int_\Omega \big[(|{n}|^2-1)^2+(|{m}|^2-1)^2+({n}\cdot{m})^2\big]\Big|_{t=0}\\
&&=\int_\Omega \big[(|{n}_0|^2-1)^2+(|{m}_0|^2-1)^2+({n}_0\cdot{m}_0)^2\big]=0,
\end{eqnarray*}
it follows from \eqref{n-m-nm-equation} and Gronwall's inequality that
$$
\int_\Omega \big[(|{n}|^2-1)^2+(|{m}|^2-1)^2+({n}\cdot{m})^2\big]=0, \ \forall t\in (0, T).
$$
This yields that $|{n}|=|{m}|=1$, and ${n}\cdot{m}=0$
in $\Omega\times (0,T)$. \qed

\medskip
Now we are ready to show the existence of local smooth solutions. 

\begin{theorem}\label{local} For any $\alpha\in (0,1)$, if $u_0\in C_0^{2+\alpha}(\Omega,\mathbb R^2)$ with ${\rm{div}}u_0=0$ and
${n}_0, {m}_0\in C^{2+\alpha}(\overline\Omega,\mathbb S^2)$ satisfies ${n}_0\cdot{m}_0=0$ in $\Omega$, then 
there is a $T>0$ depending on $\|u_0\|_{C^{2,\alpha}(\Omega)}, \|{n}_0\|_{C^{2+\alpha}(\Omega)}, \|{m}_0\|_{C^{2,\alpha}(\Omega)}$
such that there is a unique solution $(u, {n}, {m})\in C^{2+\alpha,1+\frac{\alpha}2}(\overline\Omega\times [0, T), \mathbb R^2\times \mathbb S^2\times \mathbb S^2)$ to
the initial value problem \eqref{biaxial-nlcf0} and \eqref{IBC}.
\end{theorem}
\pf For $T>0$ and $K>0$, set $Q_T=\Omega\times [0, T]$ and
\begin{eqnarray*}
X(K,T)&=&\Big\{(v, f, g)\in C^{2+\alpha,1+\frac{\alpha}2}(\overline{Q}_T, \mathbb R^2\times\mathbb R^3\times\mathbb R^3):
{\rm{div}} v=0 \ {\rm{in}}\ Q_T, \\
&&\ \ (v, f, g)=(u_0, {n}_0, {m}_0) \ {\rm{on}}\ \partial_pQ_T,
\big\|(v-u_0, f-{n}_0, g-{m}_0)\big\|_{C^{2+\alpha, 1+\frac{\alpha}2}(Q_T)}\le K\Big\}.
\end{eqnarray*}
Equip $X(K,T)$ with the norm
$$\big\|(v, f, g)\big\|_{X(K,T)}
:=\big\|(v, f, g)\big\|_{C^{2+\alpha, 1+\frac{\alpha}2}(Q_T)}, \ \forall (v, f, g)\in X(K,T),
$$
so that $(X(K,T), \|\cdot\|_{X(K,T)})$ is a Banach space. Define the solution operator
$$\mathcal{L}: X(K,T)\mapsto C^{2+\alpha, 1+\frac{\alpha}2}(Q_T)(\overline{Q}_T, \mathbb R^2\times\mathbb R^3\times\mathbb R^3)$$
as follows. For $(v, f, g)\in X(K,T)$, let $(u, n, m)=\mathcal{L}(v, f, g)$ be the unique solution to the following system:
\begin{equation}
\begin{cases}
u_t-\Delta u+\nabla P=-v\cdot\nabla v-\nabla\cdot(\nabla n\odot\nabla n+\nabla m\odot\nabla m), 
& {\rm{in}}\ Q_T,\\
\nabla\cdot u=0,  & {\rm{in}}\ Q_T,\\
n_t-\Delta n=-v\cdot\nabla f+|\nabla f|^2 f+\langle\nabla f, \nabla g\rangle f, & {\rm{in}}\ Q_T,\\
m_t-\Delta m=-v\cdot\nabla g+|\nabla g|^2 g+\langle\nabla g, \nabla f\rangle g, & {\rm{in}}\ Q_T,\\
(u, n, m)=(u_0, {n}_0, {m}_0), & {\rm{on}}\ \partial_pQ_T.
\end{cases}
\end{equation}
Now we want to show that $\mathcal{L}:X(K,T)\mapsto X(K,T)$ is a contraction map, 
provided $T>0$ is sufficiently small and $K>0$ is sufficiently large.

\smallskip
\noindent{\bf Claim 1}. {\it There exist $T>0$ and $K>0$ such that $\mathcal{L}:X(K,T)\mapsto X(K,T)$.} 

Set $$C_0\equiv 1+\|u_0\|_{C^{2+\alpha}(\Omega)}+\|{n}_0\|_{C^{2+\alpha}(\Omega)}+\|{m}_0\|_{C^{2+\alpha}(\Omega)}.$$
Assume $K\ge C_0^8$. 
It follows from the Schauder theory of parabolic equation that
\begin{equation}\label{holder0}
\begin{cases}
\|n-{n}_0\|_{C^{2+\alpha, 1+\frac{\alpha}2}(Q_T)}
\leq C\big[1+\|v\cdot\nabla f\|_{C^\alpha(Q_T)}
+\||\nabla f|^2 f\|_{C^\alpha(Q_T)}+\|\langle\nabla f, \nabla g\rangle g\|_{C^\alpha(Q_T)}\big],\\
\|m-{m}_0\|_{C^{2+\alpha, 1+\frac{\alpha}2}(Q_T))}
\leq C\big[1+\|v\cdot\nabla g\|_{C^\alpha(Q_T)}
+\||\nabla g|^2 g\|_{C^\alpha(Q_T)}+\|\langle\nabla g, \nabla f\rangle f\|_{C^\alpha(Q_T)}\big].
\end{cases}
\end{equation}
As in \cite{LLW10} Lemma 3.1, we can bound
\begin{eqnarray}\label{holder1}
&&\|v\cdot\nabla f\|_{C^\alpha(Q_T)}\nonumber\\
&&\le \|v\cdot\nabla f-u_0\cdot\nabla{n}_0\|_{C^\alpha(Q_T)}+
\|u_0\cdot\nabla{n}_0\|_{C^\alpha(\Omega)}\nonumber\\
&&\le C_0^2+\|v-u_0\|_{C^0(Q_T)}\|\nabla f\|_{C^\alpha(Q_T)}
+\|v-u_0\|_{C^\alpha(Q_T)}\|\nabla f\|_{C^0(Q_T)}\nonumber\\
&&+\|u_0\|_{C^0(Q_T)}\|\nabla (f-{n}_0)\|_{C^\alpha(Q_T)}
+\|u_0\|_{C^\alpha(Q_T)}\|\nabla (f-{n}_0)\|_{C^0(Q_T)}.
\end{eqnarray}
Since $v-u_0=f-{n}_0=0$ at $t=0$, it follows that
$$\|v-u_0\|_{C^0(Q_T)}+\|\nabla (f-{n}_0)\|_{C^0(Q_T)}
\le 2KT,$$
which, combined with the interpolation inequality, yields
that for any $0<\delta<1$, 
$$
\begin{cases} \|v-u_0\|_{C^\alpha(Q_T)}\leq C\big[\delta^{-1}\|v-u_0\|_{C^0(Q_T)}
+\delta \|v-u_0\|_{C^{2,1}_\alpha(Q_T)}\big]\leq C(\delta^{-1} T+\delta)K,\\
\|\nabla(f-{n}_0)\|_{C^\alpha(Q_T)}\leq C\big[\delta^{-1}\|\nabla(f-{n}_0)\|_{C^0(Q_T)}
+\delta \|f-{n}_0\|_{C^{2,1}_\alpha(Q_T)}\big]\leq  C(\delta^{-1} T+\delta)K.
\end{cases}
$$
Putting these estimates into \eqref{holder1}, we obtain
\begin{equation}\label{holder2}
 \|v\cdot\nabla f\|_{C^\alpha(Q_T)}\leq CK^2(T+\delta+\delta^{-1}T)+C_0^2.
\end{equation}
We can bound
\begin{eqnarray}\label{holder3}
&&\||\nabla f|^2 f\|_{C^\alpha(Q_T)}\nonumber\\
&&\le \||\nabla f|^2 f-|\nabla{n}_0|^2{n}_0\|_{C^\alpha(Q_T)}
+\||\nabla {n}_0|^2 {n}_0\|_{C^\alpha(\Omega)}\nonumber\\
&&\leq  C\|f\|_{C^{2+\alpha, 1+\frac{\alpha}2}(Q_T)}^2\|f-{n}_0\|_{C^\alpha(Q_T)}
+C(\|f\|_{C^{2+\alpha, 1+\frac{\alpha}2}(Q_T)}\\
&&+\|{n}_0\|_{C^{2+\alpha}(\Omega)})
\|\nabla(f-{n}_0)\|_{C^\alpha(Q_T)}+C_0^2\nonumber\\
&&\leq CK^2 \big(\|f-{n}_0\|_{C^\alpha(Q_T)}+C\|\nabla(f-{n}_0)\|_{C^\alpha(Q_T)}\big)+C_0^2
\nonumber\\
&&\leq CK^2(\delta +T\delta^{-1})K+C_0^2.
\end{eqnarray}
Similarly, we can also bound
\begin{eqnarray}\label{holder4}
&&\|\langle\nabla f,\nabla g\rangle f\|_{C^\alpha(Q_T)}\nonumber\\
&&\le \|\langle\nabla f,\nabla g\rangle f-\langle\nabla{n}_0, \nabla{m}_0\rangle{n}_0\|_{C^\alpha(Q_T)}
+\|\langle\nabla {n}_0,\nabla {m}_0\rangle {n}_0\|_{C^\alpha(\Omega)}\nonumber\\
&&\leq C_0^2+CK^2\big[\|f-{n}_0\|_{C^\alpha(Q_T)}+\|\nabla(f-{n}_0)\|_{C^\alpha(Q_T)}\nonumber\\
&&\ \ \ +C\|g-{m}_0\|_{C^\alpha(Q_T)}+\|\nabla(g-{m}_0)\|_{C^\alpha(Q_T)}\big]\nonumber\\
&&\leq CK^2(\delta +T\delta^{-1})K+C_0^2.
\end{eqnarray}
Substituting \eqref{holder2}, \eqref{holder3}, and \eqref{holder4} into \eqref{holder0}$_1$, we arrive at
\begin{equation}\label{holder5}
\|n-{n}_0\|_{C^{2+\alpha, 1+\frac{\alpha}2}(Q_T)}
\leq CK^2(\delta +T\delta^{-1})K+CC_0^2\le {K}^\frac13,
\end{equation}
provided $T=\delta^2$ and $\delta\approx C^{-1}K^{-\frac53}$.

By the same arguments as in \eqref{holder5}, we can also obtain that
\begin{equation}\label{holder6}
\|m-{m}_0\|_{C^{2+\alpha, 1+\frac{\alpha}2}(Q_T)}
\leq CK^2(\delta +T\delta^{-1})K+CC_0^2\le {K}^\frac13,
\end{equation}
provided $T=\delta^2$ and $\delta\approx C^{-1}K^{-\frac53}$.

For $u$, we apply the Schauder theory of linear non-stationary Stokes equations
to obtain 
\begin{equation}\label{holder7}
\|u-u_0\|_{C^{2+\alpha, 1+\frac{\alpha}2}(Q_T)}
\leq C+C\|v\cdot \nabla v\|_{C^\alpha(Q_T)}
+C\|\nabla\cdot(\nabla n\odot\nabla n+\nabla m\odot\nabla m)\|_{C^\alpha(Q_T)}.
\end{equation}
Similar to \eqref{holder2}, it is not hard to see that
\begin{eqnarray}\label{holder8}
\|v\cdot \nabla v\|_{C^\alpha(Q_T)}
&\le& \|v\cdot \nabla v-u_0\cdot\nabla u_0\|_{C^\alpha(Q_T)}+\|u_0\cdot \nabla u_0\|_{C^\alpha(\Omega)}
\nonumber\\
&\le& C_0^2+K\big[\|v-u_0\|_{C^\alpha(Q_T)}+\|\nabla(v-u_0)\|_{C^\alpha(Q_T)}\big]\nonumber\\
&\le& C_0^2+CK^2(T+\delta+\delta^{-1}T).
\end{eqnarray}
%provided $K=8C_0$, $\delta\approx K^{-2}$, and $T=\delta^2$.
For the second term in the right hand side of \eqref{holder7}, it follows from \eqref{holder5} and \eqref{holder6}
that
\begin{eqnarray}\label{holder9}
&&\|\nabla\cdot(\nabla n\odot\nabla n+\nabla m\odot\nabla m)\|_{C^\alpha(Q_T)}\nonumber\\
&&\le C\big(\|{n}_0|\|_{C^{2+\alpha}\alpha(\Omega)}^2
+ \|{m}_0\|_{C^{2+\alpha}(\Omega)}^2\big)\nonumber\\
&&\quad+C\big(\|n-{n}_0\|_{C^{2+\alpha, 1+\frac{\alpha}2}(Q_T)}^2
+\|m-{m}_0\|_{C^{2+\alpha, 1+\frac{\alpha}2}(Q_T)}^2\big)\nonumber\\
&&\leq CC_0^2+CK^\frac23\le \frac{K}3.
\end{eqnarray}
Putting \eqref{holder8} and \eqref{holder9} into \eqref{holder7}, we obtain
\begin{equation} \label{holder11}
\|u-u_0\|_{C^{2+\alpha, 1+\frac{\alpha}2}(Q_T)}+\|{n}-{n}_0\|_{C^{2+\alpha, 1+\frac{\alpha}2}(Q_T)}
+\|m-{m}_0\|_{C^{2+\alpha, 1+\frac{\alpha}2}(Q_T)}\le K.
\end{equation}
Hence $\mathcal{L}:X(K,T)\to X(K,T)$ is a bounded operator.

\medskip
\noindent{\bf Claim 2}. {\it There exist $K>0$ and $T>0$ such that $\mathcal{L}:X(K,T)\to X(K,T)$ is
a contractive map.}

For $i=1,2$ and $(v_i, f_i, g_i)\in X(K,T)$, let $(u_i, n_i, m_i)
=\mathcal{L}(v_i, f_i, g_i)$. 
Set $$(u, n, m)=(u_1-u_2, n_1-n_2, m_1-m_2).$$
Then we have that
\begin{equation}
\begin{cases} 
u_t-\Delta u+\nabla P=F\\
:=-u_1\cdot\nabla u-u\cdot\nabla u_2-\nabla\cdot\big(\nabla n_1\odot\nabla n
+\nabla n\odot\nabla n_2+\nabla m_1\odot\nabla m+\nabla m\odot\nabla m_2\big),\\
\nabla\cdot u=0,\\
n_t-\Delta n=G:=-u_1\cdot\nabla n-u\cdot\nabla n_2+|\nabla n_1|^2 n-\langle \nabla (n_1+n_2),\nabla n\rangle n_2\\
+\langle\nabla n_1,\nabla m_1\rangle m+(\langle\nabla n,\nabla m_1\rangle+\langle \nabla n_2,\nabla m\rangle) m_2,\\
m_t-\Delta m=H:=-u_1\cdot\nabla m-u\cdot\nabla m_2+|\nabla m_1|^2 m +\langle \nabla (m_1+m_2),\nabla m\rangle m_2\\
+\langle\nabla n_1,\nabla m_1\rangle n+(\langle\nabla n,\nabla m_1\rangle+\langle\nabla n_2, \nabla m\rangle) n_2.
\end{cases}
\end{equation}
Applying the Schauder theory, we obtain that
\begin{align}\label{n}
\|n\|_{C^{2+\alpha,1+\frac{\alpha}{2}}(Q_{T})}
\leq C\|{G}\|_{C^{\alpha,\frac{\alpha}{2}}(Q_{T})}
&\leq C(K^{2}+K) \|(u,n, m)\|_{C^{1+\alpha,\frac{\alpha}{2}}(Q_{T})}\nonumber\\
&\leq CK^2(\delta+T\delta^{-1}) \|(u,n, m)\|_{C^{2+\alpha,1+\frac{\alpha}{2}}(Q_{T})},
\end{align}
%--------------------(dW32)--------------------------------
\begin{align}\label{m}
\|m\|_{C^{2+\alpha,1+\frac{\alpha}{2}}(Q_{T})}
\leq C\|{H}\|_{C^{\alpha,\frac{\alpha}{2}}(Q_{T})}
&\leq C(K^{2}+K) \|(u,n, m)\|_{C^{1+\alpha,\frac{\alpha}{2}}(Q_{T})}\nonumber\\
&\leq CK^{2}(\delta+T\delta^{-1}) \|(u,n, m)\|_{C^{2+\alpha,1+\frac{\alpha}{2}}(Q_{T})},
\end{align}
and
\begin{align}\label{u}
&\|{u}\|_{C^{2+\alpha,1+\frac{\alpha}{2}}(Q_{T})}
\leq C\|{F}\|_{C^{\alpha,\frac{\alpha}{2}}(Q_{T})}\nonumber\\
&\leq CK\|u\|_{C^{1+\alpha, \frac{\alpha}2}(Q_T)}+CK \|(n, m)\|_{C^{2+\alpha, 1+\frac{\alpha}2}(Q_T)}\nonumber\\
&\leq CK(\delta+T\delta^{-1}) \|(u,n, m)\|_{C^{2+\alpha,1+\frac{\alpha}{2}}(Q_{T})}.
\end{align}
%--------------------(uW21)--------------------------------
It follows that
\begin{align*}
&\|\mathcal{L}({u}_{1}, n_{1},m_1)-\mathcal{L}(u_2, n_{2}, m_2)\|_{{X}(K,T)}
=\|({u},n, m)\|_{C^{2+\alpha,1+\frac{\alpha}{2}}(Q_{T})}\\
&\leq CK^2 (\delta+T{\delta}^{-1})\|({u},n, m)\|_{C^{2+\alpha,1+\frac{\alpha}{2}}(Q_{T})}\nonumber\\
&\le \frac{1}{2} \|(u_1, n_1, m_1)-(u_2, n_2, m_2)\|_{{X}(K,T)},
\end{align*}
provided $\delta$ and $T$ are sufficiently small. Thus 
$\mathcal{L}:{X}(K, T)\mapsto {X}(K,T)$ is a contractive map.

It follows from Claim 2 and the fixed point theorem that there exists $T>0$ depending on $C^{2,\alpha}$-norm of 
$(u_0, n_0, m_0)$ and a unique smooth solution $(u, n, m)\in C^{2+\alpha, 1+\frac{\alpha}2}(Q_T, \mathbb R^2\times 
\mathbb R^3\times \mathbb R^3)$ to the system \eqref{biaxial-nlcf0}$_{1,2,3,4}$, with the initial-boundary data
$(u_0, n_0, m_0)$. Since $|n_0|=|m_0|=1$ and $n_0\cdot m_0=0$ holds in $\Omega$, it follows from Lemma
\ref{constraint} that the solution $(u, n, m)$ satisfies the constraint $|n|=|m|=1$ and $n\cdot m=0$ in $Q_T$. This
completes the proof of Theorem \ref{local}.
\qed

\section{A priori estimate under the smallness condition}
\setcounter{equation}{0}

In this section, we will show that both interior and boundary gradient estimates of weak solutions
of \eqref{biaxial-nlcf0} under the smallness condition on $(u, n, m)$.  

\begin{lemma}\label{interior-est} For any $M>0$, there exist $\epsilon_0>0$, $\theta_0\in (0,\frac12)$ such that if $z_0=(x_0,t_0)\in\mathbb R^2\times (0,\infty)$  and $0<r_0\le\sqrt{t_0}$, $u\in  L^2_tH^1_x(P_{r_0}(z_0),\mathbb R^2)\big)$,
$n, m \in L^2_t \big([t_0-r_0^2, t_0], H^1(B_{r_0}(x_0),\mathbb S^2)\big)$ with $n\cdot m=0$ in $P_{r_0}(z_0)$,
and $P\in W^{1,0}_{\frac43}(P_{r_0}(z_0))$ is a weak solution of \eqref{biaxial-nlcf0} in $P_{r_0}(z_0)$, satisfying\begin{equation}\label{small_cond}
r_0|(u)_{z_0, r_0}|\le M;
\ \ \ \  r_0 Y(z_0, r_0;  u, n, m, p) \le \epsilon_0,
\end{equation}
then
\begin{eqnarray}\label{energy_decay}
Y(z_0,\theta_0 r_0;  u, n, m ,P)\le \theta_0^\frac13 Y(z_0, r_0;  u, n, m , P).
\end{eqnarray}
\end{lemma}
Here and henceforth, we use the following notation:
\begin{equation*}
\begin{cases}\displaystyle
(u)_{z_0, r}=\frac{1}{|P_{r}(z_0)|}\int_{P_{r}(z_0)} u; \ \ \  [P]_{x_0, r}=\frac{1}{|B_r(x_0)|}\int_{B_r(x_0)} P,\\
\displaystyle Y(z_0, r;  u, n, m, P)=Y^1(z_0, r;  u, n, m)+Y^2(z_0, r;  P),\\
\displaystyle Y^1(z_0, r;  u, n, m)= \Big[\frac{1}{|P_{r}(z_0)|}\int_{P_{r}(z_0)} (|u-(u)_{z_0, r_0}|^4
+ |\nabla n|^4+|\nabla m|^4)\Big]^\frac14,\\
\displaystyle Y^2(z_0, r;  P)=r\Big(\frac{1}{|P_r(z_0)|}\int_{P_{r}(z_0)}|P-[P]_{x_0, r}|^\frac43\Big)^\frac34.
\end{cases}
\end{equation*}

\pf First observe that by the interior $L^2$-theory of parabolic equations, it follows
from \eqref{small_cond} and \eqref{biaxial-nlcf0}$_{3,4}$ that $n, m\in W^{2,1}_2(P_r(z_0))$.
%Next we can apply the proof of Lemma 2.2 to conclude that $(u, n, m, P)$ satisfies the local energy
%inequality \eqref{LED0} in $P_r(z_0)$. 
By translation and dilation, it suffices to show \eqref{energy_decay}
for $z_0=(0,1)$ and $r_0=1$. We argue by contradiction. Suppose that it were false. Then there exists
$M>0$ such that for any
$\theta\in (0,\frac12)$ there exist
a sequence of weak solutions $(u^i, n^i, m^i, P^i)$ of \eqref{biaxial-nlcf0} in $P_1(0,1)$ such that
$$
|(u^i)_{(0,1), 1}|\le M; \ \ \ Y((0,1), 1; u^i, n^i, m^i, P^i)=\epsilon_i\rightarrow 0,
$$
but
\begin{eqnarray*}
Y((0,1), \theta; u^i, n^i, m^i, P^i)>
\theta^\frac13 Y((0,1), 1; u^i, n^i, m^i, P^i).
\end{eqnarray*}
Applying the interior $L^2$-estimate to \eqref{biaxial-nlcf0}$_{3,4}$, we obtain that
\begin{eqnarray*}
&&\||n^i_t|+|\nabla^2 n^i|\|_{L^2(P_{\frac34}(0,1))}+\||m^i_t|+|\nabla^2 m^i|\|_{L^2(P_{\frac34}(0,1))}\\
&&\le C\Big(\int_{P_{1}(0,1)} \big(|u^i-(u^i)_{(0,1),1}|^4+(M+1)(|\nabla n^i|^4+|\nabla m^i|^4)\big)\Big)^\frac14\le C\epsilon_i.
\end{eqnarray*}
Apply the $L^2$-estimate to \eqref{biaxial-nlcf0}$_{1}$, we have that
$u^i_t\in L^{2}([\frac34, 1], H^{-1}(B_\frac12)), \nabla u_i\in L^2(P_\frac12(0,1)),$ and
\begin{eqnarray*}
&&\|u^i_t\|_{L^{2}([\frac34, 1], H^{-1}(B_\frac12))}+\|\nabla u_i\|_{L^2(P_\frac12(0,1))}\\
&&\le C\Big(\int_{P_{1}(0,1)} \big(|u^i|^4+(M+1)(|\nabla n^i|^4+|\nabla m^i|^4)\big)\Big)^\frac14\le C\epsilon_i.
\end{eqnarray*}
Define the blowing up sequence in $P_1(0,1)$: 
$$\tilde{u}^i=\frac{u^i-(u^i)_{(0,1),1}}{\epsilon_i}, \ \tilde{n}^i=\frac{n^i-(n^i)_{(0,1),1}}{\epsilon_i},
\ \tilde{m}^i=\frac{m^i-(m^i)_{(0,1), 1}}{\epsilon_i}, \ \tilde{P}^i=\frac{P^i-[P^i]_{(0,1), 1}}{\epsilon_i}.$$
%where $\displaystyle(f)_2=\frac{1}{|P_2(0,4)|} \int_{P_2(0,4)} f$ denotes the average of $f$ over $P_2(0,4)$. 
Then it holds that
\begin{equation}\label{sequence-bound}
Y((0,1), 1; \tilde{u}^i, \tilde{n}^i, \tilde{m}^i, \tilde{P}^i)=1
\end{equation}
and
$$
\|\tilde{u}^i_t\|_{L^{2}([\frac34, 1], H^{-1}(B_\frac12))}^2+\int_{P_\frac12(0,1)}\big(|\nabla\tilde{u}^i|^2+|\tilde{n}^i_t|^2+|\tilde{m}^i_t|^2+|\nabla^2\tilde{n}^i|^2+|\nabla^2\tilde{m}^i|^2\big)\le C,
$$
but
\begin{equation}\label{sequence-small}
Y((0,1), \theta; \tilde{u}^i, \tilde{n}^i, \tilde{m}^i, \tilde{P}^i)>\theta^\frac13.
\end{equation}

By Aubin-Lions' compactness Lemma, we may assume that there exist $\tilde{u}\in L^2([\frac34,1], H^1(B_\frac12))$, with
$\tilde{u}_t\in L^{2}([\frac34,1], H^{-1}(B_\frac12))$,  $\tilde{n}, \tilde{m}\in W^{2,1}_2(P_\frac12(0,1), \mathbb R^3)$,
$\tilde{P}\in L^\frac43(P_\frac12(0,1))$, and $b\in\mathbb R^2$ with $|b|\le M$ such that, 
after passing to a subsequence, it holds that
\begin{equation*}
\begin{cases}
\tilde{u}^i\rightharpoonup \tilde{u}\  {\rm{in}}\  L^2([\frac34,1], H^1(B_\frac12)),  \  
\tilde{u}^i_t\rightharpoonup\tilde{u}_t\ {\rm{in}}\  L^{2}([\frac34,1], H^{-1}(B_\frac12)), \\
\tilde{P}^i \rightharpoonup\tilde{P}\ {\rm{in}}\  L^\frac43(P_\frac12(0,1)),\\
\tilde{n}^i \rightharpoonup\tilde{n},\  \tilde{m}^i\rightharpoonup \tilde{m} \ \ {\rm{in}}\ W^{2,1}_2(P_\frac12(0,1)),\\
(u^i)_{(0,1), 1}\rightarrow b,
\end{cases}
\end{equation*}
and
\begin{equation}\label{strong-l4}
\tilde{u}^i\rightarrow \tilde{u}, \ \nabla\tilde{n}^i\rightarrow \nabla\tilde{n},\ 
\nabla\tilde{m}^i\rightarrow \nabla\tilde{m} \ \ {\rm{in}}\ L^4(P_\frac12(0,1)).
\end{equation}
By the lower semicontinuity, we have that
\begin{equation}\label{l4-bound}
Y((0,1), \frac12; \tilde{u}, \tilde{n},\tilde{m},\tilde{P}) \le 1.
\end{equation}
Note that $(\tilde{u}^i, \tilde{n}^i, \tilde{m}^i, \tilde{P}^i)$ solves, in $P_1(0,1)$,
\begin{equation}\label{biaxial-nlcf1}
\begin{cases}
\tilde{u}^i_t+(\epsilon_i\tilde{u}^i+(u^i)_{(0,1), 1})\cdot\nabla \tilde{u}^i+\nabla \tilde{P}^i-\Delta \tilde{u}^i\\
=-\epsilon_i{\rm{div}}\big(\nabla \tilde{n}^i\odot \nabla \tilde{n}^i+\nabla \tilde{m}^i\odot \nabla \tilde{m}^i\big),\\
{\rm{div}} (\tilde{u}^i)=0,\\
\tilde{n}^i_t+(\epsilon_i\tilde{u}^i+(u^i)_{(0,1), 1})\cdot\nabla\tilde{n}^i-\Delta\tilde{n}^i
=\epsilon_i|\nabla\tilde{n}^i|^2{n^i}+\epsilon_i(\nabla\tilde{n}^i\cdot\nabla\tilde{m}^i) {m^i},\\
\tilde{m}^i_t+(\epsilon_i\tilde{u}^i+(u^i)_{(0,1), 1})\cdot\nabla\tilde{m}^i-\Delta\tilde{m}^i
=\epsilon_i|\nabla\tilde{m}^i|^2{m^i}+\epsilon_i(\nabla\tilde{n}^i\cdot\nabla\tilde{m}^i) {n^i}.
\end{cases}
\end{equation}
Passing to the limit in \eqref{biaxial-nlcf1}, we obtain that $(\tilde{u}, \tilde{n}, \tilde{m},\tilde{P})$ solves
in $P_1(0,1)$: 
\begin{equation}\label{limit-eqn}
\begin{cases}
\tilde{u}_t-\Delta \tilde{u}+b\cdot\nabla\tilde{u}+\nabla\tilde{P}=0,\ {\rm{div}} (\tilde{u})=0,\\
\tilde{n}_t-\Delta\tilde{n}+b\cdot\nabla\tilde{n} =0,\\
\tilde{m}_t-\Delta\tilde{m}+b\cdot\nabla\tilde{m}=0.
\end{cases}
\end{equation}
By the regularity theory of the linear Stokes equation and the heat equation,  
there exists $C=C(M)>0$ that for any $\theta\in (0,\frac12)$,
\begin{equation}\label{limit-reg}
Y((0,1), \theta; \tilde{u}, \tilde{n},  \tilde{m}, \tilde{P})
\le C\theta.
\end{equation}
This, combined with \eqref{strong-l4}, implies that 
\begin{equation}\label{sequence-reg}
Y^1((0,1),\theta;  \tilde{u}^i, \tilde{n}^i, \tilde{m}^i)
\le C\theta +\theta^{-1} o(i^{-1}),
\end{equation}
provided we choose  a sufficiently large $i$. 

To estimate $Y^2((0,1), \theta; \tilde{P}^i)$, first recall that $\tilde{P}^i$ solves the Poisson equation:
\begin{equation}\label{eqn-P}
-\Delta\tilde{P}^i=\epsilon_i{\rm{div}}^2\big(\tilde{u}^i\otimes\tilde{u}^i+\nabla \tilde{n}^i\odot \nabla \tilde{n}^i 
+\nabla \tilde{m}^i\odot \nabla \tilde{m}^i \big)\ \ {\rm{in}}\ \ B_1(0).
\end{equation}
Applying the Poincar\'e inequality and $W^{1,\frac43}$-estimate to \eqref{eqn-P} and the estimate \eqref{sequence-reg}, we have that for $i$ sufficiently large,
\begin{eqnarray}\label{P-est}
&&Y^2((0,1), \theta; \tilde{P}^i)=\theta^{-2} \Big(\int_{P_{\theta}(0,1)}\big|\tilde{P}^i-[\tilde{P}^i]_{(0,1),\theta}\big|^\frac43\Big)^\frac34\le \theta^{-1} \Big(\int_{P_{\theta}(0,1)}\big|\nabla\tilde{P}^i\big|^\frac43\Big)^\frac34\nonumber\\
&&\le C\theta^{-1}\epsilon_i\Big[\big(\int_{P_\frac34(0,1)} \big(|\tilde{u}^i|^4+|\nabla\tilde{n}^i|^4+|\nabla\tilde{m}^i|^4\big)\big)\frac14+\big(\int_{P_\frac34(0,1)}|\nabla \tilde{u}^i|^2+|\nabla^2 \tilde{n}^i|^2+|\nabla^2\tilde{m}^i|^2\big)^\frac12\Big]\nonumber\\
&&\qquad+C\theta^\frac12\Big( \int_{P_1(0,1)} \big|\tilde{P}^i-[\tilde{P}]^i_{(0,1), 1}\big|^\frac43\Big)^\frac34 \nonumber\\
&&\le C(\theta^{-1}\epsilon_i+\theta^\frac12).
\end{eqnarray}
Putting together \eqref{sequence-reg} with \eqref{P-est}, we obtain that
\begin{equation}\label{sequence-reg1}
Y((0,1), \theta; \tilde{u}^i, \tilde{n}^i, \tilde{m}^i, \tilde{P}^i)
\le C(\theta^{-1}\epsilon_i+\theta^\frac12+\theta^{-1}o(i^{-1}))<\theta^\frac13,
\end{equation}
provided $\theta$ is chosen sufficiently small and $i$ is chosen sufficiently large. 
This contradicts to \eqref{sequence-small}. Hence
the proof is complete.\qed

\bigskip
Iterating the estimates in Lemma \ref{interior-est}, we can obtain the interior regularity of weak solutions to \eqref{biaxial-nlcf0}. More precisely, we have

\begin{proposition} \label{interior-reg} Under the same assumptions as in Lemma \ref{interior-est}, $(u, n, m)\in C^\infty(P_{\frac{r_0}2}(z_0))$ and for any $l\ge 1$, it holds that
\begin{eqnarray}\label{interior-l-est} 
\big\|(u, n, m)\big\|_{C^l(P_{\frac{r_0}2}(z_0))}\le C(\epsilon_0, r_0, l).
\end{eqnarray}
\end{proposition} 
\pf First we claim that by iterating \eqref{energy_decay} finitely many times, it holds that any $z\in P_{\frac{r_0}2}(z_0)$ and for all $k\ge 1$, 
\begin{eqnarray}\label{morrey_decay1}
Y(z, \theta_0^kr_0; u, n, m, P)\le \theta_0^{\frac{k}3}
Y(z, r_0; u, n, m, P).
\end{eqnarray}
For $k=1$, \eqref{morrey_decay1} follows from Lemma \ref{interior-est}.  For $k= 2$, one can show
\eqref{morrey_decay1} as follows.  By the triangle inequality, 
we have that
\begin{eqnarray*}
(\theta_0r_0)|(u)_{z,\theta_0 r_0}|&\le& (\theta_0 r_0)|u_{z,\theta_0r_0}-u_{z,r_0}|+\theta_0 r_0|u_{z, r_0}|\\
&\le&\theta_0^{-3} Y(z, r_0; u, n, m, P)+\theta_0 r_0|u_{z, r_0}|\\
&\le& \theta_0^{-3}\epsilon_0+\theta_0 M\le M.
\end{eqnarray*}
This, combine with the fact that $Y(z, \theta_0r_0; u, n, m, P)\le \theta_0^\frac13\epsilon_0$ and Lemma \ref{interior-est}, implies that 
$$Y(z,\theta_0^2 r_0; u, n, m, P)\le \theta_0^\frac13Y(z,\theta_0 r_0; u, n, m, P)
\le \theta_0^\frac23 Y(z, r_0; u, n, m, P).
$$
For $k\ge 3$, one can prove \eqref{morrey_decay1} by an induction on $k$. We omit it here.
From \eqref{morrey_decay1}, there exists $\alpha_0\in (\frac12, 1)$ such that 
for any $0<r\le r_0$ and $z\in P_{\frac{r_0}2}(z_0)$, it holds
\begin{equation}\label{morrey_decay2}
Y(z, r; u, n, m, P)\le \big(\frac{r}{r_0}\big)^{\alpha_0} 
Y(z, r_0; u, n, m, P)\le \epsilon_0 \big(\frac{r}{r_0}\big)^{\alpha_0}.
\end{equation}
On the other hand, it is not hard to verify that $n$ and $m$ satisfies the following local energy inequality: 
\begin{equation}\label{local-nm-energy}
\int_{P_{\frac{r}2}(z)} (|n_t|^2+|m_t|^2)
\le C\int_{P_{r}(z)}[|u-u_{z, r}|^4+(M+1) |\nabla n|^4+|\nabla m|^4]
\end{equation}
holds for all $P_r(z)\subset P_{r_0}(z_0)$. Combining \eqref{morrey_decay1} with \eqref{local-nm-energy},
we obtain that
$$
\begin{cases}
\displaystyle  r^{-4}\int_{P_r(z)} |u-u_{z,r}|^4\le C\big(\frac{r}{r_0}\big)^{4\alpha_0},\\
\displaystyle  r^{-2}\int_{P_r(z)}(r^2(|n_t|^2+|m_t|^2)+(|\nabla n|^2+|\nabla m|^2))\le C\big(\frac{r}{r_0}\big)^{2\alpha_0}, 
\end{cases}
$$
for all $P_r(z)\subset P_{r_0}(z_0)$. Thus by the Morrey decay Lemma and the Campanato Lemma we conclude that
$(u, n, m)\in C^{\alpha_0,\frac{\alpha_0}2}(P_{\frac{r_0}2}(z_0))$ and
$$
\big\|(u, n, m)\big\|_{C^{\alpha_0,\frac{\alpha_0}2}(P_{\frac{r_0}2}(z_0))}\le C(M, \epsilon_0, r_0).
$$
The H\"older regularity of $(\nabla n, \nabla m)$  can be obtained, 
similar to the heat flow of harmonic maps (see \cite{LW2}). Substituting H\"older continuity of
$u,\nabla n, \nabla m$ in to \eqref{biaxial-nlcf0}$_{1,2}$, we can obtain the H\"older continuity of
$\nabla u$. The higher order regularity of $(u, n, m)$ then follows from the Schauder theory
of both the Stokes equation and parabolic equations. Moreover, the estimate \eqref{interior-l-est} holds. 
This completes the proof. \qed

\bigskip
We also need a corresponding boundary $\epsilon_0$-regularity estimate. 
For $z_0=(x_0,t_0)\in \partial \Omega\times (0, T]$ and $0<r\le r_0\in \sqrt{t_0}$, set
$$B_{r}^+(x_0)=\Omega\cap B_{r}(x_0);\  P_r^+(z_0)=(\Omega\cap B_{r}(x_0))\times [t_0-r^2, t_0],$$
and let $(u)_{z_0^+, r}$ and $[P]_{x_0^+, r}$ denote the average of $u$ over $P_{r}^+(z_0)$
and the average of $P$ over $B_{r}^+(x_0)$ respectively.  Also, denote
by $Y^1(z_0^+, r; u, n, m)$,  $Y^2(z_0^+, r; P)$, and $Y(z_0^+, r; u, n, m, P)$ by
\begin{equation*}
\begin{cases}
\displaystyle Y(z_0^+, r;  u, n, m, P)=Y^1(z_0^+, r;  u, n, m)+Y^2(z_0^+, r;  P),\\
\displaystyle Y^1(z_0^+, r;  u, n, m)= \Big[\frac{1}{|P_{r}^+(z_0)|}\int_{P_{r}(z_0)} (|u-(u)_{z_0^+, r}|^4
+|\nabla n|^4+|\nabla m|^4)\Big]^\frac14,\\
\displaystyle Y^2(z_0, r;  P)=r\Big(\frac{1}{|P_r(z_0)|}\int_{P_{r}(z_0)}|P-[P]_{x_0^+, r}|^\frac43\Big)^\frac34.
\end{cases}
\end{equation*}

\begin{lemma} \label{bdry_reg} For $\alpha\in (0,1)$ and $T>0$, let $n_0,m_0\in H^1(\Omega, \mathbb S^2)\cap C^{2+\alpha}(\partial\Omega,\mathbb S^2)$ be such that $n_0\cdot m_0=0$ in $\overline\Omega$. For any
$M>0$, there exist $\epsilon_0>0$,
$\theta_0>0$, $C_0>0$, and $r_0>0$, depending on $\partial\Omega$ and $n_0, m_0$,
such that if $u\in L^\infty_tL^2_x\cap L^2_tH^1_x(Q_T, \mathbb R^2)$,
$n, m\in L^\infty_t H^1_x(Q_T, \mathbb S^2)$ with $n\cdot m=0$ in $Q_T$, $P\in W^{1,0}_{\frac43}(Q_T)$,
is a weak solution of \eqref{biaxial-nlcf0}, with $(u, n, m)=(0,n_0, m_0)$ on $\partial\Omega$, that satisfies, for $x_0\in\partial\Omega$ and $t_0\in (r_0^2, T]$,
\begin{equation}\label{bdry_small}
r_0|(u)_{z_0^+,r_0}|\le M; \ \ \ \  r_0Y(z_0^+, r_0; u, n, m , P)\le\epsilon_0,
\end{equation}
then
\begin{eqnarray}\label{bdry_decay}
Y(z_0^+, \theta_0r_0; u, n, m, P)
\le \theta_0^\frac13\max\Big\{Y(z_0^+, r_0, u, n, m, P), \ C_0 \|(n_0, m_0)\|_{C^1(\partial\Omega\cap B_{r_0}(x_0))} r_0\Big\}.
\end{eqnarray}
\end{lemma}
\pf  Since the case for general curved boundary $\partial\Omega$ can be handled by the technique of flattening of the boundary similar to \cite{SSS}, for simplicity we will only consider the case that $x_0=0\in\partial\Omega$
and $\Omega\cap B_{r_0}(x_0)=\mathbb R^2_+\cap B_{r_0}(0):=B_{r_0}^+$. From the assumption \eqref{bdry_small} and $\nabla P\in L^\frac43(Q_T)$, we can apply the boundary $L^p$-estimate of parabolic equations to conclude that
$u\in W^{2,1}_{\frac43} (P_r^+(z_0), n\in W^{2,1}_2(P_r^+(z_0)), m\in W^{2,1}_2(P_r^+(z_0))$ for any
$0<r<r_0$ so that $(u, n, m, P)$ satisfies the boundary local energy inequality \eqref{LED1} on $P_{r_0}^+(z_0)$.  

By scaling, it suffices to prove it for $r_0=1$ and $x_0=0\in \partial\Omega$, $t_0=1$. 
We would argue by contradiction. Suppose it were false. Then there exists $M>0$ such that for any
$\theta\in (0,\frac12)$,  there exists a sequence of weak solutions $(u^k, n^k, m^k, P^k)$ to \eqref{biaxial-nlcf0} in
$P_1^+(0,1)$, with $(u^k, n^k, m^k)=(0, n_0^k, m_0^k)$ on $T_1\times [0,1]$, where $T_1=\partial B_1^+\cap \partial \mathbb R^2_+$ and $(n_0^k, m_0^k)\in C^{2+\alpha}(T_1,\mathbb S^2)^{\otimes 2}$ with $n_0^k\cdot m_0^k=0$
on $T_1$, such that 
\begin{equation}\label{bdry_small1}
|(u^k)_{(0,1)^+, 1}|\le M; \ \ \ \ Y((0,1)^+, 1; u^k, n^k, m^k, P^k)=\epsilon_k\rightarrow 0,
\end{equation}
but
\begin{eqnarray}\label{bdry_decay1}
Y((0,1)^+, \theta; u^k, n^k, m^k, P^k)
> \theta^\frac13\max\Big\{Y((0,1)^+, 1; u^k, n^k, m^k, P^k), \ k \|(n_0^k, m_0^k)\|_{C^1(T_1)}\Big\}.
\end{eqnarray}
By the $L^2$-theory of parabolic equations to \eqref{biaxial-nlcf0}$_{3,4}$,
$(n^k_t,\nabla^2 n^k, m^k_t, \nabla^2 m^k)\in L^2(P_{\frac34}^+(0,1))$, and
\begin{eqnarray}\label{n-m-t}
&&\int_{P_{\frac34}^+(0,1)} \big(|n^k_t|^2+|\nabla^2 n^k|^2+|m^k_t|^2+|\nabla^2 m^k|^2\big)\nonumber\\
&&\le C\int_{P_1^+(0,1)} \big(|u^k|^4+(M+1)(|\nabla n^k|^4+|\nabla m^k|^4))+C\|(n_0^k, m_0^k)\|_{C^1(T_1)}^4\nonumber\\
&&\le C(M, \theta). 
\end{eqnarray}
Define the blowing up sequence 
$$\big(\tilde{u}^k, \tilde{n}^k, \tilde{m}^k, \tilde{P}^k\big)
=\big(\frac{{u}^k-(u^k)_{(0,1)^+,1}}{\epsilon_k}, \frac{n^k-n_0^k(0)}{\epsilon_k}, \frac{{m}^k-m_0^k(0)}{\epsilon_k}, \frac{{P}^k-[P^k]_{(0,1)^+,1}}{\epsilon_k}\big):P_1^+\to \mathbb R^8.$$ 
From \eqref{bdry_decay1}, we have
$$
\|(n_0^k, m_0^k)\|_{C^1(T_1)}\le (\theta^{\frac13} k)^{-1}\epsilon_k\rightarrow 0.
$$
Then we would have 
\begin{equation}\label{bdry_small2}
\Big(\int_{P_{1}^+(0,1)} (|\tilde{u}^k|^4+|\nabla \tilde{n}^k|^4+|\nabla \tilde{m}^k|^4) \Big)^\frac14
+\Big(\int_{P_{1}^+(0,1)} |\tilde{P}^k|^\frac43\Big)^\frac34=1,
\end{equation}
\begin{equation}\label{n-m-t-1}
\int_{P_{\frac34}^+(0,1)} \big(|\tilde{n}^k_t|^2+|\nabla^2 \tilde{n}^k|^2+|\tilde{m}^k_t|^2+|\nabla^2 \tilde{m}^k|^2\big)\le C,
\end{equation}
and by Poincare's inequality,
$$
\int_{P_\frac34^+(0,1)} \big(|\tilde{n}^k|^4+|\tilde{m}^k|^4\big)\le C,
$$
but
\begin{eqnarray}\label{bdry_decay2}
Y((0,1)^+,\theta; \tilde{u}^k, \tilde{n}^k, \tilde{m}^k, \tilde{P}^k)>\theta^\frac23.
\end{eqnarray}
We may assume that there exists 
$\tilde{u}\in L^4(P_\frac34^+(0,1)), \tilde{n}\in W^{2,1}_2(P_\frac34^+(0,1)), 
\tilde{m}\in W^{2,1}_2(P_\frac34^+(0,1))$, $\tilde{P}\in L^\frac43(P_\frac34^+(0,1))$, $b\in\mathbb R^2$ with
$|b|\le M$, and $p, q\in \mathbb S^2$ with $p\cdot q=0$ such that after taking a
subsequence, 
\begin{equation*}
\begin{cases}
\tilde{u}^k\rightharpoonup \tilde{u}& \  \ {\rm{in}}\ \ L^4(P_\frac34^+(0,1)),\\
(\tilde{n}^k, \tilde{m}^k)\rightharpoonup (\tilde{n}, \tilde{m})& \  \ {\rm{in}}\ \ W^{2,1}_2(P_\frac34^+(0,1)),\\
\tilde{P}^k\rightharpoonup \tilde{P}& \  \ {\rm{in}}\ \ L^\frac43(P_\frac34^+(0,1)),\\
(\nabla\tilde{n}^k, \nabla\tilde{m}^k)\rightarrow(\nabla\tilde{n}, \nabla\tilde{m})& \ \ {\rm{in}}\ \ L^4(P_\frac34^+(0,1)),\\
(u^k)_{(0,1)^+, 1}\rightarrow b,\ \  (n_0^k, m_0^k)\rightarrow (p,q)\ \ {\rm{in}}\ \ C^1(T_1).
\end{cases} 
\end{equation*}
By the lower semicontinuity, we have that
\begin{equation}\label{bdry_small3}
Y((0,1)^+, \frac34; \tilde{u}, \tilde{n}, \tilde{m}, \tilde{P})\le 1,
\end{equation}
and
\begin{equation}\label{n-m-t-2}
\int_{P_{\frac34}^+(0,1)} \big(|\tilde{n}_t|^2+|\nabla^2 \tilde{n}|^2+|\tilde{m}_t|^2+|\nabla^2 \tilde{m}|^2\big)\le C.
\end{equation}
Since $(\tilde{u}^k, \tilde{n}^k, \tilde{m}^k, \tilde{P}^k)$ solves in $P_\frac34^+(0,1)$ the system
\begin{equation}\label{biaxial-nlcf3}
\begin{cases}
\tilde{u}^k_t+(\epsilon_i\tilde{u}^k+(u^k)_{(0,1)^+, 1})\cdot\nabla \tilde{u}^k+\nabla \tilde{P}^k-\Delta \tilde{u}^k\\
=-\epsilon_k{\rm{div}}\big(\nabla \tilde{n}^k\odot \nabla \tilde{n}^k+\nabla \tilde{m}^k\odot \nabla \tilde{m}^k\big),\\
{\rm{div}} (\tilde{u}^k)=0,\\
\tilde{n}^k_t+(\epsilon_k\tilde{u}^k+(u^k)_{(0,1)^+, 1})\cdot\nabla\tilde{n}^k-\Delta\tilde{n}^k
=\epsilon_k|\nabla\tilde{n}^k|^2{n^k}+\epsilon_k(\nabla\tilde{n}^k\cdot\nabla\tilde{m}^k) {m^k},\\
\tilde{m}^k_t+(\epsilon_k\tilde{u}^k+(u^k)_{(0,1)^+, 1})\cdot\nabla\tilde{m}^k-\Delta\tilde{m}^k
=\epsilon_k|\nabla\tilde{m}^k|^2{m^k}+\epsilon_k(\nabla\tilde{n}^k\cdot\nabla\tilde{m}^k) {n^k},
\end{cases}
\end{equation}
we obtain, after sending $k\rightarrow\infty$, that $(\tilde{u}, \tilde{n}, \tilde{m},
\tilde{P})$ solves the linear system \eqref{limit-eqn} in $P_\frac34^+(0,1)$, along with the boundary condition:
$$
(\tilde{u}, \tilde{n}, \tilde{m})=(0, 0, 0) \  \ {\rm{on}}\ \ T_\frac34\times \big[1-(\frac34)^2, 1\big].
$$
By the boundary regularity of heat equations, we have that
\begin{equation*}
\int_{P_\theta^+(0,1)} \big(|\nabla \tilde{n}|^4+|\nabla\tilde{m}|^4\big) 
\le C\theta^4\int_{P_\frac34^+(0,1)} \big(|\nabla \tilde{n}|^4+|\nabla\tilde{m}|^4\big)\le C\theta^4,
\ \forall\theta\in (0,\frac12).
\end{equation*}
Hence we have that
\begin{equation}\label{limit-reg1}
\int_{P_\theta^+(0,1)} \big(|\nabla \tilde{n}^k|^4+|\nabla\tilde{m}^k|^4\big) 
\le C\theta^4+o(i^{-1}),
\ \forall\theta\in (0,\frac12).
\end{equation}
To estimate $(\tilde{u}, \tilde{P})$, we recall the boundary local energy inequality \eqref{LED1}
for $(u^k, n^k, m^k, P^k)$ implies that $(\tilde{u}^k, \tilde{n}^k, \tilde{m}^k, \tilde{P}^k)$ satisfies
\begin{eqnarray}\label{LED2}
&&\sup_{t\in[\frac34, 1]}\int_{B_\frac12^+} |\tilde{u}^k(t)|^2+\int_{P_\frac12^+(0,1)}|\nabla\tilde{u}^k|^2\nonumber\\
&&\le C\int_{P_\frac34^+(0,1)} \big[(\epsilon_k|\tilde{u}^k|+M) (|\tilde{u}^k|^2+|\nabla\tilde{n}^k|^2+|\nabla\tilde{m}^k|^2)+(|\nabla \tilde{u}^k|+|\tilde{P}^k)|\tilde{u}^k|\nonumber\\
&&\ \ +|\nabla\tilde{n}^k||\tilde{n}^k_t|+|\nabla\tilde{m}^k||{\tilde m}^k_t|\big]
+\inf_{s\in [\frac34,1]} \int_{B_\frac34^+} \big(|\tilde{u}^k|^2+|\nabla\tilde{n}^k|^2+|\nabla\tilde{m}^k|^2\big)(s).
\end{eqnarray}
This, combined with \eqref{bdry_small2} and \eqref{n-m-t-1}, implies that 
\begin{equation}\label{LH-est}
\sup_{t\in[\frac34, 1]}\int_{B_\frac12^+} |\tilde{u}^k(t)|^2+\int_{P_\frac12^+(0,1)}|\nabla\tilde{u}^k|^2\le C.
\end{equation} 
Observe that the equation \eqref{biaxial-nlcf1} provides the following estimate
\begin{equation}\label{u-t}
\big\|\tilde{u}^k_t\big\|_{L^2([1-(\frac12)^2, 1], H^{-1}(B_\frac12^+))}\le C.
\end{equation}
Thus we have
\begin{equation}\label{limit-H1-bound}
\|\tilde{u}\|_{L^\infty_tL^2_x\cap L^2_tH^1_x(P_\frac12^+(0,1))}
+\|\tilde{u}_t\|_{L^2_tH^{-1}_x(P_\frac12^+(0,1))}\le C,
\end{equation}
and, by Aubin-Lions' compactness Lemma, 
\begin{equation}\label{u-l4-conv}
\tilde{u}^k\rightarrow \tilde{u} \ \ {\rm{in}}\ \ L^4(P_\frac12^+(0,1)).
\end{equation}
Since $(\tilde{u}, \tilde{P})$ solves the Stokes equation
$$\begin{cases}
\tilde{u}_t-\Delta\tilde{u}+\nabla\tilde{P}=0, \ \nabla\cdot\tilde{u}=0 & {\rm{in}}\ \ P_\frac34^+(0,1),\\
\tilde{u}=0 & \ {\rm{on}}\ \ T_\frac34\times [1-(\frac34)^2, 1],
\end{cases}
$$
it follows from Lemma 3.2 of \cite{SSS}  that $\tilde{u}\in W^{2,1}_{q, \frac43}(P_\frac12^+(0,1))$
and $\tilde {P}\in W^{1,0}_{q, \frac43}(P_\frac12^+(0,1))$ for any $1<q<\infty$, and
\begin{equation}
\big\|\tilde{u}\big\|_{W^{2,1}_{q, \frac43}(P_\frac12^+(0,1))} +\big\|\tilde {P}\big\|_{W^{1,0}_{q, \frac43}(P_\frac12^+(0,1))}
\le C(q) \big(\|\nabla\tilde{u}\|_{L^\frac43(P_\frac34^+(0,1))}+ \|\tilde{P}\|_{L^\frac43(P_\frac34^+(0,1))}\big)\le C(q).
\end{equation} 
This, combined with the Sobolev embedding theorem, implies that
$\tilde{u}\in C^{\alpha,\frac{\alpha}2}(P_\frac12^+(0,1))$, where $\alpha=2-\frac{2}{q}-\frac{2}{\frac43}=\frac12-\frac{2}{q}$
and $q>2$, and
\begin{equation}\label{u-bound}
\|\tilde{u}\|_{L^\infty(P_r^+(0,1))} \le C(\alpha) r^\alpha,  \ \forall \alpha\in (0,\frac12), \ 0<r\le \frac12.
\end{equation} 
This, combined with \eqref{u-l4-conv}, yields that for $k$ sufficiently large,
\begin{equation}\label{u-l4-decay}
\int_{P_\theta^+(0,1)}|\tilde{u}^k|^4\le C\theta^\alpha+o(i^{-1}), \ \forall\theta\in (0,\frac12).
\end{equation}
To estimate $\tilde{P}^k$, let $\hat{u}^k\in W^{2,1}_\frac43(\mathbb R^2_+\times (0,1))$ and $\hat{P}^k\in W^{1,0}_\frac43(\mathbb R^2_+\times (0,1))$ be the unique solution of
$$\begin{cases}
\hat{u}^k_t-\Delta\hat{u}^k+\nabla\hat{P}^k=\epsilon_k\nabla\cdot\big(\tilde{u}^k\otimes\tilde{u}^k
+\nabla \tilde{n}^k\odot \nabla \tilde{n}^k 
+\nabla \tilde{m}^k\odot \nabla \tilde{m}^k\big)\chi_{P_\frac34^+(0,1)} \ \ \ \ {\rm{in}}\ \ \ \ \ \mathbb R^2_+\times (0,1), \\
\nabla\cdot\hat{u}^k=0   \qquad\qquad\qquad\qquad\qquad\qquad\qquad\ \  \ {\rm{in}}\  \ \ \ \mathbb R^2_+\times (0,1),\\
\hat{u}^k=0  \qquad\qquad\qquad\qquad\qquad\qquad\qquad\qquad \ {\rm{on}}\ \ (T_\frac34\times [1-(\frac34)^2, 1])\cup  (\mathbb R^2_+\times\{0\}).
\end{cases}
$$
Applying Lemma 3.1 of \cite{SSS},  we have
\begin{eqnarray}\label{auxi1-est}
&&\big\|\hat{u}^k\big\|_{W^{2,1}_\frac43(\mathbb R^2_+\times (0,1))}
+\big\|\nabla\hat{P}^k\big\|_{L^\frac43(\mathbb R^2_+\times (0,1))}\nonumber\\
&&\le C\epsilon_k\big\|\nabla\cdot\big(\tilde{u}^k\otimes\tilde{u}^k
+\nabla \tilde{n}^k\odot \nabla \tilde{n}^k 
+\nabla \tilde{m}^k\odot \nabla \tilde{m}^k\big)\big\|_{L^\frac43(P_\frac34^+(0,1))}\nonumber\\
&&\le C\epsilon_k\Big(\|\tilde{u}^k\|_{L^4(P_\frac34^+(0,1))}\|\nabla\tilde{u}^k\|_{L^2(P_\frac34^+(0,1))}
+\|\nabla\tilde{n}^k\|_{L^4(P_\frac34^+(0,1))}\|\nabla^2\tilde{n}^k\|_{L^2(P_\frac34^+(0,1))}\nonumber\\
&&\qquad+
\|\nabla\tilde{m}^k\|_{L^4(P_\frac34^+(0,1))}\|\nabla^2\tilde{m}^k\|_{L^2(P_\frac34^+(0,1))}\Big)\le C\epsilon_k.
\end{eqnarray}
Let $v^k=\tilde{u}^k-\hat{u}^k$ and $Q^k=\tilde{P}^k-\hat{P}^k$. Then 
$$
v^k_t-\Delta v^k+\nabla Q^k=0, \ \nabla\cdot v^k=0 \ \ {\rm{in}}\ \ P_\frac34^+(0,1);\ \
v^k=0 \ \ {\rm{on}}\ \ T_\frac34\times [1-(\frac34)^2,1].
$$
By Lemma 3.1 of \cite{SSS}, we have that
\begin{eqnarray*}
\big\|\nabla Q^k\big\|_{L^\frac43_tL^{q}_x(P_\frac12^+(0,1))}\le C(q)\big(\|\nabla v^k\|_{L^\frac43(P_\frac34^+(0,1))}
+\|Q^k\|_{L^\frac43(P_\frac34^+(0,1))}\big)\le C, \ \forall q>1.
\end{eqnarray*}
Now we can apply H\"older and Poincar\'e's inequality to conclude that for any $\theta\in (0,\frac12)$, 
\begin{eqnarray}\label{P-43-decay}
Y^2((0,1)^+, \theta; P)&=&\theta^{-2} \Big(\int_{P_\theta^+(0,1)} \big|\tilde{P}^k-[\tilde{P}^k]_{(0,1)^+, \theta}\big|^\frac43\Big)^{\frac34}\nonumber\\
&\le& C\theta^{-1}\Big(\int_{P_\theta^+(0,1)}
|\nabla\tilde{P}^k|^\frac43\Big)^\frac34\nonumber\\
&\le& C\theta^{-1}\Big[\big(\int_{P_\frac12^+(0,1)}|\nabla \hat{P}^k|^\frac43\big)^\frac34
+\big(\int_{P_\theta^+(0,1)}|\nabla Q^k|^\frac43\big)^\frac34\Big]\nonumber\\
&\le& C\theta^{-1}\epsilon_k + C\theta^{\frac12-\frac2{q}} \|\nabla Q^k\|_{L^\frac43_tL^q_x(P_\theta^+(0,1))}\nonumber\\
&\le& C(\theta^{-1}\epsilon_k+\theta^\frac25),
\end{eqnarray}
provided $q\ge 5$. Putting \eqref{P-43-decay}, \eqref{u-l4-decay}, and \eqref{limit-reg1} together, we obtain
$$
Y((0,1)^+, \theta; \tilde{u}^k, \tilde{n}^k, \tilde{m}^k, \tilde{P}^k)
\le C(\theta^{-1}\epsilon_k+\theta^\frac25+o(k^{-1}))<\theta^\frac13,
$$
provided that $k$ is chosen sufficiently large and $\theta$ is chosen to be sufficiently small.  
This contradicts to \eqref{bdry_decay2}. 
\qed

\bigskip
Similar to Proposition \ref{interior-reg}, after iterating the estimates in Lemma \ref{bdry_reg}
we can obtain the boundary regularity of weak solutions to \eqref{biaxial-nlcf0}. More precisely, we have

\begin{proposition} \label{bdry_reg1} Under the same assumptions as in Lemma \ref{bdry_reg}, $(u, n, m)\in C^\infty\big(P_{\frac{r_0}2}^+(z_0)\big)\cap C^{2+\alpha, 1+\frac{\alpha}2}\big(\overline{P_{\frac{r_0}2}^+(z_0)}\big)$ and 
\begin{eqnarray}\label{bdry_est1} 
\big\|(u, n, m)\big\|_{C^{2+\alpha,1+\frac{\alpha}2}(P_{\frac{r_0}2}(z_0))}\le C\big(\epsilon_0, r_0, \|(n_0, m_0)\|_{C^{2,\alpha}(T_{r_0}(x_0)}\big).
\end{eqnarray}
\end{proposition} 
\pf The proof is similar to that of Proposition \ref{interior-reg}. First by iterating the argument of Lemma \ref{bdry_reg}
finitely many times and applying Propostion \ref{interior-reg}, we obtain that there exists $\alpha_1\in (0,1)$ such that for any $0<r\le r_0$ and $z\in P_{\frac{r_0}2}^+(z_0)$,
$$ 
Y(z^+, r; u, n, m , P)\le C\big(\frac{r}{r_0}\big)^{\alpha_1} \max\big\{Y(z_0, r_0; u, n, m, P), 
\ C_0\|(n_0,m_0)\|_{C^{2,\alpha}(T_{r_0}(x_0))}\big\}.
$$
Recall that $(n, m)$ satisfies the local energy inequality near the boundary:
$$
\int_{P_{\frac{r}2}^+(z)} (|n_t|^2+|m_t|^2)
\le C\int_{P_{r}^+(z)}[|u-u_{z^+, r}|^4+(M+1) |\nabla n|^4+|\nabla m|^4], \ \forall z\in P_{\frac{r_0}2}^+(z_0),
\ 0<r\le \frac{r_0}2.
$$
Thus we obtain that 
$$
\begin{cases}
\displaystyle  r^{-4}\int_{P_r^+(z)} |u-u_{z^+,r}|^4\le C\big(\frac{r}{r_0}\big)^{4\alpha_1},\\
\displaystyle  r^{-2}\int_{P_r^+(z)}(r^2(|n_t|^2+|m_t|^2)+(|\nabla n|^2+|\nabla m|^2))\le C\big(\frac{r}{r_0}\big)^{2\alpha_1}, 
\end{cases}
$$
for all $z\in  P_{\frac{r_0}2}(z_0), \ 0<r\le\frac{r_0}2$. Hence $(u, n, m)\in C^{\alpha_1,\frac{\alpha_1}2}(P_{\frac{r_0}2}^+(z_0))$.
Since $(n_0, m_0)\in C^{2+\alpha}(T_{r_0}(x_0))$, 
the higher order regularity of $(u, n, m)$ can be obtained by first obtaining the $C^{2+\alpha, 1+\frac{\alpha}2}$-regularity of
$(n, m)$ via suitable modifications of the boundary regularity of the heat flow of harmonic map, and then establishing the $C^{2+\alpha, 1+\frac{\alpha}2}$-regularity of $u$ via the boundary regularity of Stokes equations with H\"older continuous force functions. Moreover, the estimate \eqref{bdry_est1}  holds. \qed

\bigskip
As a consequence of both Proposition \ref{interior-reg} and Proposition \ref{bdry_reg1}, we have the following $\epsilon_0$-regularity.
\begin{corollary} \label{global_reg1} For $\alpha\in (0,1)$ and $T>0$, let $n_0,m_0\in H^1(\Omega, \mathbb S^2)\cap C^{2+\alpha}(\partial\Omega,\mathbb S^2)$ be such that $n_0\cdot m_0=0$ in $\overline\Omega$. There exist $\epsilon_0>0$
and $r_0>0$, depending on $\partial\Omega$ and $n_0, m_0$,
such that if $u\in L^\infty_tL^2_x\cap L^2_tH^1_x(Q_T, \mathbb R^2)$,
$n, m\in L^\infty_t H^1_x(Q_T, \mathbb S^2)$ with $n\cdot m=0$ in $Q_T$, $P\in W^{1,0}_{\frac43}(Q_T)$,
is a weak solution of \eqref{biaxial-nlcf0}, with $(u, n, m)=(0,n_0, m_0)$ on $\partial\Omega$, that satisfies, 
for $x_0\in\overline\Omega$ and $t_0\in (r_0^2, T]$,
\begin{equation}\label{global_small}
\Big(\int_{\Omega\times [0, T]\cap P_{r_0}(z_0)}(|u|^4+|\nabla n|^4+|\nabla m|^4) \Big)^\frac14
+\Big(\int_{\Omega\times [0, T]\cap P_{r_0}(z_0)} |\nabla P|^\frac43\Big)^\frac34
\le\epsilon_0,
\end{equation}
then
\begin{itemize}
\item [(a)] if $B_{r_0}(x_0)\subset\Omega$,  $(u, n, m)\in C^\infty(\Omega\times (0,T)\cap P_{{r_0}}(z_0))$,
and 
$$\big\|(u, n, m)\big\|_{C^l(P_{\frac{r_0}2}(z_0))}\le C(\epsilon_0, l), \forall l\ge 1.$$
\item [(b)] if $B_{r_0}(x_0)\cap\partial\Omega\not=\emptyset$, 
$(u, n, m)\in C^\infty(\Omega\times (0,T)\cap P_{{r}}(z_0))\cap C^{2+\alpha, 1+\frac{\alpha}2}(\overline\Omega\times (0,T)\cap P_r(z_0))$ for all $0<r<r_0$, and
$$\big\|(u, n, m)\big\|_{C^{2+\alpha, 1+\frac{\alpha}2}(\overline\Omega\times (\delta,T)\cap
P_{r}(z_0))}\le C(\epsilon_0,\delta, \|(n_0, m_0)\|_{C^{2+\alpha}(\partial\Omega)}, r), \forall 0<r<r_0.$$
\end{itemize}
\end{corollary}
\pf (a) It is readily seen that the condition \eqref{global_small} implies that
$$r_0|(u)_{z_0,r_0}|\le \epsilon_0, \ \ r_0Y(z_0, r_0; u, n, m, P)\le \epsilon_0$$
so that Proposition \ref{interior-est} yields the conclusion of (a) holds.\\
(b) In this case, it is also easy to see that
$$r_0(u)_{z_0^+, r_0}\le \epsilon_0, \ \ r_0 Y(z_0^+, r_0; u, n, m, P)\le \epsilon_0$$
so that Proposition \ref{bdry_reg1} and Proposition \ref{interior-est} yield
the conclusion of (b) holds. 
\qed

\medskip
As a direct consequence of Corollary 5.5 and Ladyzhenskaya's inequality, we have the following global regularity
result.
\begin{corollary}\label{global_reg2} For $\alpha\in (0,1)$ and $T>0$, let $n_0,m_0\in H^1(\Omega, \mathbb S^2)\cap C^{2+\alpha}(\partial\Omega,\mathbb S^2)$ be such that $n_0\cdot m_0=0$ in $\overline\Omega$. If $u\in L^\infty_tL^2_x\cap L^2_tH^1_x(Q_T, \mathbb R^2)$,
$n, m\in L^\infty_t H^1_x(Q_T, \mathbb S^2)\cap L^2_tH^2_x(Q_T, \mathbb S^2)$ with $n\cdot m=0$ in $Q_T$, 
and $P\in W^{1,0}_\frac43(Q_T)$,
is a weak solution of \eqref{biaxial-nlcf0}, with $(u, n, m)=(0,n_0, m_0)$ on $\partial\Omega$,
then $(u, n, m)\in C^\infty(\Omega\times (0, T))\cap C^{2+\alpha, 1+\frac{\alpha}2}(\overline\Omega\times (0,T))$.
\end{corollary}
\pf Applying Ladyzhenskaya's inequality, we have that $u\in L^4(Q_T), \nabla n\in L^4(Q_T),$ and
$\nabla m\in L^4(Q_T)$. By the absolute continuity, we see that there exists $r_0>0$ such that 
the smallness condition \eqref{global_small} holds. Hence the conclusion of Corollary \ref{global_reg2}
follows from Corollary \ref{global_reg1}. This completes the proof.  \qed

\section{Existence of global weak solutions}
\setcounter{equation}{0}

The goal of this section is to establish the global existence of weak solutions of \eqref{biaxial-nlcf0}. 
With the estimates obtained from the previous sections, we can adapt the scheme from \cite{LLW10} to
achieve this. More precisely, this involves establishing a uniform lower bound of time intervals of short time smooth solutions depends only on the local energy profile radius of the initial data $(u_0, n_0, m_0)$. Then we characterize finite singular time in term of certain nontrivial global harmonic maps from $\mathbb S^2$ to the manifold $\mathcal{N}$, introduced in the section one.

\medskip
\noindent{\it Proof of Theorem \ref{LH}}. Since this can be done similar to the scheme of \cite{LLW10},   we only sketch it here. It is divided into several steps.

\noindent{\it Step 1} (lower bound of time intervals). By Lemma \ref{approx}, there exists 
$u_0^k\in C^\infty_0(\Omega, \mathbb R^2)$ with $\nabla\cdot u_0^k=0$, 
$(n_0^k, m_0^k)\in C^\infty(\Omega,\mathbb S^2\times\mathbb S^2)\cap C^{2+\alpha}(\overline\Omega,\mathbb S^2\times\mathbb S^2)$, with $(n_0^k, m_0^k)=(n_0, m_0)$ on $\partial\Omega$ and $n_0^k\cdot m_0^k=0$
in $\Omega$, such that
\begin{equation}\label{strong-approx}
\lim_{k\rightarrow \infty}\big(\big\|u_0^k-u_0\big\|_{L^2(\Omega)}+\big\|(n_0^k-n_0,  m_0^k-m_0)\big\|_{H^1(\Omega)}\big)=0.
\end{equation}

By Theorem \ref{local}, there exist $T_k>0$ and a unique smooth solution $(u^k, n^k, m^k)
\in \big(C^\infty(Q_{T_k}, \mathbb R^2)\cap C^{2+\alpha, 1+\frac{\alpha}2}(\overline{Q_{T_k}},\mathbb R^2)\big)
\times \big(C^\infty(Q_{T_k}, \mathbb S^2)\cap C^{2+\alpha, 1+\frac{\alpha}2}(\overline{Q_{T_k}},\mathbb S^2)\big)\times \big(C^\infty(Q_{T_k}, \mathbb S^2)\cap C^{2+\alpha, 1+\frac{\alpha}2}(\overline{Q_{T_k}},\mathbb S^2)\big)$ 
with $n^k\cdot m^k=0$ in $Q_{T_k}$ of \eqref{biaxial-nlcf0}, along with the initial and boundary data $(u_0^k, n_0^k, m_0^k)$. 

For any small $\epsilon_0>0$, it follows from \eqref{strong-approx} that there exists $r_0>0$ such that
$$\mathcal{E}(k,r_0; u_0^k, n_0^k, m_0^k)\equiv
\max_{x\in\overline\Omega} \int_{\Omega\cap B_{2r_0}(x)}(|u_0^k|^2+|\nabla n_0^k|^2+|\nabla m_0^k|^2)
\le \frac{\epsilon_0^2}{2}, \ \ \forall k\ge 1.$$
We now claim that there exists a $c_0>0$ such that $T_k\ge c_0r_0^2$ holds for all $k$. To see this,
let $t_0\in (0, T_k]$ be the maximal time such that 
$$ 
\begin{cases}
\mathcal{E}(k, r_0; u^k(t), n^k(t), m^k(t))\le \epsilon_0^2, \ \forall 0\le t\le t_0,\\
\mathcal{E}(k, r_0; u^k(t_0), n^k(t_0), m^k(t_0))= \epsilon_0^2.
\end{cases}
$$
As in \cite{LLW10} (see Lemma 5.1 and Lemma 5.2), we can apply Ladyzhenskaya's inequality
and Lemma \ref{global_energy_ineq} to deduce that 
\begin{equation}\label{high_est1}
\begin{cases}
\displaystyle\int_{Q_{t_0}} (|\nabla u^k|^2+|\nabla^2 n^k|^2+|\nabla^2 m^k|^2) \le CE_0,\\
\displaystyle\int_{Q_{t_0}} (|u^k|^4+|\nabla n^k|^4+|\nabla m^k|^4)\le C\epsilon_0^2E_0,\\
\displaystyle\int_{Q_{t_0}} (|u^k_t|^\frac43+|n^k_t|^2+|m^k_t|^2) \le CE_0.
\end{cases}
\end{equation}
Here $E_0=\int_{\Omega} (|u_0|^2+|\nabla n_0|^2+|\nabla m_0|^2)$ denotes the energy of initial data.
 Note that in the proof of Theorem \ref{local}, we also have the estimate (see also \cite{LLW10} Lemma 4.4)
 \begin{eqnarray}\label{pressure_est1}
&& \big\|\nabla P^k\big\|_{L^\frac43(Q_{t_0})}\nonumber\\
 &&\le C\big[\|u^k\|_{L^4(Q_{t_0})}\|\nabla u^k\|_{L^2(Q_{t_0})}
 + \|\nabla n^k\|_{L^4(Q_{t_0})}\|\nabla^2 n^k\|_{L^2(Q_{t_0})}+\|\nabla m^k\|_{L^4(Q_{t_0})}
 \|\nabla^2 m^k\|_{L^2(Q_{t_0})}\big]\nonumber\\
 &&\le C\epsilon_0^\frac12 E_0^\frac34.
 \end{eqnarray}
With the estimates \eqref{high_est1} and \eqref{pressure_est1}, we can apply the local energy inequalities
in both Lemma \ref{interior_local} and Lemma \ref{bdry_local} to deduce that  
\begin{eqnarray*}
&&\mathcal{E}(k, r_0, u^k(t_0), n^k(t_0),  m^k(t_0))-\mathcal{E}(k, 2r_0) \\
&&\le C(\frac{t_0}{r_0^2})^\frac14\Big[\|u^k\|_{L^4(Q_{t_0})}^3+\|u^k\|_{L^4(Q_{t_0})}\big(\|\nabla u^k\|_{L^2(Q_{t_0})}
+\|\nabla n^k\|_{L^4(Q_{t_0})}^2
+\|\nabla m^k\|_{L^4(Q_{t_0})}^2\big)\\
&&\ \ +\|n^k_t\|_{L^2(Q_{t_0})}\|\nabla n^k\|_{L^4(Q_{t_0})}
+\|m^k_t\|_{L^2(Q_{t_0})}\|\nabla m^k\|_{L^4(Q_{t_0})}
+\|\nabla P^k\|_{L^\frac43(Q_{t_0})}\|u^k\|_{L^\infty_tL^2(Q_{t_0})}\Big]\\
&&\le C(\frac{t_0}{r_0^2})^\frac14\epsilon_0^\frac12.
\end{eqnarray*}
Thus if $\mathcal{E}(k, r_0, u^k(t_0), n^k(t_0),  m^k(t_0))=2\epsilon_0^2,$
then we must have that $\epsilon_0^2\le C(\frac{t_0}{r_0^2})^\frac14\epsilon_0^\frac12$, which implies
$t_0\ge C\epsilon_0^6 r_0^2$ so that $T_k\ge C\epsilon_0^6 r_0^2$  and the claim is true. 

From \eqref{high_est1}, we can apply Corollary \ref{global_reg2} to conclude that 
\begin{equation}\label{globaL_est3}
\big\|(u^k, n^k, m^k)\big\|_{C^{2+\alpha, 1+\frac{\alpha}2}(\Omega\times [\delta, {T_k}]}
\le C\big(\epsilon_0, \delta, \|(n_0, m_0)\|_{C^{2+\alpha}(\partial\Omega)}\big), \ \forall k\ge 1.
\end{equation} 

Hence there exist $T_0\ge c_0 r_0^2$, $u\in L^\infty_tL^2_x\cap L^2_tH^1_x(Q_{T_0},\mathbb R^2)$,
$n,m\in L^\infty_tH^1_x\cap L^2_tH^2_x(Q_{T_0},\mathbb \mathbb S^2)$ with $n\cdot m=0$ in $Q_{T_0}$,
and $P\in W^{1,0}_{\frac43}(Q_{T_0})$ such that, after passing to a subsequence, 
$$
\begin{cases}
u^k\rightharpoonup u \ {\rm{in}}\ W^{1,0}_2(Q_{T_0}), \ \ u^k\rightarrow u  \ {\rm{in}}\ L^4(Q_{T_0}), \\
(n^k, m^k) \rightharpoonup (n, m) \ {\rm{in}}\ W^{2,1}_2(Q_{T_0}), \ 
(\nabla n^k, \nabla m^k)\rightarrow (\nabla n, \nabla m) \ {\rm{in}} \  L^4(Q_{T_0}),\\
P^k\rightharpoonup P \ {\rm{in}}\ W^{1,0}_{\frac43}(Q_{T_0}),
\end{cases}
$$
and for any compact set $K\subset\subset\Omega$, $\delta>0$,  and $l\ge 1$, 
$$(u^k, n^k, m^k) \rightarrow(u, n, m)\ \ {\rm{in}}\ \  C^{l}(K\times [\delta, T_0]),$$
and
$$(u^k, n^k, m^k) \rightarrow(u, n, m)\ \ {\rm{in}}\ \  
C^{2+\alpha, 1+\frac{\alpha}2}(\overline\Omega\times [\delta, T_0]).
$$
Hence $(u, n, m)\in C^\infty(\Omega\times (0, T_0))\cap C^{2+\alpha, 1+\frac{\alpha}2}(\overline\Omega\times [\delta, T_0))$ for any $\delta>0$.
Moreover, by the lower semicontinuity $(u, n, m, P)$ satisfies the energy inequality: for any $0<t\le T_0$, 
\begin{equation}\label{global_energy_ineq1}
\int_\Omega (|u|^2+|\nabla n|^2+|\nabla m|^2)+2\int_{Q_t} (|\nabla u|^2+|n_t+u\cdot\nabla n|^2+|m_t+u\cdot\nabla m|^2)\le E_0.
\end{equation} 
With \ref{global_energy_ineq1}, we can further show that as $t\downarrow 0^+$, 
$$(u(t), n(t), m(t))\rightarrow (u_0, n_0, m_0) \ {\rm{in}}\ L^2(\Omega)\times H^1(\Omega)\times H^1(\Omega).$$

\noindent{\it Step 2} (characterization of the first singular time). Suppose $T_1$ is the maximal time for the smooth solution $(u, n, m)$ obtained in step 1. Then we must have
\begin{equation} \label{concentration1}
\limsup_{t\uparrow T_1} \max_{x\in\overline\Omega} \int_{\Omega\cap B_r(x)}
(|u(t)|^2+|\nabla n(t)|^2+|\nabla m(t)|^2)\ge \epsilon_0^2,\ \ \forall r>0.
\end{equation}
Again, it follows from \eqref{global_energy_ineq1} and the weak continuity of $(u(t), \nabla n(t),\nabla m(t))$ in $L^2(\Omega)$ that 
$$(u(t), n(t),m(t))\rightharpoonup (u(T_1), n(T_1), m(T_1)) \ {\rm{in}}\ L^2(\Omega)\times H^1(\Omega)\times H^1(\Omega), \ {\rm{as}}\ t\uparrow T_1.$$
Moreover, one can verify that $(u(T_1), n(T_1), m(T_1))=(0, n_0, m_0)$ on $\partial\Omega$. 

Next we will use $(u(T_1), n(T_1), m(T_1))$ and $(0, n_0, m_0)$ as the initial and boundary data to repeat the above
procedure  to obtain a continuation of $(u, n, m)$ beyond $T_1$ as a weak solution
of \eqref{biaxial-nlcf0}. At any further singular time, we repeat this procedure. We will prove
that there are at most finitely many such singular times, and afterwards we will 
construct an eternal weak solution.

It follows from \eqref{concentration1} that for any $r>0$, 
\begin{eqnarray*}
\int_\Omega (|u|^2+|\nabla n|^2+|\nabla m|^2)(T_1) 
&\le&\liminf_{t\rightarrow T_1} \int_\Omega  (|u|^2+|\nabla n|^2+|\nabla m|^2)(t)\\
&& -\limsup_{t\rightarrow T_1}\max_{x\in\overline\Omega}\int_{\Omega\cap B_r(x)}(|u|^2
+|\nabla n|^2+|\nabla m|^2)(t)\\
&\le& E_0-\epsilon_0^2.
\end{eqnarray*} 
From this, we see that the number of finite singular times must be bounded by
$L = [E_0/\epsilon_0^2]$. If $0 < T_L < \infty$ is the last singular time, then we must have
$$E(T_L ) =\int_\Omega (|u|^2+|\nabla n|^2+|\nabla m|^2)(T_L)\le\epsilon_0^2.$$
If we use $(u, n, m)(T_L)$ and $(0, n_0, m_0)$
as the initial and boundary value to construct a weak solution 
$(u, n, m)$ to \eqref{biaxial-nlcf0} as before, then $(u, n, m)$ will be an
eternal weak solution.

Now, similar to \cite{LLW10},  we can perform a blowing up argument at any singular time $T_i$ for
$1\le i\le L$. For simplicity, we sketch it for $T_1$.  For any $0<\epsilon_1<\epsilon_0$, there exist
 $r_k\rightarrow 0$, and $0<t_0<t_k\uparrow T_1$ such that
$$\epsilon_1^2=\mathcal{E}(t_k; r_k)=\max_{x\in\overline\Omega} 
\int_{\Omega\cap B_{r_k}(x)}(|u|^2+|\nabla n|^2+|\nabla m|^2)(t_k).$$
This, combined with the claim in Step 1, implies that there exist $\theta_0\in (0,1)$
and $x_k\in\overline\Omega$ such that
$$\int_{\Omega\cap B_{2r_k}(x_k)}(|u|^2+|\nabla n|^2+|\nabla m|^2)(t_k-\theta_0 r_k^2)
\ge \frac12 \mathcal{E}(t_k-\theta_0^2 r_k^2; 2r_k)\ge\frac{\epsilon_1^2}{2}.$$
Hence we can deduce that
$$
\begin{cases}
\displaystyle\int_{\Omega\times [t_0, t_k]} (|u_t|^\frac43+|\nabla P|^\frac43+|n_t|^2+|\nabla^2 n|^2+|m_t|^2+|\nabla^2 m|^2)\le C,\\
%\displaystyle\int_{\Omega\times [t_0, t_k]}(|u_t|^\frac43+|\nabla P|^\frac43)\le C,\\
\displaystyle\int_{\Omega\times [t_0, t_k]} (|u|^4+|\nabla n|^4+|\nabla m|^4)\le C\epsilon_1^2.
\end{cases}
$$
Set  $\Omega_k=r_k^{-1}(\Omega\setminus \{x_k\})$ and $Q_k=\Omega\times[\frac{t_0-t_k}{r_k^2}, 0]$.
Define a blowing up sequence
$$(u^k, n^k, m^k, P^k)(x,t)=(r_k u, n, m, r_k^2 P)(x_k+r_k x, t_k+r_k^2 t): Q_k \to \mathbb R^8.$$
Then $(u^k, n^k, m^k, P^k)$ solves \eqref{biaxial-nlcf0} in $Q_k$, and
\begin{equation}\label{blow_bound}
\begin{cases}
\displaystyle\int_{Q_k} (|u^k_t|^\frac43+|\nabla P^k|^\frac43+
|n^k_t|^2+|\nabla^2 n^k|^2+|m^k_t|^2+|\nabla^2 m^k|^2)\le C,\\
\displaystyle\int_{Q_k} (|u^k|^4+|\nabla n^k|^4+|\nabla m^k|^4)\le C\epsilon_1^2,\\
\displaystyle\int_{\Omega_k\cap B_2(0)} (|u^k|^2+|\nabla n^k|^2+|\nabla m^k|^2)(-\theta_0) \ge \frac{\epsilon_1^2}{2},\\
\displaystyle\sup_{(x,t)\in Q_k} \int_{\Omega_k\cap B_1(x)} (|u^k|^2+|\nabla n^k|^2+|\nabla m^k|^2)(t) \le \epsilon_1^2.
\end{cases}
\end{equation} 
Assume $x_k\rightarrow x_0\in\overline\Omega$. We divide it into two cases:
\begin{itemize}
\item [a)] $x_0\in\Omega$. Then $\Omega_k\rightarrow \mathbb R^2$. We may assume that there exists
$(u^\infty, n^\infty, m^\infty)\in C^\infty\cap W^{1,0}_2(\mathbb R^2\times (-\infty, 0], \mathbb R^2)\times 
C^\infty\cap L^{\infty} H^1(\mathbb R^2\times (-\infty, 0], \mathbb S^2)\times C^\infty\cap L^{\infty}H^1(\mathbb R^2\times (-\infty, 0], \mathbb S^2)$, with $n_\infty\cdot m_\infty=0$, such that 
$$
\begin{cases}
(u^k, n^k, m^k)\rightarrow (u^\infty, n^\infty, m^\infty) & {\rm{in}}\  C^2_{\rm{loc}}(\mathbb R^2\times (-\infty, 0]),\\
(u^k, \nabla n^k, \nabla m^k)\rightarrow (u^\infty, \nabla n^\infty, \nabla m^\infty) & {\rm{in}}\  L^2_{\rm{loc}}(\mathbb R^2\times (-\infty, 0]).
\end{cases}
$$
As in \cite{LLW10}, we have that for any $R>0$,
\begin{eqnarray*}
&&\int_{P_R} (|\nabla u^\infty|^2+|n^\infty_t+u^\infty\cdot\nabla n^\infty|^2+|m^\infty_t+u^\infty\cdot\nabla m^\infty|^2)\\
&&\le \lim_{k\rightarrow\infty} \int_{P_R} (|\nabla u^k|^2+|n^k_t+u^k\cdot\nabla n^k|^2
+|m^k_t+u^k\cdot\nabla m^k|^2)\\
&&=\lim_{k\rightarrow\infty} \int_{P_{R r_k}(x_k, t_k)} (|\nabla u|^2+|n_t+u\cdot\nabla n|^2
+|m_t+u\cdot\nabla m|^2)=0.
\end{eqnarray*}
It follows that $u^\infty=n^\infty_t=m^\infty_t=0$ in $\mathbb R^2\times (-\infty, 0]$ and hence
$(n^\infty, m^\infty)$ is time independent. Moreover, it follows from
\eqref{blow_bound} that 
$$0<\int_{\mathbb R^2}(|\nabla n^\infty|^2+|\nabla m^\infty|^2)<\infty.$$ 
Hence $(n^\infty, m^\infty)$ is
a nontrivial smooth harmonic map from $\mathbb R^2\to \mathcal{N}=\big\{(p,q)\in \mathbb S^2\times \mathbb S^2: \ p\cdot q=0\big\}$ with finite energy:
$$\begin{cases}
-\Delta n^\infty=|\nabla n^\infty|^2 n^\infty +\langle \nabla n^\infty, \nabla m^\infty\rangle m^\infty,\\
-\Delta m^\infty=|\nabla m^\infty|^2 m^\infty +\langle \nabla n^\infty, \nabla m^\infty\rangle n^\infty.
\end{cases}
$$
It is well-known that the conformal invariance of harmonic maps in dimension two (see \cite{LWbook}) that there exists
a universal positive constant ${\mathcal{C}}_0$  given by  \eqref{energy_sp} such that
$$\int_{\mathbb R^2}(|\nabla n^\infty|^2+|\nabla m^\infty|^2)\ge {\mathcal{C}}_0.$$
\item [b)] $x_0\in\partial\Omega$. Then either (i) $\frac{|x_k-x_0|}{r_k}\rightarrow \infty$ or
(ii) $\frac{|x_k-x_0|}{r_k}\rightarrow c\ge 0$. If (i) happens, then $\Omega_k\rightarrow \mathbb R^2$ and it can be reduced to the case a). Now we will rule out (ii) as follows.  For, otherwise, assume
$\frac{x_k-x_0}{r_k}\rightarrow (0,c)$ and $\Omega_k\rightarrow \mathbb R^2_{c}=\big\{x=(x_1, x_2)\in\mathbb R^2: \ x_2\ge c\big\}$. Since $u^k=0$, $n^k=n_0(x_k+r_k\cdot)$, and $m^k=m_0(x_k+r_k\cdot)$ on $\partial\Omega_k$, 
as in the case a) we can show that $(u^k, n^k, m^k)\rightarrow (0, n^\infty, m^\infty)$ 
in $C^2_{\rm{loc}}(\mathbb R^2_c)$ such that $(n^\infty, m^\infty):\mathbb R^2_c\to \mathcal{N}$ is a
nontrivial smooth harmonic map of finite energy, with $(n^\infty, m^\infty)= (p,q)$ on $\partial \mathbb R^2_c$
for  a point $(p,q)\in\mathcal{N}$. This contradicts the nonexistence theorem by Lemaire \cite{Lemaire}. 
Thus b) is proven.
\end{itemize}

So far we have shown parts (i), (ii), (iii) of Theorem \ref{LH}. For (iv), observe the energy inequality for
$(u, n, m)$ implies that there exists $t_k\uparrow\infty$ such that
$$
\begin{cases}\displaystyle
\int_\Omega (|u|^2+|\nabla n|^2+|\nabla m|^2)(t_k)\le E_0,\\
\displaystyle\int_{\Omega} (|\nabla u|^2+|n_t+u\cdot \nabla n|^2+|m_t+u\cdot \nabla m|^2)(t_k)\rightarrow 0.
\end{cases}
$$
Denote $(u^k, n^k, m^k)=(u, n, m)(t_k)$. 
Since $u^k=0$ on $\partial\Omega$, we see that $u^k\rightarrow 0$ in $H^1(\Omega)$. Moreover,
$(n^k, m^k):\Omega\to{\mathcal{N}}$ satisfies
$$\begin{cases}
\Delta n^k+|\nabla n^k|^2 n^k +\langle \nabla n^k, \nabla m^k\rangle m^k
=f_k:=-(n_t+u\cdot \nabla n)(t_k),\\
\Delta m^k+|\nabla m^k|^2 m^k +\langle \nabla n^k, \nabla m^k\rangle n^k
=g_k:=-(m_t+u\cdot \nabla m)(t_k).
\end{cases}
$$
Hence $(n^k, m^k)$ is a sequence of approximated harmonic maps from $\Omega$ to $\mathcal{N}$
that has their tension fields bounded in $L^2(\Omega)$. Hence (iv) follows from the energy identity
result of approximated harmonic maps in dimension two, see \cite{Qing} and \cite{LW1}. It is readily seen
that (v) of Theorem \ref{LH} follows from (i), (ii), (iii), and (iv) of Theorem \ref{LH}. 

Part (vi) on the uniqueness of the constructed weak solutions will be proven in next section.
\qed

\section{Uniqueness}
\setcounter{equation}{0}
In this section, we will prove (vi) of Theorem \ref{LH} on the uniqueness of weak solutions in a certain class.
First we need the following Lemma. 

\begin{lemma}\label{H2-regularity} For $\alpha\in (0,1)$ and $T>0$, let $n_0,m_0\in H^1(\Omega, \mathbb S^2)\cap C^{2,\alpha}(\partial\Omega,\mathbb S^2)$ be such that $n_0\cdot m_0=0$ in $\overline\Omega$. If
$(u, {n}, {m})\in L^2_tH^1_x\cap L^\infty_t L^2_x(Q_T,\mathbb R^2)\times L^\infty_t H^1_x(Q_T,\mathbb S^2)\times  L^\infty_t H^1_x(Q_T,\mathbb S^2)$ with $n\cdot m=0$ in $Q_T$, and $P\in W^{1,0}_\frac43(Q_T)$, is a weak solution of \eqref{biaxial-nlcf0} along with the initial and boundary condition $(u_0, n_0, m_0)$, 
that enjoys the following properties:
\begin{itemize}
\item [(a)] $E(u(t), {n}(t), {m}(t); \Omega)$ is monotone nonincreasing for $t\in [0, \infty)$, and
\item [(b)] ${n}_t+u\cdot\nabla{n},\ {m}_t+u\cdot\nabla{m}\in L^2([0,T), L^2(\Omega))$
for any $0<T<\infty$.
\end{itemize}
Then there exist $\epsilon_0>0$, $r_0>0$, $C_0>0$ depending on $\partial\Omega$ and
$(u_0, n_0, m_0)$ such that $(u, n, m)\in L^2_tH^1_x(Q_{r_0^2})
\times L^2_tH^2_x(Q_{r_0^2})\times L^2_tH^2_x(Q_{r_0^2})$,
and  $(u, n,m)\in C^{2+\alpha, 1+\frac{\alpha}2}(\overline\Omega\times [\delta, r_0^2])$ for any
$0<\delta<r_0^2$. Moreover,
\begin{equation}\label{gradient_estimate}
\big|(u,\nabla n, \nabla m)\big|(x, t)\le \frac{C_0 \epsilon_0}{\sqrt{t}}, \ \forall (x,t)\in\overline\Omega\times(0, r_0^2].
\end{equation}
\end{lemma}
\pf From the condition (a) and the lower semicontinuity
$$E(u_0, n_0, m_0; \Omega)\le \liminf_{t\downarrow 0^+} E(u(t), {n}(t), {m}(t); \Omega),$$
we see that
$$
E(u_0, n_0, m_0; \Omega)= \lim_{t\downarrow 0^+} E(u(t), {n}(t), {m}(t); \Omega).
$$
Hence for any small $\epsilon_0>0$, there exists $r_0>0$ and $t_0=r_0^2$ such that
\begin{equation}\label{global_small3}
\sup_{(x,t)\in \overline\Omega\times [0, t_0]} E\big(u(t), {n}(t), {m}(t); \Omega\cap B_{r_0}(x)\big)\le \epsilon_0^2.
\end{equation}
For almost all $t\in (0, t_0)$, from the condition (b) we can view that $(n(t), m(t)):\Omega\to\mathcal{N}$ is an approximated harmonic map
with $L^2$-tension field:
$$\begin{cases}
\Delta n(t)+|\nabla n(t)|^2 n(t) +\langle \nabla n(t), \nabla m(t)\rangle m(t)
=f(t):=-(n_t+u\cdot \nabla n)(t),\\
\Delta m(t)+|\nabla m(t)|^2 m(t) +\langle \nabla n(t), \nabla m(t)\rangle n(t)
=g(t):=-(m_t+u\cdot \nabla m)(t)
\end{cases}
\ \ {\rm{in}}\ \ \Omega,
$$
subject to the boundary condition $(n(t), m(t))\big|_{\partial\Omega}=(n_0, m_0)$. 
By the higher order Sobolev regularity  theory on approximated harmonic maps in dimension two 
(see Wang \cite{Wang}, Sharp-Topping \cite{SP}, and Moser \cite{Moser}), we conclude that under
the condition \eqref{global_small3}, 
$(n(t), m(t))\in H^2(\Omega)$ and
\begin{eqnarray}\label{H2-estimate}
&&\big\|(n(t), m(t))\big\|_{H^2(\Omega)} \nonumber\\
&&\le C(r_0, \epsilon_0) 
\Big[\|(\nabla n(t),\nabla m(t))\|_{L^2(\Omega)} 
+\big\|\big((n_t+u\cdot \nabla n)(t), (m_t+u\cdot \nabla m)(t)\big)\big\|_{L^2(\Omega)}\Big].
\end{eqnarray}
Now we can apply Corollary \ref{global_reg2} to conclude that $(u, n, m)\in L^2_tH^1_x(Q_{r_0^2})
\times L^2_tH^2_x(Q_{r_0^2})\times L^2_tH^2_x(Q_{r_0^2})$,
and  $(u, n,m)\in C^{2+\alpha, 1+\frac{\alpha}2}(\overline\Omega\times [\delta, r_0^2])$ for any
$0<\delta<r_0^2$ along with the estimate \eqref{gradient_estimate}. This completes the proof.
\qed

\medskip
Based on Lemma \ref{H2-regularity}, we can give a proof of (vi) of Theorem \ref{LH} by suitable modifications of
\cite{LW2} as follows.  

\smallskip
\noindent{\it Proof of Theorem \ref{LH} (iv)}: For $i=1,2$, let $(u_i, n_i, m_i, P_i)$ be two solutions of 
\eqref{biaxial-nlcf0} in $Q_T$ with $(u_i, n_i, m_i)=(u_0, n_0, m_0)$ on $\partial_p Q_T$. Then by Lemma \ref{H2-regularity}, there exist $\epsilon_0>0$ and $t_0>0$ such that 
\begin{equation}\label{gradient_estimate1}
\max_{i=1,2}\big|(u_i,\nabla n_i, \nabla m_i)\big|(x, t)\le \frac{C_0 \epsilon_0}{\sqrt{t}}, \ \forall (x,t)\in\overline\Omega\times(0, t_0].
\end{equation}
Define $u=u_1-u_2, n=n_1-n_2, m=m_1-m_2, P=P_1-P_2$. Then $(u, n, m, P)$ satisfies, in $Q_{t_0}$,
\begin{equation}\label{eqn_diff}
\begin{cases} 
u_t-\Delta u+\nabla P=F\\
:=-u_1\cdot\nabla u-u\cdot\nabla u_2-\nabla\cdot\big(\nabla n_1\odot\nabla n
+\nabla n\odot\nabla n_2+\nabla m_1\odot\nabla m+\nabla m\odot\nabla m_2\big),\\
\qquad \nabla\cdot u=0,\\
n_t-\Delta n=G:=-u_1\cdot\nabla n-u\cdot\nabla n_2+|\nabla n_1|^2 n-\langle \nabla (n_1+n_2),\nabla n\rangle n_2\\
\qquad \qquad \qquad \qquad +\langle\nabla n_1,\nabla m_1\rangle m+(\langle\nabla n,\nabla m_1\rangle+\langle \nabla n_2,\nabla m\rangle) m_2,\\
m_t-\Delta m=H:=-u_1\cdot\nabla m-u\cdot\nabla m_2+|\nabla m_1|^2 m +\langle \nabla (m_1+m_2),\nabla m\rangle m_2\\
\qquad \qquad \qquad \qquad +\langle\nabla n_1,\nabla m_1\rangle n+(\langle\nabla n,\nabla m_1\rangle+\langle\nabla n_2, \nabla m\rangle) n_2.
\end{cases}
\end{equation}
Set
$$\mathcal{A}(t)=\sqrt{t}\sum_{i=1}^2\big(\|u_i(t)\|_{L^\infty(\Omega)}+\|\nabla n_i(t)\|_{L^\infty(\Omega)}
+\|\nabla m_i(t)\|_{L^\infty(\Omega)}\big),$$
$$
\mathcal{B}(t)=\sum_{i=1}^2\big(\|u_i(t)\|_{L^2(\Omega)}+\|\nabla n_i(t)\|_{L^2(\Omega)}
+\|\nabla m_i(t)\|_{L^2(\Omega)}\big),
$$
and for $0<\delta<1$,
$$
\mathcal{D}_\delta(t)=t^{\frac{1-\delta}{2}}\sum_{i=1}^2\big(\|u_i(t)\|_{L^{\frac{2}{\delta}}(\Omega)}+\|\nabla n_i(t)\|_{L^{\frac{2}{\delta}}(\Omega)}
+\|\nabla m_i(t)\|_{L^{\frac{2}{\delta}}(\Omega)}\big).
$$
Then it follows from \eqref{gradient_estimate1}  that 
\begin{equation}\label{gradient_estimate2}
\mathcal{A}(t)\le C\epsilon_0, \ \mathcal{B}(t)\le CE_0^\frac12, 
\ \mathcal{D}_\delta(t) \le C\epsilon^{1-\delta}, \ \forall 0<t\le t_0.
\end{equation}
Also set
$$
\Phi(t)=\sup_{0\le s\le t} \Big[
\|u(s)\|_{L^2(\Omega)}+\|(\nabla n(s),\nabla m(s))\|_{L^2(\Omega)}
+s^{-\frac{\delta}2}\|(n(s), m(s))\|_{L^{\frac{2}{\delta}}(\Omega)}\Big].
$$
Utilizing Duhamel's formula and various estimates on both the heat kernel and the Oseen-kernel, 
we can argue as in the proof of Lemma 3.1 of \cite{LW2}  to obtain:
\begin{eqnarray}\label{phi1}
\big\|(n(t), m(t))\big\|_{L^{\frac{2}{\delta}}(\Omega)}
\le Ct^{\frac{\delta}2}\Big[\sup_{0\le s\le t} \mathcal{D}_\delta(s)+(\sup_{0\le s\le t}\mathcal{A}(s))
(\sup_{0\le s\le t}\mathcal{B}(s))\Big] \Phi(t),
\end{eqnarray}
\begin{eqnarray}\label{phi2}
\big\|(\nabla n(t), \nabla m(t))\big\|_{L^{2}(\Omega)}
\le C\Big[\sup_{0\le s\le t} \mathcal{D}_\delta(s)+(\sup_{0\le s\le t}\mathcal{A}(s))
(\sup_{0\le s\le t}\mathcal{B}(s))\Big] \Phi(t),
\end{eqnarray}
and
\begin{eqnarray}\label{phi3}
\big\|u(t)\big\|_{L^{2}(\Omega)}
\le C\Big(\sup_{0\le s\le t} \mathcal{D}_\delta(s)\Big) \Phi(t).
\end{eqnarray}
Adding \eqref{phi1}, \eqref{phi2}, and \eqref{phi3}, we obtain that for $0<t\le t_0$, 
$$\Phi(t)\le C\Big[\sup_{0\le s\le t} \mathcal{D}_\delta(s)+(\sup_{0\le s\le t}\mathcal{A}(s))
(\sup_{0\le s\le t}\mathcal{B}(s))\Big] \Phi(t)\le \frac12\Phi(t),$$
provided $\epsilon_0>0$ is chosen so small that
$$C\Big[\sup_{0\le s\le t} \mathcal{D}_\delta(s)+(\sup_{0\le s\le t}\mathcal{A}(s))
(\sup_{0\le s\le t}\mathcal{B}(s))\Big]\le C[\epsilon_0^{1-\delta}+\epsilon_0]\le \frac12.$$
Therefore we must have $\Phi(t)\equiv 0$ for $0<t\le t_0$. Hence
$(u_1, n_1, m_1)\equiv (u_2, n_2, m_2)$ in $\Omega\times [0,t_0]$. We can repeat the same argument to
show that $(u_1, n_1, m_1)\equiv (u_2, n_2, m_2)$ beyond $t_0$ and until $t=T$. 
Hence the proof is complete.
\qed

\bigskip
\noindent{\bf Acknowledgments.} The first author is partially supported by NSF of China (No. 11571117 and 11371152), by Guangdong Provincial Natural Science Foundation (No.2016A030313451) and by the Science and Technology Program of Guangzhou (No.201707010040)
The second author is supported by the Chinese Council of Scholarship. The third author is partially supported by NSF DMS 1764417. This work was initiated during the visits of first and second author at Purdue University, both of whom would express their gratitudes to Department of Mathematics for its warm hospitality.


\bigskip
\bigskip
\begin{thebibliography}{99}

\bibitem{Ball} J. M. Ball, {\it Mathematics and liquid crystals}. Mol. crystals and liquid crystals,
647 (2017) 1-27.


\bibitem{Ericksen} J. Ericksen, {\em Hydrostatic theory of liquid crystal.}
Arch. Ration. Mech. Anal. 9, 371-378 (1962).

\bibitem{Ericksen1} J. L. Ericksen. {\em Inequalities in liquid crystal theory.}
Physics of Fluids (1958-1988), 9(6):1205-1207, 1966.

\bibitem{Leslie} F. Leslie, {\em Some constitutive equations for liquid crystals.}
Arch. Ration. Mech. Anal. 28, 265-283.

\bibitem{GV1} E. Govers, G. Vertogen, {\it Elastic continuum theory of biaxial nematics}.
Phys. Rev. A, Vol. 30, No. 4 (1984).

\bibitem{GV2} E. Govers, G. Vertogen, {\it Erratum: Elastic continuum theory of biaxial nematics}
[Phys. Rev. A, Vol. 30, No. 4 (1984)], Phys. Rev. A, Vol. 31, No. 3 (1985).

\bibitem{GV3} E. Govers, G. Vertogen,  {\it FLUID DYNAMICS OF BIAXIAL NEMATICS}.
Physica 133A (1985) 337-344.


\bibitem{HKL} R. Hardt, D. Kinderleherer, F. H. Lin, {\it Existence and partial regularity of static liquid crystal configurations}. Comm. Math. Phys. 105 (1986), no. 4, 547-570.

\bibitem{Hong} M. C.  Hong, {\em Global existence of solutions of the simplified Ericksen-Leslie system in dimension two.} Calc. Var. Partial Differential Equations 40, no. 1-2 (2011), 15-36.

\bibitem{HLLW} T. Huang, F. H. Lin, C. Liu, C. Y. Wang, {\em Finite time singularity of the nematic liquid crystal flow in dimension three.} Arch. Ration. Mech. Anal. 221 (2016), no. 3, 1223-1254.

\bibitem{LLWWZ}Chen-Chih Lai, Fanghua Lin, Changyou Wang, Juncheng Wei, Yifu Zhou,
{\it Finite time blow-up for the nematic liquid crystal flow in dimension two.}  arXiv:1908.10955.


\bibitem{Luck1} S. Luckhaus, {\it Partial H\"older continuity for minima of certain energies among maps into a Riemannian manifold}. Indiana Univ. Math. J. 37 (1988), no. 2, 349-367.

\bibitem{Luck2} S. Luckhaus, {\it Convergence of minimizers for the p-Dirichlet integral.}
Math. Z. 213 (1993), no. 3, 449?456.


\bibitem{LLC} F.  Leslie, L. Laverty, T. Carlsson,
{\it Continuum theory for biaxial nematic liquid crystals.}
Quart. J. Mech. Appl. Math. 45 (1992), no. 4, 595-606.

\bibitem{HW} J. Hineman, C. Y. Wang, {\em Well-posedness of nematic liquid crystal flow in $L^3_{\rm{uloc}}(\mathbb R^3)$}. Arch. Ration. Mech. Anal. 210 (2013), no. 1, 177-218. 

\bibitem{HLW} J. R. Huang, F. H. Lin, C. Y. Wang, {\em Regularity and existence of global solutions to the Ericksen-Leslie system in $\mathbb R^2$}. Comm. Math. Phys. 331 (2014), no. 2, 805-850.

\bibitem{Lin} F. H. Lin, {\em Nonlinear theory of defects in nematic liquid crystals; phase transition and flow phenomena.} Comm. Pure Appl. Math. 42 (1989), no. 6, 789-814.

\bibitem{LLW10} F. H. Lin, J. Y. Lin, C. Y. Wang, {\em Liquid crystal flows in two dimensions.}
Arch. Ration. Mech. Anal. 197 (2010), no. 1, 297-336.

%\bibitem{LW} F. H. Lin, C. Y. Wang, {\em  }. 

\bibitem{LW1} F. H. Lin, C. Y. Wang, {\em Energy identity of harmonic map flows from surfaces at finite singular time.} Calc. Var. Partial Differential Equations 6 (1998), no. 4, 369-380.

\bibitem{LW2} F. H. Lin, C. Y. Wang, {\em On the uniqueness of heat flow of harmonic maps and hydrodynamic flow of nematic liquid crystals.} Chin. Ann. Math. Ser. B 31 (2010), no. 6, 921-938.

\bibitem{LW3} F. H. Lin, C. Y. Wang, {\em Global existence of weak solutions of the nematic liquid crystal flow in dimension three.} Comm. Pure Appl. Math. 69 (2016), no. 8, 1532-1571.

\bibitem{LW4} F. H. Lin, C. Y. Wang, {\em Recent developments of analysis for hydrodynamic flow of nematic liquid crystals.} Philos. Trans. R. Soc. Lond. Ser. A Math. Phys. Eng. Sci. 372 (2014), no. 2029, 20130361, 18 pp.

\bibitem{LWbook} F. H. Lin, C. Y. Wang, The analysis of harmonic maps and their heat flows. World Scientific Publishing Co. Pte. Ltd., Hackensack, NJ, 2008. xii+267.

\bibitem{LL1} F. H. Lin, C. Liu, {\em Nonparabolic dissipative systems modeling the flow of liquid crystals.}
 CPAM XLVIII, 501-537 (1995).
 
\bibitem{LL2} F. H. Lin, C. Liu, {\em Partial regularity of the dynamic system modeling the flow of liquid crystals.}
DCDS 2(1), 1-22 (1998).

\bibitem{LWZZ} Q. Liu, C. Y. Wang, X. T. Zhang, J. F. Zhou, {\em On optimal boundary control of Ericksen-Leslie system in dimension two.} Calc. Var. Partial Differential Equations 59 (2020), no. 1, Paper No. 38, 64 pp.


\bibitem{Lemaire} L. Lemaire, Lemaire, {\em Applications harmoniques de surfaces riemanniennes.}
J. Differ. Geom. 13(1), 51-78 (1978).

\bibitem{Morrey} C. B. Morrey, Multiple integrals in the calculus of variations. Die Grundlehren der mathematischen Wissenschaften, Band 130 Springer-Verlag New York, Inc., New York 1966.

\bibitem{Moser} R. Moser, {\em An $L^p$ regularity theory for harmonic maps.}
Trans. Amer. Math. Soc. 367 (2015), no. 1, 1-30.

\bibitem{Qing} J. Qing, {\em On singularities of the heat flow for harmonic maps from surfaces into spheres}. 
Comm. Anal. Geom. 3(1-2), 297-315 (1995).

\bibitem{SSS} G. A. Seregin, T. N. Shilkin, V. A. Solonnikov, {\em Boundary partial regularity for the Navier-Stokes equations}. J. Math. Sci. 132(3), 339-358 (2006). 

\bibitem{Struew} M. Struwe, {\em On the evolution of harmonic mappings of Riemannian surfaces}. 
Comment. Math. Helvetici 60, 558-581 (1985).

\bibitem{SU1} R. Schoen, K. Uhlenbeck, {\em A regularity theory for harmonic mappings}. J. Diff. Geom.
17 (1982), 307-335. 

\bibitem{SP} B. Sharp, P. Topping, {\em Decay estimates for Rivi\'ere's equation, with applications to regularity and compactness}. Trans. Amer. Math. Soc. 365 (2013), no. 5, 2317-2339.

\bibitem{Wang} C. Y. Wang, {\em A remark on harmonic map flows from surfaces.}
Differential Integral Equations 12 (1999), no. 2, 161-166.

\bibitem{XZ} X. Xu, Z. F. Zhang, {\em Global regularity and uniqueness of weak solution for the 2-D liquid crystal flows.} J. Differential Equations 252 (2012), no. 2, 1169-1181.

\end{thebibliography}
\end{document}